\documentclass[review]{sn-jnl}
\normalbaroutside

\makeatletter
\def\ps@pprintTitle{%
 \let\@oddhead\@empty
 \let\@evenhead\@empty
 \def\@oddfoot{}%
 \let\@evenfoot\@oddfoot}
\makeatother

\makeatletter
\g@addto@macro\normalsize{%
  \setlength\abovedisplayskip{.6em}
  \setlength\belowdisplayskip{.6em}
  \setlength\abovedisplayshortskip{.6em}
  \setlength\belowdisplayshortskip{.6em}
}

\usepackage{amsmath,amscd,wrapfig,amsfonts,amssymb,mathtools,mathrsfs,breqn}
\usepackage{booktabs,multirow,lscape,cool,bm}
\usepackage{url,parskip,caption,bigints}
\newtheorem{remark}{Remark}

\theoremstyle{thmstyleone}%
\newtheorem{theorem}{Theorem}
\newtheorem{lemma}{Lemma}
\newtheorem{corollary}{Corollary}

\theoremstyle{thmstyletwo}

\theoremstyle{thmstylethree}%
\newtheorem{definition}{Definition}
\begin{document}

\title{On the adaptive Levin method}

\author*[1]{\fnm{Shukui} \sur{Chen}}
\author*[1,2]{\fnm{Kirill} \sur{Serkh}}
\author*[2]{\fnm{James} \sur{Bremer} \email{bremer@math.toronto.edu}}

\affil*[1]{\orgdiv{Department of Computer Science}, \orgname{University of Toronto}}

\affil[2]{\orgdiv{Department of Mathematics}, \orgname{University of Toronto}}

\abstract{The Levin method is a well-known technique for evaluating oscillatory integrals, which
operates by solving a certain ordinary differential equation in order to construct an
antiderivative of the integrand.   It was long believed that this approach suffers 
from ``low-frequency breakdown,'' meaning that the accuracy of the calculated value 
of the integral  deteriorates when the integrand is only slowly oscillating.
Recently presented experimental evidence, however, suggests that if a Chebyshev spectral method
is used to discretize the differential equation and the resulting linear
system is solved via a truncated singular value decomposition, then no 
low-frequency breakdown occurs.  Here, we provide a proof that 
this is the case,  and  our proof applies not only when the integrand is slowly oscillating,
but even in the case of stationary points.   
Our result puts adaptive schemes based on the Levin method
on a firm theoretical foundation and accounts for their behavior  in the 
presence of stationary points.  
%
%
We go on to point out that  by combining an adaptive Levin scheme with 
phase function methods for ordinary differential equations, a large class of oscillatory 
integrals involving special functions, including products of such functions and the compositions of such functions
with slowly-varying functions, can be 
easily evaluated without  the need for symbolic computations. Finally, we present the results of numerical experiments 
which illustrate the consequences of  our analysis and demonstrate the properties of 
the  adaptive Levin method.

}

\maketitle

\begin{section}{Introduction}

The Levin method, which was introduced in \cite{Levin}, is a classical technique  for evaluating integrals of the form
\begin{equation}
\int_a^b \exp(i g(x)) f(x)\,dx,
\label{introduction:int}
\end{equation}
where $f$ is a slowly varying, $g$ is real-valued and slowly varying and $g'$ is of large magnitude.
It operates by solving the ordinary differential equation
\begin{equation}
p'(x) + i g'(x) p(x) = f(x)
\label{introduction:ode}
\end{equation}
in order to find a function $p$ such that
\begin{equation}
\frac{d}{dx} \left(p(x) \exp(i g(x))\right) = f(x) \exp(i g(x)).
\end{equation}
The value of the integral (\ref{introduction:int}) is then equal to
\begin{equation}
p(b) \exp(ig(b)) - p(a) \exp(ig(a)).
\label{introduction:diff}
\end{equation}
Under the assumptions made above,  (\ref{introduction:ode}) admits a slowly-varying solution $p_0$,
which means that (\ref{introduction:int}) can be computed rapidly.
Moreover, although the differential operator
\begin{equation}
L\left[p\right](x) = p'(x) + i g'(x) p(x)
\label{introduction:diffop}
\end{equation}
appearing on the left-hand side of (\ref{introduction:ode}) has a one-dimensional nullspace
consisting of all multiples of the function $\eta(x) = \exp(-i g(x))$, 
the equation (\ref{introduction:ode}) can still be solved accurately
via a spectral collocation method if some care is taken.  
In particular,  as long as the chosen grid of collocation nodes is dense enough
to resolve $f$ and $g'$,  but not sufficiently dense to resolve $\eta$, 
the matrix discretizing the differential operator $L$ will be well-conditioned
and inverting it will result in a high-accuracy approximation
of the slowly varying solution $p_0$.

When the magnitude of $g'$ is sufficiently  small,
any grid of collocation nodes dense enough to resolve $f$ and $g'$ 
will also  resolve $\eta$.  As a consequence, the corresponding 
spectral discretization of $L$ will have a small singular value and numerical difficulties
can arise in the course of solving (\ref{introduction:ode}).
This phenomenon is known as  ``low-frequency breakdown,''   and it was long believed
to be a significant limitation of Levin methods.
The recent articles  \cite{LevinLi,LiImproved} present experimental evidence showing that,
when Chebyshev spectral methods are used to discretize the differential
equation (\ref{introduction:ode}) and the resulting linear system is solved
via a truncated singular value decomposition, no low-frequency breakdown seems to occur.

Here, we prove that this is the case.   We first show that, when $g'$ is nonvanishing, the Levin equation
admits a well-behaved solution that can be approximated by 
a polynomial expansion at a cost which depends on a measure of the complexity of $f$, $g'$  and 
the inverse function $g^{-1}$ of $g$, as well as the ratio of the maximum and minimum values of $g'$,
 but not on the magnitude of $g'$.   This result improves on the those presented
in both original Levin paper  \cite{Levin} and its follow-up \cite{Levin3}
by handling the case in which $g'$ is of positive, but arbitrarily small, magnitude
and by explicitly showing that the cost of representing the well-behaved solution of the Levin equation
via a polynomial expansion is independent of the magnitude of $g'$
(note that  the assumptions made in both \cite{Levin} and \cite{Levin3} fail to hold when $g'$ 
is of sufficiently small positive magnitude).
We then consider the case in which $g'$ is of small magnitude, possibly with  zeros.
In this event, $g$ need not be invertible, but we show that the Levin equation nonetheless admits 
a well-behaved solution, and that this solution can be approximated by a polynomial expansion
at a cost which decreases with the magnitude of $g'$.
We use these two results to prove that,  when the Levin equation is discretized
via a Chebyshev spectral collocation method   and the resulting linear system is solved
using a truncated singular value decomposition, high accuracy is obtained
regardless of the magnitude of $g'$ and whether or not $g'$ has zeros.

One implication of the absence of low-frequency breakdown is that
the Levin method can be used as the basis of an adaptive integration scheme
--- since adaptive subdivision reduces the effective magnitude of $g'$,
such an approach requires that the method remain accurate regardless
of the magnitude of $g'$.           
An adaptive Levin scheme was proposed by Moylan \cite{Moylan} 
before the absence of low-frequency breakdown was observed  in \cite{LevinLi,LiImproved}.
Moylan's scheme is based on  a more general version of the Levin method introduced 
in \cite{Levin2}. In particular,  the integrand is of the form
\begin{equation}
\int_a^b  \mathbf{w}(x) \cdot \mathbf{f}(x)\, dx
\label{introduction:int2}
\end{equation}
with  $\mathbf{f}\colon[a,b]\to\mathbb{C}^n$ slowly-varying and
$\mathbf{w}\colon[a,b]\to\mathbb{C}^n$ a solution of a system of
ordinary differential equations
\begin{equation}
\mathbf{w}'(x) = \mathbf{A}(x)\mathbf{w}(x)
\label{introduction:system}
\end{equation}
whose coefficient matrix $\mathbf{A}\colon[a,b]\to\mathbb{C}^{n\times n}$ is also slowly-varying.
The scheme of \cite{Levin2} operates by  finding a slowly-varying vector field
$\mathbf{F}\colon [a,b]\to\mathbb{C}^{n}$  satisfying  the system of differential equations
\begin{equation}
\mathbf{F}'(x) + \mathbf{A}(x)^t \mathbf{F}(x) = \mathbf{f}(x)
\label{introduction:systemode}
\end{equation}
and evaluating  (\ref{introduction:int2}) via the formula
\begin{equation}
\mathbf{F}(b) \cdot \mathbf{w}(b) -
\mathbf{F}(a) \cdot \mathbf{w}(a).
\end{equation}
Moylan's algorithm adaptively decomposes the domain of an oscillatory integral and uses the above formulation of the
Levin method to evaluate it over each subinterval.  An implementation of this scheme using Wolfram's Mathematica 
package is described in  \cite{Moylan}, and experimental evidence showing  that it works well in practice in the low-frequency regime, 
and even in the case of stationary points, is presented.  However, all of the estimates discussed in \cite{Moylan} break down in the 
low-frequency regime.  Moreover, the implementation of the adaptive Levin method described in
\cite{Moylan} relies heavily on symbolic and arbitrary precision computations, including in the solution of the
linear system which arises when (\ref{introduction:systemode}) is discretized.  Since that system
typically becomes highly ill-conditioned at low frequencies, it is unclear how the scheme would fare if the  calculations 
were performed in double precision arithmetic.

One ostensible advantage of the Moylan scheme over an adaptive scheme based on the original Levin formulation
\cite{Levin} is the great generality of the integrands which can be expressed  as 
solutions of systems of the form  (\ref{introduction:system}) with slowly-varying coefficient matrices.
However,  most oscillators of interest are solutions of scalar differential  equations of the form 
\begin{equation}
y^{(n)}(t) + q_{n-1}(t) y^{(n-1)}(t) + \cdots + q_1(t) y'(t) + q_0(t) y(t) = 0
\label{introduction:scalarode}
\end{equation}
whose coefficients are slowly-varying, and it is well known that  such equations admit slowly-varying 
phase functions.  More explicitly, under mild conditions on the
coefficients $q_0,\ldots,q_{n-1}$, there exist slowly-varying
 $\psi_1,\ldots,\psi_n$ such that 
\begin{equation}
\left\{
\exp\left(\psi_j(t)\right) : j=1,\ldots,n
\right\}
\label{introduction:phasebasis}
\end{equation}
is a basis in the space of solutions of (\ref{introduction:scalarode}).
This observation underlies the WKB method and almost all other
techniques  for the asymptotic approximation of solutions of ordinary differential equations
in the high-frequency regime (see, for instance, \cite{Miller}, \cite{Wasov} and \cite{SpiglerPhase1,SpiglerPhase2,SpiglerZeros}).  
See also  \cite{BremerRokhlin},   which contains a careful proof of the existence of slowly-varying phase functions for second order differential
equations that can be extended to the case of higher order scalar equations without too much difficulty.
A numerical algorithm for constructing slowly-varying functions for a scalar equations of the form
(\ref{introduction:scalarode}) in time independent of frequency is described in \cite{Scalar1}.
Specialized methods for the case of second order linear ordinary differential equations
are described in \cite{BremerPhase} and \cite{BremerPhase2}.  Because many of the oscillators
of most interest satisfy second order linear ordinary differential equations (including the Bessel functions,
Jacobi polynomials, spheroidal wave functions, and Hermite polynomials), the experiments
of this article focus on this case and we use the algorithm of \cite{BremerPhase2} in order
to construct the necessary slowly-varying phase functions.

It follows from the existence of slowly-varying phase functions for scalar
differential equations  that an adaptive method based on Levin's original formulation suffices to evaluate
a huge class of  oscillatory integrals.  Moreover,  when explicit exponential representations of oscillators are available, 
integrals involving the product of two oscillators or the composition of an oscillator with 
a slowly-varying function can be evaluated in a straightforward fashion.
This is in contrast to Moylan's adaptive scheme, which requires knowledge
of a system of differential equations satisfied by whatever combination of oscillators
is being considered.  Moylan's scheme relies heavily on Mathematica's symbolic capabilities to derive the
 necessary  differential equations.

The remainder of this article is structured as follows.
Section~\ref{section:preliminaries} discusses the requisite mathematical and
numerical preliminaries.    Section~\ref{section:leveq} contains our analysis showing
that the Levin equation always admits a well-behaved solution.
In Section~\ref{section:levnum}, we use the results of Section~\ref{section:leveq}
to establish that no accuracy is lost when a Chebyshev collocation scheme
is used to discretize the Levin equation and the resulting system of linear equations
is solved via a truncated singular value decomposition.
In Section~\ref{section:algorithm}, we give a detailed description
of the  adaptive Levin method.  
In Section~\ref{section:experiments}, we present the
results of numerical experiments conducted to demonstrate the properties of
the adaptive Levin method.  Finally, we close with a few brief remarks in
Section~\ref{section:conclusions}.

\end{section}

\begin{section}{Preliminaries}
\label{section:preliminaries}

\begin{subsection}{Notation and conventions}

We use capital script letters such as $\mathscr{A}$ for matrices,
and vectors are indicated with an overline as in $\overline{x}$.
We denote by  $L^\infty\left(\Omega\right)$ the Banach space of functions which
are essentially bounded on the Lebesgue measurable subset  $\Omega$ of $\mathbb{R}$,
and by $C^n\left([a,b]\right)$ the Banach space of 
functions $(a,b) \mapsto \mathbb{C}$ whose derivatives through
order $n$ are uniformly continuous and hence admit continuous extensions to $[a,b]$.  
We denote by $M\left([a,b]\right)$ the Banach space of complex Radon measures on 
$[a,b]$, which is the dual of $C^0\left([a,b]\right)$,
and  by $\left|\mu\right|$ the  total variation of $\mu \in M\left(\left[a,b\right]\right)$.
The notation $\|x\|_X$ is used for the norm of an element $x$ of one of the Banach spaces $X$ 
mentioned above.  
The Euclidean norm of a vector $\overline{x} \in \mathbb{C}^n$ is $\left\|\overline{x}\right\|$,
and $\left\|\mathscr{A}\right\|$ is the 
$\mathbb{C}^n \to \mathbb{C}^n$ Euclidean operator norm of the $n\times n$ matrix $\mathscr{A}$.

We denote by  $C^\infty\left([a,b]\right)$ 
the set of functions which are infinitely differentiable and whose derivatives
of all orders are uniformly continuous on $(a,b)$,
and use  $S\left(\mathbb{R}\right)$
for the  Schwartz space of infinitely differentiable functions whose derivatives of all orders are rapidly decaying
functions.  The space of tempered distributions is  $S'\left(\mathbb{R}\right)$.
We often write $\left\langle\varphi,f\right\rangle$
for the action of the tempered distribution $\varphi$ on the function $f \in S\left(\mathbb{R}\right)$,
although we occasionally use other notations for this when their meanings are clear.
The support of $\varphi \in S'\left(\mathbb{R}\right)$ is the complement of the union of all
open sets $U$ such that 
$\left\langle  \varphi,f \right\rangle=0$
whenever $f \in S\left(\mathbb{R}\right)$ has support contained in $U$.
The order of a tempered distribution $\varphi$ is the  least nonnegative integer $n$
with the property that for each compact set $K \subset \mathbb{R}$ there exists a constant $M_K$ 
such that 
\begin{equation}
\left|\left\langle \varphi, f\right\rangle\right|
\leq 
M_K \sup_{0 \leq k \leq n}\, \sup_{x \in K}  \left| D^k f(x) \right|
\end{equation}
for every function $f \in S\left(\mathbb{R}\right)$ whose support is contained in $K$
(we note that every tempered distribution has  finite order).
The space of tempered distributions of order $n$ which are supported on $[a,b]$ can, of course, be identified with 
the dual space of $C^n\left([a,b]\right)$.

We use the convention
\begin{equation}
\widehat{f}(\xi) = \frac{1}{2\pi}\int_{-\infty}^\infty f(x) \exp(-i \xi x)\, dx
\end{equation}
for the Fourier transform of $f \in S\left(\mathbb{R}\right)$ (this
is slightly nonstandard, but quite convenient when discussing the Levin equation).
Of course, $\widehat{f} \in S\left(\mathbb{R}\right)$ and we have
\begin{equation}
f(x) = \int_{-\infty}^\infty \widehat{f}(\xi) \exp(i\xi x)\,d\xi.
\end{equation}
%
We extend the Fourier transform to a mapping  $S'(\mathbb{R}) \to S'\left(\mathbb{R}\right)$ in the usual way,
via the formula
$\left\langle \widehat{\varphi}, g \right\rangle =  \left\langle \varphi, \widehat{g} \right\rangle$.
If $\varphi$ is compactly supported on $[a,b]$, then its Fourier transform
coincides with  the function $\widehat{\varphi}(\xi) = \left\langle \varphi, \eta_\xi \right\rangle$,
where $\eta_\xi(x) = 1/(2\pi) \exp(-ix\xi)$.  The Schwartz-Paley-Wiener theorem 
asserts that $\widehat{\varphi}$ is entire in this case, and gives bounds on its growth
at infinity.  Likewise, if $\widehat{\varphi}$ is compactly supported on $[a,b]$, then
$\varphi$ is the entire function defined via the formula  $\varphi(x) = \left\langle \widehat{\varphi}, \nu_x \right\rangle$
with  $\nu_x(\xi) = \exp(i \xi x)$.    If the Fourier transform of $\varphi$ is supported on the 
finite interval $[-c,c]$, then we say that $\varphi$ has  bandlimit $c$.  
Note that we do not require that $c$ be the smallest positive real number with this property.


The notation $x \lesssim y$ indicates that there is some constant $C$ not depending on $y$ such that
$x \leq C y$.  
We say that
%
\begin{equation}
f(z) = \mathcal{O}\left(g(z)\right)\ \ \mbox{as} \ \ z \to \infty
\end{equation}
provided  there exist $M > 0$ and $C >0$ such that 
\begin{equation}
\left|f(z)\right| \leq C \left|g(z)\right| \ \ \mbox{for all}  \ \ z > M.
\end{equation}

\end{subsection}

\begin{subsection}{Chebyshev interpolation}
\label{section:preliminaries:chebyshev}

We use $T_n$ to denote the Chebyshev polynomial of degree $n$ and 
\begin{equation}
-1 = x^{\mbox{\tiny cheb}}_{1,k} <  x^{\mbox{\tiny cheb}}_{2,k} <  \ldots <  x^{\mbox{\tiny cheb}}_{k,k}=1
\label{preliminaries:chebnodes}
\end{equation}
for the $k$-point grid of Chebyshev extremal nodes on the interval $[-1,1]$.  That is,
\begin{equation}
x^{\mbox{\tiny cheb}}_{j,k} = \cos\left(\pi \frac{k-j}{k-1}\right), \ \ \ \ j=1,\ldots,k.
\end{equation}
It is well known that 
\begin{equation}
\begin{aligned}
&\frac{2}{k-1}\sum_{j=1}^{k}\mbox{}^{''} T_n\left(x^{\mbox{\tiny cheb}}_{j,k}\right) T_m\left(x^{\mbox{\tiny cheb}}_{j,k}\right)\\
=
&\begin{cases}
1 & \mbox{if}\  n=1,\ldots,k-2 \ \mbox{and}\  \mbox{mod}(m, 2(k-1)) = n\\
2 & \mbox{if}\  n=0, (k-1) \ \mbox{and}\  \mbox{mod}(m, 2(k-1)) = n\\
0 & \mbox{otherwise},
\end{cases}
\end{aligned}
\label{preliminaries:cheborth}
\end{equation}
where the double prime symbol after the summation indicates that the first and last term of the series
are weighted by $1/2$ (this formula can be found, for instance, in Chapter~4 of \cite{Mason}).
If $f \in C^\infty\left([-1,1]\right)$, then we use
$\mathsf{P}_k\left[f\right]$ to denote the Chebyshev expansion of the form
\begin{equation}
\sum_{n=0}^{k-1} a_n T_n(x)
\label{preliminaries:chebexp1}
\end{equation}
which agrees with $f$ at the nodes (\ref{preliminaries:chebnodes}).  It follows from
(\ref{preliminaries:cheborth}) that the coefficients in (\ref{preliminaries:chebexp1}) are
given by 
\begin{equation}
a_n =
\begin{cases}
\begin{aligned}  \frac{1}{k-1}\sum_{j=1}^{k}\mbox{}^{''} T_n\left(x_j^{\mbox{\tiny cheb}}\right) f\left(x_j^{\mbox{\tiny cheb}}\right)\end{aligned}  
& \mbox{if}\ n = 0\ \mbox{or}\ n=k-1     \\[1.3em]
\begin{aligned}   \frac{2}{k-1}\sum_{j=1}^{k}\mbox{}^{''} T_n\left(x_j^{\mbox{\tiny cheb}}\right) f\left(x_j^{\mbox{\tiny cheb}}\right)\end{aligned}  & \mbox{if}\ n =1,\ldots,k-2.
\end{cases}
\end{equation}
The $k^{th}$ order Chebyshev spectral differential matrix is the $k\times k$ matrix $\mathscr{D}_k$ 
which maps the vector
\begin{equation}
\left(
\begin{array}{c}
g\left(x_{1,k}^{\mbox{\tiny cheb}}\right)\\
g\left(x_{2,k}^{\mbox{\tiny cheb}}\right)\\
\vdots\\
g\left(x_{k,k}^{\mbox{\tiny cheb}}\right)\\
\end{array}
\right)
\end{equation}
of values of an expansion of the form
\begin{equation}
g(x) = \sum_{n=0}^{k-1} a_n T_n(x)
\label{preliminaries:chebexp0}
\end{equation}
at the extremal Chebyshev nodes to the vector
\begin{equation}
\left(
\begin{array}{c}
g'\left(x_{1,k}^{\mbox{\tiny cheb}}\right)\\
g'\left(x_{2,k}^{\mbox{\tiny cheb}}\right)\\
\vdots\\
g'\left(x_{k,k}^{\mbox{\tiny cheb}}\right)\\
\end{array}
\right)
\end{equation}
of values of the derivative of (\ref{preliminaries:chebexp0}) at the extremal Chebyshev nodes.

If  $f \in C^\infty\left([-1,1]\right)$ is infinitely differentiable, then it admits a uniformly convergent Chebyshev series
\begin{equation}
f(x) = \sum_{n=0}^\infty  b_n T_n(x)
\label{preliminaries:chebexp2}
\end{equation}
such that
\begin{equation}
\left|b_n\right| = \mathcal{O}\left(\frac{1}{n^m}\right)
\ \ \mbox{as}\ \ n\to\infty
\end{equation}
for all $m \geq 1$.
From (\ref{preliminaries:cheborth}) it is clear that  the coefficients $\{b_n\}$ in (\ref{preliminaries:chebexp2})
are related to the coefficients $\{a_n\}$ in the expansion of $\mathsf{P}_k\left[f\right]$ 
through the formula
\begin{equation}
a_n = b_n + b_{n+2(k-1)} +  b_{n+4(k-1)} + \ldots\ .
\label{preliminaries:chebrel}
\end{equation}
Since $\left|T_n(x)\right| \leq 1$ and  $\left|T_n'(x)\right| \leq n^2$ for all $x \in[-1,1]$,
it follows from (\ref{preliminaries:chebrel}) that 
\begin{equation}
\left\| \mathsf{P}_k\left[f\right] - f \right\|_{L^\infty\left([-1,1]\right)}
\leq 2 \sum_{n=k}^\infty \left|b_n\right|
\label{preliminaries:chebound1}
\end{equation}
and
\begin{equation}
\left\| \mathsf{P}_k\left[f\right]' - f' \right\|_{L^\infty\left([-1,1]\right)}
\leq 2 \sum_{n=k}^\infty n^2 \left|b_n\right|.
\label{preliminaries:chebound2}
\end{equation}
Obviously, if $f$ is a polynomial of degree less than $k$, then 
\begin{equation}
\mathsf{P}_k\left[f\right](x) = f(x)
\end{equation}
and
\begin{equation}
\mathsf{P}_k\left[f'\right](x) = f'(x) = \mathsf{P}_k\left[f\right]'(x).
\end{equation}
%


%

%
\end{subsection}

\begin{subsection}{Approximation by bandlimited functions}
\label{sectioon:preliminaries:approx}
 
It will often be necessary for us to approximate a function $f$ given on the interval $[-1,1]$
via a ``well-behaved'' bandlimited function $f_b$.
For us, well-behaved means
that the $L^\infty\left(\mathbb{R}\right)$ norms of 
$f_b$, its Fourier transform and the derivative of its Fourier transform
are bounded by small constant multiples of  the $L^\infty\left([-1,1]\right)$ norm of $f$.
Moreover, it is  desirable to minimize the bandlimit of $f_b$ subject to these constraints.

Under various regularity assumptions on $f$, it is possible to prove the existence
of an approximate satisfying these requirements.   The following theorem,  which  is a slightly modified 
version of Theorem~1 in  \cite{Boyd}, is an example of a result of this type which
holds under relatively weak regularity conditions on $f$.

\begin{theorem}
\label{preliminaries:approx:thm1}
Suppose that $f:[-1,1] \to \mathbb{C}$ admits an infinitely differentiable
extension to an open neighborhood of $[-1,1]$.
Then for each positive integer $m$ and each real number $c>1$,
there exist a constant $k$ depending on $m$ but not $c$,
and  a function $f_b \in S\left(\mathbb{R}\right)$ 
such that 
\begin{enumerate}
\item
$\begin{aligned}
\widehat{f_b}\end{aligned}$ 
is supported on $[-c-2,c+2]$,

\item
$\begin{aligned}
\left\| f_b - f \right\|_{L^\infty\left(\left[-1,1\right]\right)} \leq \frac{k}{c^m},
\label{approx:conclusion1}
\end{aligned}$

\item
$\begin{aligned}
\left\|f_b(x)\right\|_{L^\infty\left(\mathbb{R}\right)} \leq 2\|f\|_{L^\infty\left([-1,1]\right)} + \frac{k}{c^m},
\end{aligned}$

\item
$\begin{aligned}
\left\|\widehat{f_b}(\xi)\right\|_{L^\infty\left(\mathbb{R}\right)}\leq 4 \|f\|_{L^\infty\left([-1,1]\right)}
\end{aligned}$
and

\item
$\begin{aligned}
\left\|\widehat{f_b}'(\xi)\right\|_{L^\infty\left(\mathbb{R}\right)}\leq 4 \|f\|_{L^\infty\left([-1,1]\right)}.
\end{aligned}$

\end{enumerate}
\end{theorem}

\begin{proof}
We let $M = \left\|f\right\|_{L^\infty\left([-1,1]\right)}$ and 
choose $0 < \delta<1$ such that $f$ is infinitely
differentiable on  $[-1-\delta,1+\delta]$ and bounded in magnitude by $2M$ there.
We take $T_1$ to be an infinitely differentiable window function $\mathbb{R} \to \mathbb{R}$
such that   $\left|T_1(x)\right| \leq 1 $ for all $x\in \mathbb{R}$,
$T_1(x) = 1$ for all $x \in [-1,1]$ and  $T_1(x)  =0$ for all 
$\left|x\right| \geq 1+\delta$.
%
One can construct this  function quite easily
using the infinitely differentiable ramp function
\begin{equation}
H(x) = 
\begin{cases}
\frac{1}{2} \left( 1 + \mbox{erf}\left(\frac{x}{\sqrt{1-x^2}}\right) \right) & x \in [-1,1]\\
0 & |x| > 1
\end{cases}
\label{approx:h}
\end{equation}
suggested in \cite{Boyd}.  Then $f_1(x) = f(x) T_1(x)$ is an element of the space $S\left(\mathbb{R}\right)$,
and so is its Fourier transform $\widehat{f_1}$.
Consequently, there exists a constant $k_1$ (depending on $m$) such that
\begin{equation}
\sup_{|\xi|\geq 1}\left|\widehat{f_1}(\xi)\right| \leq \frac{k_1}{\left|\xi\right|^{m+1}}.
\label{approx:mdecay}
\end{equation}
Since   $T_1$ is bounded in magnitude by $1$,
\begin{equation}
\left|f_1(x)\right| \leq 2M\ \ \mbox{for all} \ \ x \in \mathbb{R}.
\label{approx:f1bound}
\end{equation}
It follows from this that 
\begin{equation}
\left|\widehat{f_1}(\xi) \right|
= 
\left|\frac{1}{2\pi}\int_{-1-\delta}^{1+\delta}
 f_1(x) \exp(-i\xi x)\, dx\right|
\leq 
 \frac{2}{\pi} (1+\delta)M
\leq \frac{4}{\pi} M
\label{approx:f1fbound}
\end{equation}
and
\begin{equation}
\left|\widehat{f_1}'(\xi) \right|
=
\left|
\frac{-i}{2\pi}\int_{-1-\delta}^{1+\delta}
f_1(x) x \exp(-i\xi x)\, dx
 \right|
\leq 
\frac{2}{\pi} \left(1+\delta\right)^2 M
\leq
\frac{8}{\pi} M
\label{approx:f1fpbound}
\end{equation}
for all $\xi\in\mathbb{R}$.

Now we let $T_2$ be an infinitely differentiable window function such that 
$\left|T_2(x)\right| \leq 1$ and $\left|T_2'(x)\right| \leq 1$ 
for all $x \in \mathbb{R}$, $T_2(x) = 1$ for all $x\in [-c,c]$
and  $T_2(x) = 0 $ for all $|x|>c+2$.  Such a function
can be easily constructed using the ramp function $H$ defined
in  (\ref{approx:h}). We define
$f_b$ by 
\begin{equation}
\widehat{f_b}(\xi) = \widehat{f_1}(\xi) T_2(\xi),
\end{equation}
so that the first of the properties listed above is clearly satisfied.
From (\ref{approx:mdecay}) and the definition of $T_2$, it is clear that
for all $x \in [-1,1]$,
\begin{equation}
\begin{aligned}
\left|f_1(x)-f_b(x) \right|
&=
\left|
\frac{1}{2\pi} \int_{-\infty}^\infty
\widehat{f_1}(\xi)\left(1 - T_2(\xi) \right) \exp(i\xi x)\,d\xi
\right| \\
&\leq
\frac{1}{\pi}
 \int_{|\xi|>c} \left|\widehat{f_1}(\xi)\right|\,d\xi  \\
&\leq
\frac{2 k_1}{m \pi\,  c^{m}}.
\end{aligned}
\end{equation}
This inequality establishes the second of the properties of the function $f_b$ listed
above and, by combining it with (\ref{approx:f1bound}), we obtain the third.
It follows from (\ref{approx:f1fbound}), (\ref{approx:f1fpbound}) and the properties of $T_2$ that
\begin{equation}
\begin{aligned}
 \left|\widehat{f_b}(\xi) \right|
= 
\left|\widehat{f_1}(\xi) T_2(\xi) \right|
\leq \frac{4}{\pi} M
\end{aligned}
\end{equation}
and 
\begin{equation}
\begin{aligned}
 \left|\widehat{f_b}'(\xi) \right|
= 
\left|\widehat{f_1}'(\xi) T_2(\xi) + \widehat{f_1}(\xi) T_2'(\xi) \right|
\leq \frac{4}{\pi} M + \frac{8}{\pi} M 
= \frac{12}{\pi} M
\end{aligned}
\end{equation}
for all $\xi \in \mathbb{R}$.  This establishes the last two of the properties
of $f_b$ listed above.
\end{proof}

More precise results of this type can be given under stronger regularity
assumptions on $f$.   However, the necessary arguments are rather technical and beyond the scope
of this paper.
Accordingly,  now that we have established the existence of a suitable bandlimited
approximate  under the weak
condition that $f$ admits an infinitely differentiable extension in a neighborhood
of $[-1,1]$,  we prefer to simply introduce the following definition.

\begin{definition}
\label{preliminaries:approx:def}
Suppose that  $f:[-1,1] \to \mathbb{C}$ admits
an  infinitely differentiable extension to an open interval containing 
$[-1,1]$. Then, for each real number $0 <\epsilon <1$, we denote by $c_f(\epsilon)$ 
the smallest positive real number $c$ such that there exists a function $f_b \in S\left(\mathbb{R}\right)$
of bandlimit $c$ with the following properties:

\begin{enumerate}
\item
$\begin{aligned}
\left\|f - f_b \right\|_{L^\infty\left(\left[-1,1\right]\right)} 
\leq \epsilon \left\|f\right\|_{L^\infty\left([-1,1]\right)},
\end{aligned}$

\item
$\begin{aligned}
\left\|f_b\right\|_{L^\infty\left(\mathbb{R}\right)} \leq 4 \left\|f\right\|_{L^\infty\left([-1,1]\right)},
\end{aligned}$

\item
$\begin{aligned}
\left\|\widehat{f_b}\right\|_{L^\infty\left(\mathbb{R}\right)} \leq 4  \left\|f\right\|_{L^\infty\left([-1,1]\right)}
\end{aligned}$ and

\item
$\begin{aligned}
\left\|\widehat{f_b}'\right\|_{L^\infty\left(\mathbb{R}\right)} \leq 4  \left\|f\right\|_{L^\infty\left([-1,1]\right)}.
\end{aligned}$

\end{enumerate}
\end{definition}

The following is an immediate consequence of  Theorem~\ref{preliminaries:approx:thm1}.
\begin{corollary}
If $f:[-1,1] \to \mathbb{C}$ admits an infinitely differentiable extension
to a neighborhood of $[-1,1]$, then for every positive integer $m$,
\begin{equation}
c_f(\epsilon) = \mathcal{O}\left(\left(\frac{1}{\epsilon}\right)^{\frac{1}{m}}\right)
\ \ \mbox{as} \ \ \epsilon\to 0.
\end{equation}
\end{corollary}

Under the assumption that $f:[-1,1] \to \mathbb{C}$ admits an extension which is bounded
and analytic on a horizontal strip containing the real axis, then 
we clearly have
\begin{equation}
c_f(\epsilon) = \mathcal{O}\left(\log\left(\frac{1}{\epsilon}\right)\right)
\ \ \mbox{as} \ \ \epsilon\to 0.
\end{equation}
We conjecture that this also holds when
$f$ admits an extension which is analytic in an open neighborhood of $[-1,1]$,
but a proof of this is beyond the scope of the present article.

\end{subsection}

\begin{subsection}{Legendre expansions of bandlimited functions}

It is a consequence of the Schwartz-Paley-Wiener theorem that  
if $\widehat{\varphi}$ has compact support, then $\varphi$
is an entire function satisfying certain growth conditions at $\infty$.  
It follows from this and standard results in approximation
theory (which can be found in \cite{Davis}, for example) that the Legendre
expansion of $\varphi$ decays superexponentially.  Here, 
for the sake of concreteness, we give a simple
bound on the  coefficients in the Legendre expansion of $\varphi$ which depends on the order $n$ 
of $\widehat{\varphi}$ and its bandlimit $c$.    Throughout this subsection, we use $P_k$ to denote
the Legendre polynomial of degree $k$ and $j_\nu$ to denote the spherical Bessel function
of order $\nu$.  Definitions of these functions can be found, for instance, in \cite{DLMF}.

\begin{lemma}
 \label{preliminaries:measure:lemma1}
If $\varphi$  is a tempered distribution of order $n$ which is supported on $[a,b]$,
then there exist complex numbers $c_0,\ldots,c_{n-1}$ and a complex Radon measure $\mu$ on
$[a,b]$ such that 
\begin{equation}
\left\langle\varphi,f\right\rangle = \sum_{k=0}^{n-1} c_k f^{(k)}(a) + \int_a^b f(x)\,d\mu(x).
\label{preliminaries:measure:form}
\end{equation}
for all $f \in C^n\left([a,b]\right)$.
\end{lemma}

\begin{proof}
Because the space of tempered distributions of order $n$ which are supported on $[a,b]$
can be identified with the dual space of $C^n\left([a,b]\right)$, it suffices to show that 
any element of the dual of $C^n\left([a,b]\right)$ is of the form (\ref{preliminaries:measure:form}).
The case $n=0$ follows by a trivial application of the Riesz representation theorem, 
so we suppose that $n \geq 1$.  
Because any function $f \in C^n\left([a,b]\right)$ can be written as 
\begin{equation}
f(x) = \sum_{k=0}^{n-1} \frac{f^{(k)}(a)}{k!}  (x-a)^k + 
\int_a^x  \frac{f^{(n)}(u)}{(n-1)!} (u-a)^{n-1} \,du,
\end{equation}
the map
\begin{equation}
f \mapsto \left(f(a),f'(a),\ldots,f^{(n-1)}(a),f^{(n)}\right)
\end{equation}
is an isomorphism from $C^n\left([a,b]\right)$ to $\mathbb{C}^n \times C\left([a,b]\right)$. 
The result now follows from this and the observations
that the dual of $C\left([a,b]\right)$ is the space of complex Radon measures $M([a,b])$
and $\mathbb{C}^n$ is its own dual space.
\end{proof}

\begin{lemma}
 \label{preliminaries:measure:lemma2}
If the Fourier transform of the tempered distribution 
$\varphi$ is of order $n$ and has support on the interval $[a,b]$, then
$\varphi$ coincides with an entire function of the form
\begin{align}
\varphi(x) = p_{n-1}(x) + x^n \int_a^b e^{i \xi x}\, d\mu(\xi),
\end{align}
where $p_{n-1}$ is a polynomial of degree at most $n-1$ and $\mu \in M\left(\left[a,b\right]\right)$.

\end{lemma}

\begin{proof}
The tempered distribution $\varphi$ is given by
\begin{equation}
\varphi(x) = \left\langle \widehat{\varphi}, \xi_x \right\rangle,
\end{equation}
where $\xi_x(\xi) = \exp(i x \xi)$.  The result follows from this and Lemma~\ref{preliminaries:measure:lemma1}.
\end{proof}


\begin{lemma}
\label{preliminaries:measure:lemma3}
For all real-valued $\xi$ and  nonnegative integers $k$,
\begin{equation}
\left|\int_{-1}^1 \exp(i \xi x) P_k(x)\ dx\right| \leq \frac{2 \left|\frac{\xi}{2} \right|^k}{\Gamma(k+1)}.
\end{equation}
\end{lemma}
\begin{proof}
The formula
\begin{equation}
\int_{-1}^1 \exp(i \xi x) P_k(x)\ dx = 2i^k j_k(\xi),
\end{equation}
%
can be found in Section~7.8  of \cite{HTFII}, among many other sources.
Combining it with the well-known inequality
\begin{equation}
\left|j_k(z)\right| \leq \frac{\left|\frac{z}{2}\right|^k}{\Gamma(k+1)},
\label{preliminaries:measure:crudejk}
\end{equation}
which is a special case of Formula~10.14.4 in \cite{DLMF},
yields the conclusion of the lemma.
\end{proof}

\begin{theorem}
\label{preliminaries:measure:theorem}
Suppose that the Fourier transform of $\varphi \in S'\left(\mathbb{R}\right)$ is a tempered distribution
of order $n$ supported on the interval $[-c,c]$ with $c\geq1$.
Then $\varphi$ is an entire function 
and the coefficients in the Legendre expansion
\begin{equation}
\varphi(x) = \sum_{m=0}^\infty a_m P_m(x)
\label{preliminaries:measure:lege}
\end{equation}
of $\varphi$ satisfy 
\begin{equation}
\left|a_m\right| \lesssim 
\frac{\left(\frac{c}{2}\right)^{m+n}}{\Gamma(m-n+1)}.
\label{preliminaries:measure:bound}
\end{equation}
for all $m \geq n$.
\end{theorem}
\begin{proof}
%
By Lemma~\ref{preliminaries:measure:lemma2},
\begin{equation}
\begin{aligned}
\int_{-1}^1 \varphi(x) P_m(x)\, dx  =&\int_{-1}^1 p_{n-1}(x) P_m(x)\, dx \\
+ &\int_{-c}^c \left(\int_{-1}^1 x^n  \exp(i\xi x) P_m(x)\,dx\right)\,d\mu(\xi),
\end{aligned}
\end{equation}
where $p_{n-1}$ is a polynomial of degree at most $n-1$ and $\mu \in M\left(\left[-c,c\right]\right)$.
For $m\ge n$, 
\begin{equation}
\int_{-1}^1 p_{n-1}(x) P_m(x)\, dx = 0,
\end{equation}
and so 
\begin{equation}
\left|\int_{-1}^1 \varphi(x) P_m(x)\, dx\right|  \le
\left|\mu\right|([-c,c])  \max_{ \xi\in [-c,c]}
\left|\int_{-1}^1 \exp(i \xi x)  x^n P_m(x)\,dx\right|.
\label{preliminaries:measure:ineq}
\end{equation}
We now observe that
\begin{equation}
x^n P_m(x) = \sum_{k=m-n}^{m+n} b_k P_k(x),
\label{preliminaries:measure:xnpm}
\end{equation}
where
\begin{equation}
b_k = \sqrt{k+\frac{1}{2}}\int_{-1}^1 x^n P_m(x) P_k(x)\,dx.
\end{equation}
By combining (\ref{preliminaries:measure:xnpm}) and Lemma~\ref{preliminaries:measure:lemma3}, we see
that
\begin{equation}
\begin{aligned}
\left|\int_{-1}^1 \exp(i\xi x)  x^n P_m(x)\,dx\right|
&\leq 
\sum_{k=m-n}^{m+n} \frac{2 \left|b_k\right| \left|\frac{\xi}{2}\right|^k}{\Gamma(k+1)}\\
&\leq
\frac{(4n+2)  \max\{\left|b_k\right|\}\,  \left(\frac{c}{2}\right)^{m+n}}{\Gamma(m-n+1)},
\end{aligned}
\label{eq:101}
\end{equation}
for all $|\xi| \leq c$.   
Together with (\ref{preliminaries:measure:ineq}), this gives
us (\ref{preliminaries:measure:bound}).
\end{proof}


\begin{remark}
\label{preliminaries:measure:remark1}
Using the well-known approximation
\begin{equation}
j_\nu(z) \approx
\frac{1}{z\sqrt{2e}} \left( \frac{ez}{2(\nu+\frac{1}{2})} \right)^{\nu+1}
=
\left( \frac{e}{2\sqrt{2e}(\nu+\frac{1}{2})} \right)
\left( \frac{ez}{2(\nu+\frac{1}{2})} \right)^{\nu}
\label{proof:jk}
\end{equation}
in lieu of the rather crude bound (\ref{preliminaries:measure:crudejk}), we see that
\begin{equation}
\begin{aligned}
&\left|\int_{-1}^1 \exp(i x \xi)  x^n P_m(x)\,dx\right|\\
=
&\left|
\sum_{k=m-n}^{m+n}
 2 i^k b_k j_k(\xi)
\right|
\\
&\leq
(4n+2) \max\{\left|b_k\right|\}
\sum_{k=m-n}^{m+n}\left|j_k(\xi)\right|\\
&\approx
(4n+2) \max\{\left|b_k\right|\}
\sum_{k=m-n}^{m+n}
\left( \frac{e}{2\sqrt{2e}(k+\frac{1}{2})} \right)
\left( \frac{ec}{2(k+\frac{1}{2})} \right)^{k}
\end{aligned}
\end{equation}
for all $|\xi| \leq c$.    When $m > ec/2  + n -\frac{1}{2}$, 
\begin{equation}
\left( \frac{ec}{2(k+\frac{1}{2})} \right) < 1
\end{equation}
and  the preceding sum is bounded by a multiple of 
\begin{equation}
\frac{2n}{m-n+\frac{1}{2}}
\left(\frac{e c}{2\left(m-n+\frac{1}{2}\right)}\right)^{m-n}.
\end{equation}
So we expect the coefficients in the Legendre expansion of $\varphi$ to decay superexponentially as soon
as this condition is met.
\end{remark}

\end{subsection}

\begin{subsection}{Truncated singular value decompositions}

If $\mathscr{A}$ is a complex-valued $n\times n$ matrix, then any decomposition
of the form
\begin{equation}
\mathscr{A} = \left(
\begin{array}{ccccc}
\overline{u_1} & \overline{u_2} & \cdots & \overline{u_n}
\end{array}
\right)
\left(
\begin{array}{ccccccc}
\sigma_1 & \\
         & \sigma_2 & \\
         &          & \ddots & \\
         &          &        &  \sigma_n
\end{array}
\right)
\left(
\begin{array}{ccccc}
\overline{v_1} & \overline{v_2} & \cdots & \overline{v_n} 
\end{array}
\right)^*,
\label{proof:svd}
\end{equation}
where $\sigma_1,\ldots,\sigma_n \in \mathbb{R}$ and both $\{\overline{u_1},\ldots,\overline{u_n}\}$ and $\{\overline{v_1},\ldots,\overline{v_n}\}$ 
are  orthonormal bases in $\mathbb{C}^n$,
is known as a singular value decomposition  of $\mathscr{A}$.  The quantities $\sigma_1,\ldots,\sigma_n$
are uniquely determined up to ordering, 
and they are known as the singular values of $\mathscr{A}$.  It is 
conventional to arrange them in descending order, and we will assume that this 
is the case with all singular value decompositions that we consider.

A truncated singular value decomposition  of $\mathscr{A}$ is any approximation
of the form
\begin{equation}
\mathscr{A} \approx 
\left(
\begin{array}{ccccc}
\overline{u_1} & \overline{u_2} & \cdots & \overline{u_k}
\end{array}
\right)
\left(
\begin{array}{ccccccc}
\sigma_1 & \\
         & \sigma_2 & \\
         &          & \ddots & \\
         &          &        &  \sigma_k \\ 
\end{array}
\right)
\left(
\begin{array}{ccccc}
\overline{v_1} & \overline{v_2} & \cdots & \overline{v_k}
\end{array}
\right)^*,
\label{proof:tsvd}
\end{equation}
where (\ref{proof:svd}) is a singular value decomposition of $A$ and $1 \leq k \leq n$.
Typically, some desired precision $\epsilon > 0$ is specified and $k$ is taken to be the
least integer between $1$ and $n-1$ such that
$\sigma_{k+1} <  \epsilon$, if such an integer exists, or $k=n$ otherwise.
We say that (\ref{proof:tsvd}) is a singular value decomposition
which has been truncated at precision $\epsilon$ when $k$
has been chosen in this fashion.  We call the vector
\begin{equation}
\overline{x}=
\left(
\begin{array}{ccccc}
\overline{v_1} & \overline{v_2} & \cdots & \overline{v_k}
\end{array}
\right)
\left(
\begin{array}{ccccccc}
\frac{1}{\sigma_1} & \\
         & \frac{1}{\sigma_2} & \\
         &          & \ddots & \\
         &          &        &  \frac{1}{\sigma_k} \\ 
\end{array}
\right)
\left(
\begin{array}{ccccc}
\overline{u_1} & \overline{u_2} & \cdots & \overline{u_k}
\end{array}
\right)^* y
\end{equation}
the solution of the linear system $Ax=y$ obtained from the 
truncated singular value decomposition (\ref{proof:tsvd}).

The following theorem implies that, when a linear system admits an approximate
solution with a modest norm, and it is solved numerically using a truncated
singular value decomposition, the computed solution will have both
a small residual and a modest norm. The proof can be found in~
\cite{zhao2023approximation}.

\begin{theorem}
\label{thm:tsvd}
Suppose that $\mathscr{A}\in \mathbb{R}^{m \times n}$, where $m\ge n$, and let
$\sigma_1 \ge \sigma_2 \ge \cdots \ge \sigma_n$ be the singular values 
of $\mathscr{A}$. Suppose that $\overline{x}\in \mathbb{R}^n$ satisfies
  \begin{align}
\mathscr{A}\overline{x} = \overline{b}.
  \end{align}
Let $\epsilon > 0$, and suppose that 
  \begin{align}
\hat{\overline{x}}_k = (\mathscr{A}+\mathscr{E})^\dagger_k (\overline{b}
  +\overline{e}),
  \end{align}
where $(\mathscr{A}+\mathscr{E})^\dagger_k$ is the pseudo-inverse of the
$k$-TSVD of $\mathscr{A}+\mathscr{E}$, so that 
  \begin{align}
\hat \sigma_{k} \ge \epsilon \ge \hat \sigma_{k+1},
  \end{align}
where $\hat \sigma_k$ and $\hat \sigma_{k+1}$ are the $k$th and $(k+1)$th
largest singular values of $\mathscr{A}+\mathscr{E}$, and where
$\mathscr{E}\in \mathbb{R}^{m\times n}$ and $\overline{e}\in \mathbb{R}^m$, with
$\|\mathscr{E}\|_2 < \epsilon/2$. Then
  \begin{align}
\|\hat{\overline{x}}_k\|_2 \le \frac{1}{\hat \sigma_k}
(2\epsilon\|\overline{x}|_2 +
\|\overline{e}\|_2) + \|\overline{x}\|_2.
  \end{align}
and
  \begin{align}
\|\mathscr{A} \hat{\overline{x}}_k - \overline{b}\|_2 \le 5\epsilon
\|\overline{x}\|_2 + \frac{3}{2}
\|\overline{e}\|_2.
  \end{align}
\end{theorem}

We will make use of the following simplified version of this theorem.

\begin{corollary}
\label{preliminaries:numerical:svd}
Suppose that $\epsilon > 0$,  $\mathscr{A}$ is an $n\times n$ matrix with complex entries, and that
\begin{equation}
\mathscr{A}\overline{x} = \overline{y} + \overline{\delta y}
\end{equation}
for some $\overline{x}$, $\overline{y}$ and $\overline{\delta y}$ in $\mathbb{C}^n$ with
\begin{equation}
\left\| \overline{\delta y} \right\|  \lesssim \epsilon \left\|\mathscr{A}\right\| \left\|\overline{x}\right\|.
\end{equation}
Suppose further that the linear system
\begin{equation}
\mathscr{A}\overline{x} = \overline{y}
\end{equation}
is solved in finite precision arithmetic using a singular value decomposition which is truncated at precision 
$\epsilon \left\|\mathscr{A}\right\|$, and that $\overline{z}$ is the computed solution.
Then 
\begin{equation}
\left\| \overline{z} \right\| \lesssim
\left\|\overline{x}\right\|
\end{equation}
and
\begin{equation}
\left\| \mathscr{A}\overline{z} -\overline{y} \right\| \lesssim
\epsilon \left\|\mathscr{A}\right\|\left\|\overline{x}\right\|.
\end{equation}

%
%
%
%
\end{corollary}

\end{subsection}

\end{section}

\begin{section}{Analysis of the Levin equation}
\label{section:leveq}

This sections contains our analysis of the Levin equation.
The two principal results are Theorems~\ref{leveq:thm1} and \ref{leveq:thm2}.
Theorem~\ref{leveq:thm1}
generalizes the classical result of \cite{Levin} by
showing that, whenever $g'$ is nonzero over the interval
 and the ratio of the maximum value of $g'$ to the minimum value of $g'$ is small,
the Levin equation admits a well-behaved solution which can be approximated by 
a polynomial expansion at a cost which is independent of the magnitude
of $g'$.    Theorem~\ref{leveq:thm2} shows that, on an interval in which $g'$ is of small
magnitude, the Levin equation admits a well-behaved solution  which can be represented
by a polynomial expansion at a cost which decreases with the magnitude of $g'$.
It applies whether or not $g'$ has zeros in the interval.
Throughout, we use the ``$\lesssim$ notation'' to suppress constants which do not depend on the magnitude of $g'$.
We begin with the following lemma, which applies in the simple case when $g'$ is a nonzero constant.

\begin{lemma}
\label{leveq:lemma1}
Suppose that $f:[-1,1]\to\mathbb{C}$ admits an infinitely differentiable extension
to a neighborhood of $[-1,1]$ and $W \neq 0$.  Then for each $0 < \epsilon < 1$,
there  exists a function $p_b \in S\left(\mathbb{R}\right)$ such that 
\begin{enumerate}
\item
$\begin{aligned}
\widehat{p_b}
\end{aligned}$
is a tempered distribution of order one supported on $\left[-c_f\left(\epsilon\right),c_f\left(\epsilon\right)\right]$,

\item
$\begin{aligned}
\left|p_b'(x) + iW p_b(x) \right| \leq \epsilon \|f\|_{L^\infty\left([-1,1]\right)}
\end{aligned}$
for all $x \in [-1,1]$, 

\item
$\begin{aligned}
\left\|p_b\right\|_{L^\infty\left([-1,1]\right)}
\lesssim
\min\left\{1,\frac{1}{|W|}\right\} \|f\|_{L^\infty\left(\left[-1,1\right]\right)}
\end{aligned}$
and

\item
$\begin{aligned}
\left\|p_b'\right\|_{L^\infty\left([-1,1]\right)}
\lesssim 
\min\left\{1,\frac{1}{|W|}\right\}
\|f\|_{L^\infty\left(\left[-1,1\right]\right)}.
\end{aligned}$

\end{enumerate}
\end{lemma}
\begin{proof}
We let $f_b \in S\left(\mathbb{R}\right)$ be a function with bandlimit $W_0=c_f(\epsilon)$ such that conditions (1)-(5) in 
Definition~\ref{preliminaries:approx:def} hold, and  define the function $p_b$ via the formula
\begin{equation}
p_b(x) = \mbox{p.v.}\, \int_{-W_0}^{W_0} \frac{\widehat{f_b}(\xi)}{i(W+\xi)} \exp(i\xi x)\,d\xi.
\label{leveq:101}
\end{equation}
That is, $p_b$ is the  inverse Fourier transform of the product of the tempered distribution 
$T$ defined via
\begin{equation}
\left\langle T, \varphi\right\rangle =  \mbox{p.v.}\, \int_{-\infty}^\infty \frac{\varphi(\xi)}{i\left(W+\xi\right)}\, d\xi
\label{leveq:102}
\end{equation}
and the infinitely differentiable function $\widehat{f_b}$.  
It is clear that $\widehat{p_b}$  is a tempered distribution
of order one which is supported on $[-W_0,W_0]$, so the first of the conditions
listed above is satisfied.

Since
\begin{equation}
p_b'(x) + iW p_b(x) = 
\int_{-W_0}^{W_0} \frac{i \left(W+\xi\right)\widehat{f_b}(\xi)}{i(W+\xi)}  \exp(i\xi x)\,d\xi
= 
f_b(x)
\label{leveq:103}
\end{equation}
and
\begin{equation}
\|f - f_b\|_{L^\infty\left(\left[-1,1\right]\right)} \leq \epsilon \|f\|_{L^\infty\left(\left[-1,1\right]\right)},
\label{leveq:104}
\end{equation}
the second of the properties of $p_b$ listed above holds.

To establish the third of the properties of $p_b$ listed above, we first assume that $|W| \leq W_0+1$.
Then we have 
\begin{equation}
\begin{aligned}
p_b(x) &= \mbox{p.v.}\, \int_{-W_0}^{W_0} \frac{\widehat{f_b}(\xi)}{i(W+\xi)}\exp(ix\xi) \,d\xi\\
&= \mbox{p.v.}\, \int_{-W-2W_0-1}^{-W+2W_0+1} \frac{\widehat{f_b}(\xi)}{i(W+\xi)}\exp(ix\xi)\,d\xi\\
&= \, \int_{-W-2W_0-1}^{-W+2W_0+1} \frac{\widehat{f_b}(\xi)\exp(ix\xi) - \widehat{f_b}(-W)\exp(-ixW)}{i(W+\xi)}\,d\xi\\
\end{aligned}
\label{leveq:105}
\end{equation}
for all $x \in \mathbb{R}$, where the last equality follows from the fact that $1/(W+\xi )$ is odd about
$\xi=-W$.  By the mean value theorem, there exists a function $\eta(\xi)$ such that 
\begin{equation}
\begin{aligned}
p_b(x) = \, \int_{-W-2W_0-1}^{-W+2W_0+1} \frac{G'(\eta(\xi)) (W+\xi)}{i(W+\xi)}\, d\xi
=
-i  \int_{-W-2W_0-1}^{-W+2W_0+1} G'(\eta(\xi))\, d\xi
\end{aligned}
\label{leveq:106}
\end{equation}
for all $x \in \mathbb{R}$, where
\begin{equation}
G(\xi) = \widehat{f_b}(\xi) \exp(ix\xi).
\label{leveq:107}
\end{equation}
Since
\begin{equation}
G'(\xi) = \widehat{f_b}'(\xi) \exp(ix\xi) + i x\widehat{f_b}(\xi) \exp(ix\xi),
\label{leveq:108}
\end{equation}
we have 
\begin{equation}
\left| p_b(x) \right|
\leq 
2 (2W_0 + 1)
\left(\left\|\widehat{f_b}'\right\|_{L^\infty\left(\mathbb{R}\right)}+ \left\|\widehat{f_b}\right\|_{L^\infty\left(\mathbb{R}\right)}\right)
\label{leveq:109}
\end{equation}
for all $x \in [-1,1]$ (note that we are now only considering  $x$ in $[-1,1]$).  Now making use of properties (4) and (5) 
in Definition~\ref{preliminaries:approx:def} gives us 
\begin{equation}
\left \|p_b\right\|_{L^\infty\left([-1,1]\right)}
\leq
16 (2 W_0 +1 )  \left\|f \right\|_{L^\infty\left([-1,1]\right)}.
\label{leveq:110}
\end{equation}
  An analogous argument which makes use of the fact that
\begin{equation}
p_b'(x) 
= \, \int_{-W-2W_0-1}^{-W+2W_0+1} \frac{H(\xi) - H(-W)}{i(W+\xi)}\, d\xi,
\label{leveq:111}
\end{equation}
where
\begin{equation}
H(\xi) = i \xi \widehat{f_b}(\xi) \exp(i x \xi ),
\label{leveq:112}
\end{equation}
shows that 
\begin{equation}
\begin{aligned}
\left \|p_b'\right\|_{L^\infty\left([-1,1]\right)}
&\leq
8 (2 W_0 + 1) \left(
1 + 2(|W|+2W_0+1)
\right) 
\left\|f\right\|_{L^\infty\left([-1,1]\right)}\\
&\leq 
8
\left(2W_0+1\right)
\left(
6W_0 + 5
\right)
\left\|f\right\|_{L^\infty\left([-1,1]\right)}\\
\end{aligned}
\label{leveq:113}
\end{equation}
when $|W| \leq W_0+1$.

We now suppose that  $|W| > W_0+1$.  Then
\begin{equation}
\begin{aligned}
\left|p_b(x)\right| &= \left|\int_{-W_0}^{W_0} \frac{\widehat{f_b}(\xi) \exp(i\xi x)}{i(W + \xi)}\, d\xi\right| 
&\leq \frac{2W_0}{\left|W\right|-\left|W_0\right|}\left\|\widehat{f_b}\right\|_{L^1\left(\mathbb{R}\right)}
\end{aligned}
\label{leveq:114}
\end{equation}
and
\begin{equation}
\begin{aligned}
\left|p_b'(x)\right| &= \left|\int_{-W_0}^{W_0} \frac{\widehat{f_b}(\xi) \exp(i\xi x)i\xi}{i(W + \xi)}\, d\xi\right| 
&\leq \frac{2W_0^2}{\left|W\right|-\left|W_0\right|}\left\|\widehat{f_b}\right\|_{L^1\left(\mathbb{R}\right)}
\end{aligned}
\label{leveq:115}
\end{equation}
hold for all $x \in \mathbb{R}$.  From the above inequalities, the fact that 
\begin{equation}
\frac{1}{|W|-W_0} \leq \frac{W_0 +1}{|W|}
\label{leveq:116}
\end{equation}
and property (4) in Definition~\ref{preliminaries:approx:def}, we see that
\begin{equation}
\left\|p_b\right\|_{L^\infty\left(\mathbb{R}\right)}
\leq 
 \frac{8W_0 (W_0 +1)}{|W|}  \left\| f\right\|_{L^\infty\left([-1,1]\right)}
\label{leveq:117}
\end{equation}
and
\begin{equation}
\left\|p_b'\right\|_{L^\infty\left(\mathbb{R}\right)}
\leq 
 \frac{8W_0^2 (W_0 +1)}{|W|}  \left\| f\right\|_{L^\infty\left([-1,1]\right)}
\label{leveq:118}
\end{equation}
whenever $|W| > W_0+1$.    
The third of the conditions on $p_b$ listed in the conclusion of the lemma follows from
the combination of (\ref{leveq:110}) and  (\ref{leveq:117}), while
the fourth is obtained by  combining (\ref{leveq:113}) and  (\ref{leveq:118}).
\end{proof}

In accordance with Remark~\ref{preliminaries:measure:remark1},
the sequence $\{a_m\}$ of Legendre coefficients of the approximate solution $p_b$ of Levin's equation
appearing in Lemma~\ref{leveq:lemma1} decays superexponentially once $m > e/2 c_f(\epsilon) + 1/2 $.
Consequently, $p_b$ can be represented to a fixed relative precision
via a polynomial expansion at a cost which depends
on the complexity of $f$ but not the magnitude of $W$.

We now move on to the case in which $g'$ is nonconstant, but with no zeros on the interval $[-1,1]$.
To that end, we suppose that $f:[-1,1] \to \mathbb{C}$ and $g:[-1,1] \to \mathbb{R}$ admit infinitely
differentiable  extensions to a neighborhood of $[-1,1]$, and that the extension of $g'$ is nonzero
in an open interval containing $[-1,1]$.  We also let 
\begin{equation}
G_0 = \min_{-1 \leq x \leq 1} \left|g'(x)\right|,\ \ \ 
W = \frac{1}{2} \int_{-1}^1 g'(x)\, dx
\label{leveq:199}
\end{equation}
and define $u:[-1,1] \to [-1,1]$  via the formula
\begin{equation}
u(x) = -1 + \frac{1}{W} \int_{-1}^x g'(y)\ dy.
\end{equation}
Because $g'$ is nonzero in a neighborhood of $[-1,1]$, $u$ is invertible and its inverse
extends to an open neighborhood of $[-1,1]$.  
Finally, we define $h:[-1,1] \to [-1,1]$ via the formula
\begin{equation}
h(z) = \frac{f\left(u^{-1}\left(z\right)\right)}{u'\left(u^{-1}\left(z\right)\right)}
=
f\left(u^{-1}\left(z\right)\right) \frac{du^{-1}}{dz}(z).
\end{equation}
%

With the preceding notations and assumptions, we have the following:
\begin{theorem}
\label{leveq:thm1}
For every $0 < \epsilon < 1$,  there exists a function $p_b:[-1,1] \to \mathbb{C}$ such that
\begin{enumerate}
\item
the Fourier transform of $p_b\left(u^{-1}(z)\right)$ 
is a tempered distribution of order $1$ supported on the interval
$\left[-c_h(\epsilon),c_h(\epsilon)\right]$,

\item
$\begin{aligned}
\left|p_b'(x) + i g'(x) p_b(x) - f(x) \right|
\leq \epsilon 
\frac{|W|}{G_0}
\left\|f\right\|_{L^\infty\left([-1,1]\right)}
\end{aligned}$
for all $x \in [-1,1]$,

\item
$\begin{aligned}
\left\|p_b\right\|_{L^\infty\left([-1,1]\right)}
\lesssim
\frac{|W|}{G_0}\min\left\{1,\frac{1}{|W|}\right\}
\left\|f\right\|_{L^\infty\left([-1,1]\right)}
\end{aligned}$
and

\item
$\begin{aligned}
\left\|p_b'\right\|_{L^\infty\left([-1,1]\right)}
\lesssim
\frac{|W|}{G_0}\min\left\{1,\frac{1}{|W|}\right\}
\left\|f\right\|_{L^\infty\left([-1,1]\right)}.
\end{aligned}$

\end{enumerate}
\end{theorem}

\begin{proof}
By introducing the new variable  $z = u(x)$, we transform the Levin equation
\begin{equation}
p'(x) + i g'(x) p (x)= f(x),\ \ \ \ \ -1 < x < 1,
\label{leveq:201}
\end{equation}
into the simplified form
\begin{equation}
p'(z) + i W p(z) = h(z), \ \ \ \ \ -1 < z < 1.
\label{leveq:202}
\end{equation}
Our assumptions on $f$ and $g$, and  the method we used to construct $u$ ensure
that $h$ admits an infinitely differentiable extension to a neighborhood of $[-1,1]$.
Applying Lemma~\ref{leveq:lemma1} to (\ref{leveq:202}) shows that,
for all $0 < \epsilon < 1$, there exists an entire  function $p_1$ such that
\begin{enumerate}
\item
$\begin{aligned}
\widehat{p_1}
\end{aligned}$
is a tempered distribution of order $1$ supported on the interval $\left[-c_{h}\,(\epsilon),c_{h}\,(\epsilon)\right]$,

\item
$\begin{aligned}
\left|p_1'(z) + i W p_1(z) - h(z)\right|
\leq \epsilon  \|h\|_{L^\infty\left([-1,1]\right)}
\end{aligned}$
for all $z \in [-1,1]$,

\item
$\begin{aligned}
\left\|p_1\right\|_{L^\infty\left([-1,1]\right)}
\lesssim
\min\left\{1,\frac{1}{|W|}\right\}
\|h\|_{L^\infty\left(\left[-1,1\right]\right)}
\end{aligned}$ and

\item
$\begin{aligned}
\left\|p_1'\right\|_{L^\infty\left([-1,1]\right)}
\lesssim
\min\left\{1,\frac{1}{|W|}\right\}
\|h\|_{L^\infty\left(\left[-1,1\right]\right)}.
\end{aligned}$

\end{enumerate}
%

%
%

We now define $p_b$ via the formula $p_b(x) = p_1(u(x))$.
It is clear that the first condition on $p_b$ listed above is satisfied.
We observe that, since $g'(x)=Wu'(x)$,
\begin{equation}
\begin{aligned}
\left\|h\right\|_{L^\infty\left([-1,1]\right)}
\leq \frac{|W|}{G_0}\left\|f\right\|_{L^\infty\left([-1,1]\right)}.
\end{aligned}
\label{leveq:204}
\end{equation}
Combining (\ref{leveq:204}) and properties of $p_1$ listed above yields
\begin{equation}
\left\|p_b\right\|_{L^\infty\left([-1,1]\right)}
\lesssim
\frac{|W|}{G_0} \min\left\{1,\frac{1}{|W|}\right\}
\|f\|_{L^\infty\left(\left[-1,1\right]\right)},
\label{leveq:205}
\end{equation}
\begin{equation}
\left\|p_b'\right\|_{L^\infty\left([-1,1]\right)}
\lesssim
\frac{|W|}{G_0} \min\left\{1,\frac{1}{|W|}\right\}
\|f\|_{L^\infty\left(\left[-1,1\right]\right)}
\label{leveq:206}
\end{equation}
and
\begin{equation}
\left| p_b(x) - i g'(x) p_b(x) - f(x) \right| \leq \epsilon \frac{|W|}{G_0} \|f\|_{L^\infty\left([-1,1]\right)},
\label{leveq:207}
\end{equation}
the latter of which holds for all $x \in [-1,1]$.  This suffices to establish the theorem.
\end{proof}


We emphasize that Theorem~\ref{leveq:thm1} does not imply that 
the approximate solution $p_b$ of the Levin equation
has a Fourier transform which is a tempered distribution with compact support.
Instead, we have
\begin{equation}
p_b(x) = \left\langle
\widehat{p_1}, \eta_x
\right\rangle
\ \ \mbox{with} \ \ \eta_x(\xi) = \exp(i u(x) \xi),
\label{leveq:208}
\end{equation}
where $\widehat{p_1}$ is a tempered distribution of order one with compact support
(c.f., Formula~(A3) in the original Levin paper \cite{Levin}).
Since $p_b$ is the composition of the entire function $p_1$ and
the function $u$, which is infinitely differentiable
on a neighborhood of $[-1,1]$, the magnitudes of the
coefficients $\{a_m\}$ in its Legendre expansion decay faster 
than $m^{-k}$ for any positive integer $k$.
Moreover, (\ref{leveq:208}) implies that  $p_b$ 
can be approximated to a fixed relative precision via a Legendre expansion
at a cost which depends on the complexity of $f$ and $g^{-1}$, but not
on the magnitude of $g'$.  
To see this, we first note that 
$h$ is defined in terms of $f$ and the normalized version $u$ of $g$, so that
$c_h(\epsilon)$ does not depend on the magnitude of $g'$.  The bandlimit of $p_1$
then depends on $c_h(\epsilon)$, and $p_b$ is obtained by the composition
of $p_1$ with the normalized version $u$ of $g$.


It would be of interest to develop better estimates 
on the rate of decay of the Legendre coefficients of $p_b$, 
as we did in the case when $g'$ is constant,
but  useful results of this type seem to be fairly complicated.
Indeed, the difficulty in  estimating the 
complexity of $p_b$ {\it a priori} is one of the principal motivations for the adaptive
version of the Levin algorithm that we introduce in this article.


We close this section with the following theorem which shows that, if $g'$ is of small magnitude,
then the Levin equation admits a solution which can be well-approximated by a bounded, bandlimited
function $p_b$. 

\begin{theorem}
\label{leveq:thm2}
Suppose that $f:[-1,1] \to \mathbb{C}$ and $g:[-1,1]\to\mathbb{R}$ admit extensions
to infinitely differentiable functions in a neighborhood of $[-1,1]$, and that 
\begin{equation}
G_1 = \max_{-1\leq x \leq 1} \left|g'(x)\right| < \frac{1}{2}.
\end{equation}
Let  $0 < \epsilon < 1$ be given, and define the integer $n$ via the formula
\begin{equation}
n = \left\lfloor 
\frac{\log(\epsilon)}{\log(2G_1)}
\right\rfloor.
\label{leveq:300}
\end{equation}
Then there exists $p_b \in C^\infty\left(\mathbb{R}\right)$ such that
\begin{enumerate}
\item
$\begin{aligned}\widehat{p_b}\end{aligned}$
is a tempered distribution of order $1$ supported on the interval  
\begin{equation}
\left[-c_f(\epsilon)- n c_{g'}(\epsilon), c_f(\epsilon) + n c_{g'}(\epsilon)\right],
\end{equation}

\item
$\begin{aligned}\left|p_b'(x) + i g'(x) p_b(x) - f(x) \right| \leq 
2\epsilon
\left(1+\frac{G_1}{1-2G_1}\right)
  \left\|f\right\|_{L^\infty\left([-1,1]\right)}
\end{aligned}$ 
for all $x \in [-1,1]$,

\item
$\begin{aligned}
\left\|p_b\right\|_{L^\infty\left([-1,1]\right)} \leq \frac{2}{1-2G_1} \left\|f\right\|_{L^\infty\left([-1,1]\right)}
\end{aligned}$ 
and

\item
$\begin{aligned}
\left\|p_b'\right\|_{L^\infty\left([-1,1]\right)} \leq 
4 \left(1+ \frac{G_1}{1-2G_1} \right)
\left\|f\right\|_{L^\infty\left([-1,1]\right)}.
\end{aligned}$ 

\end{enumerate}
\end{theorem}
\begin{proof}
We let $f_b$ and $g_b'$ be bandlimited approximates of $f$ and $g'$ which 
satisfy conditions (1)-(5) in Definition~\ref{preliminaries:approx:def}.
We next define
$A:L^\infty\left([-1,1]\right) \to L^\infty\left([-1,1]\right)$ via the formula
\begin{equation}
A\left[\varphi\right](x) = -i \int_0^x g_b'(y) \varphi(y)\,dy.
\label{leveq:301}
\end{equation}
Obviously, $\|A\|_\infty$  is bounded by $\|g_b'\|_{L^\infty\left([-1,1]\right)}$,
which is less than  $(1+\epsilon)G_1 \leq 2 G_1 < 1$ by property (1) in Definition~\ref{preliminaries:approx:def}.
Now we let $G_2 = (1+\epsilon)G_1$,
\begin{equation}
h(x) = \int_0^x f_b(y)\,dy
\label{leveq:302}
\end{equation}
and
\begin{equation}
p_b(x) = \sum_{k=0}^n A^k \left[ h \right](x),
\label{leveq:303}
\end{equation}
where $A^k$ denotes the repeated application of the operator $A$ and $n$ is defined by (\ref{leveq:300}).
Since $\|f_b\|_{L^\infty\left([-1,1]\right)} \leq (1+\epsilon) \|f\|_{L^\infty\left([-1,1]\right)}$ by
property (1) of Definition~\ref{preliminaries:approx:def}, we have
\begin{equation}
\begin{aligned}
\left\|p_b\right\|_{L^\infty\left([-1,1]\right)} 
&\leq
\left(\sum_{k=0}^n  \left\|A\right\|_\infty^k \right) \left\|h\right\|_{L^\infty\left([-1,1]\right)}\\
&\leq \frac{1+G_2^{n+1}}{1-G_2}  \|f\|_{L^\infty\left([-1,1]\right)}\\
&\leq \frac{1+\epsilon}{1-G_2}  \|f\|_{L^\infty\left([-1,1]\right)},
\end{aligned}
\label{leveq:304}
\end{equation}
which implies the third property of $p_b$ listed above.
Since $h'(x) = f_b(x)$ and 
\begin{equation}
\frac{d}{dx} A^k\left[h\right](x)  = -i g_b'(x) A^{k-1}\left[h\right](x),
\label{leveq:305}
\end{equation}
 we have 
\begin{equation}
\begin{aligned}
p_b'(x) &= f_b(x) - i g_b'(x) \sum_{k=0}^{n-1} A^k\left[h\right](x)\\
&= f_b(x)- i g_b'(x) \sum_{k=0}^{n} A^k\left[h\right](x) + i g_b '(x) A^n\left[h\right](x) \\
&= f_b(x)- ig_b'(x) p_b(x) + i g_b '(x) A^n\left[h\right](x) 
\end{aligned}
\label{leveq:306}
\end{equation}
for all $x\in [-1,1]$.  From this and  (\ref{leveq:304}), we see that
\begin{equation}
\begin{aligned}
\left\|p_b'(x)\right\|_{L^\infty\left([-1,1]\right)}
&\leq 
\left\|f_b\right\|_{L^\infty\left([-1,1]\right)}
+ 
\left\|g'\right\|_{L^\infty\left([-1,1]\right)} 
\left\|p_b\right\|_{L^\infty\left([-1,1]\right)} \\
&\ \ \ \ \ +
G_2^{n+1} \left\|f_b\right\|_{L^\infty\left([-1,1]\right)}
\\
&\leq 
\left((1+\epsilon) 
+ G_1 \frac{1+\epsilon}{1-G_2}
+\epsilon (1+\epsilon)
\right)\left\|f\right\|_{L^\infty\left([-1,1]\right)}\\
&\leq 
\left((1+\epsilon) ^2+ \frac{2G_1}{1-G_2}
\right)\left\|f\right\|_{L^\infty\left([-1,1]\right)},
\end{aligned}
\label{leveq:307}
\end{equation}
from which  the fourth of the properties of $p_b$ listed in the conclusion of the theorem
follows.  Similarly,  (\ref{leveq:306}) implies that 
\begin{equation}
\begin{aligned}
\left| p_b'(x) + i g_b'(x) p_b(x) - f(x) \right| 
&\leq
\left| p_b'(x) + i g_b'(x) p_b(x) - f_b(x) \right|  + 
\left| f_b(x) - f(x) \right|  
\\
&\leq
\left|
i g_b '(x) A^n\left[h\right](x) 
\right|
+\left| f_b(x) - f(x) \right|  
\\
&\leq
G_2^{n+1}  \left\|f\right\|_{L^\infty\left([-1,1]\right)}
+ \epsilon \left\|f\right\|_{L^\infty\left([-1,1]\right)}\\
&\leq
2 \epsilon \left\|f\right\|_{L^\infty\left([-1,1]\right)}
\end{aligned}
\label{leveq:308}
\end{equation}
for all $x \in [-1,1]$.  
Combining (\ref{leveq:308}) with  (\ref{leveq:304}), (\ref{leveq:307}) and using the properties of the approximates $f_b$ 
and $g_b$, we see that 
\begin{equation}
\begin{aligned}
&\left| p_b'(x) + i g'(x) p_b(x) - f(x) \right|\\
\leq
&\left| p_b'(x) + i g_b'(x) p_b(x) - f(x)   \right| + \left| ig'(x) p_b(x) - i g_b'(x) p_b(x) \right|\\
\leq
&2\epsilon  \left\|f\right\|_{L^\infty\left([-1,1]\right)}+
\left\|p_b\right\|_{L^\infty\left([-1,1]\right)}
\left\|g'-g_b'\right\|_{L^\infty\left([-1,1]\right)}
\\
\leq
&\left(
2 
+
\left(\frac{1+\epsilon}{1-G_2}
\right)
\left(
 \left\|g'\right\|_{L^\infty\left([-1,1]\right)}\right)
\right)
\epsilon 
\left\|f\right\|_{L^\infty\left([-1,1]\right)}
\end{aligned}
\label{leveq:310}
\end{equation}
for all $x \in [-1,1]$.     The second of the properties of $p_b$
listed above follows from (\ref{leveq:310}).

It remains only to establish the first of the properties of $p_b$ listed above.
To do so, we first observe that the Fourier transform of $h$ is 
\begin{equation}
\mbox{p.v.}\, \frac{\widehat{f_b}(\xi)}{i \xi} 
+ 
\frac{1}{2} \delta(\xi) 
\left(\int_{0}^\infty f_b(x)\,dx  - \int_{-\infty}^0 f_b(x)\,dx,
\right),
\end{equation}
which is a tempered distribution of order one supported  on the interval $\left[-c_\epsilon(f),c_\epsilon(f)\right]$.
Next, we suppose that $\psi$ is a $C^\infty\left(\mathbb{R}\right)$ function
whose Fourier transform is a tempered distribution of order $1$ supported on $[-c,c]$.
The Fourier transform of the  convolution of $\psi$ with the function $g_b' \in S\left(\mathbb{R}\right)$ 
is a tempered distribution of order $0$ supported on $\left[-c-c_{g'},c+c_{g'}\right]$.
Since  $A\left[\psi\right]$ is obtained by integrating this convolution, 
its  Fourier transform 
 is a tempered distribution of order $1$ with bandlimit
$\left[-c-c_{g'}(\epsilon),c+c_{g'}(\epsilon)\right]$.  It now follows by induction that
the Fourier transform of $A^{k}\left[h\right]$ is a tempered
distribution of order $1$ with bandlimit 
$\left[-c_f(\epsilon)-k c_{g'}(\epsilon),c_f(\epsilon)+ kc_{g'}(\epsilon)\right]$.
Combining this with  (\ref{leveq:303}) yields
the first property of $p_b$ listed above.
\end{proof}

The bandlimit of the function $p_b$ 
whose existence is established in Theorem~\ref{leveq:thm2} decreases with 
magnitude of $g'$.
Consequently, on any interval on which $g'$ is sufficiently small, we expect to
be able to  represent $p_b$ to a fixed relative accuracy with a polynomial expansion
whose number of terms is bounded.  Moreover, this is the case whether or not $g'$
has zeros in the interval.

\end{section}

\begin{section}{Numerical aspects of the Levin method}
\label{section:levnum}

In this section, we show that high-accuracy can be obtained when
the  Levin equation (\ref{introduction:ode})  is discretized via a Chebyshev spectral collocation method and 
the resulting linear system is solved with a truncated singular value decomposition, regardless of the magnitude of $g'$
and whether or not it has zeros.
To that end, we suppose that $f:~[-1,1]\to\mathbb{C}$ and $g:[-1,1]\to\mathbb{R}$
admit continuously differentiable extensions to a neighborhood of $[-1,1]$, 
and that $0 < \epsilon <1$.  Moreover, we let
\begin{equation}
G_0 = \min_{-1 \leq x \leq 1} \left|g'(x)\right|,\ \ \ G_1 = \max_{-1 \leq x \leq 1} \left|g'(x)\right|\ \ \
\mbox{and}\ \ \   W = \frac{1}{2} \int_{-1}^1 g'(x)\, dx.
\label{levnum:1}
\end{equation}

We first consider the case in which $G_0>0$ and invoke
Theorem~\ref{leveq:thm1} to see  that there exists a bandlimited function $p_b$ such
that
\begin{equation}
\left|p_b'(x) +i g'(x) p_b(x) - f(x) \right| \leq 
\epsilon \frac{|W|}{G_0} \|f\|_{L^\infty\left([-1,1]\right)}
\ \ \mbox{for all}\ \  x \in [-1,1],
\label{levnum:2}
\end{equation}
\begin{equation}
\|p_b\|_{L^\infty\left([-1,1]\right)} \lesssim \frac{|W|}{G_0} \min\left\{1,\frac{1}{|W|}\right\}
\|f\|_{L^\infty\left([-1,1]\right)}
\label{levnum:3}
\end{equation}
and
\begin{equation}
\|p_b'\|_{L^\infty\left([-1,1]\right)} \lesssim \frac{|W|}{G_0} \min\left\{1,\frac{1}{|W|}\right\}
\|f\|_{L^\infty\left([-1,1]\right)}.
\label{levnum:4}
\end{equation}
Moreover, it is clear from the discussion following Theorem~\ref{leveq:thm1} that 
\begin{equation}
p_b(x) = \sum_{n=0}^\infty a_n T_n(x),
\label{levnum:5}
\end{equation}
where $|a_n|$ is  bounded by a rapidly decaying function of $n$ which is independent  of the magnitude of $g'$.

We claim that there is an integer $k$  such that
\begin{equation}
\begin{aligned}
&\ 
\left\|
\mathsf{P}_k\left[p_b\right]' + i \mathsf{P}_k \left[g'\right] \mathsf{P}_k \left[p_b\right]
- \left( p_b' + i g' p_b\right)
\right\|_{L^\infty\left([-1,1]\right)}
\\[0.5em]
\leq &\ 
\epsilon \|p_b'\|_{L^\infty\left([-1,1]\right)}
+
 4 \epsilon G_1 \|p_b\|_{L^\infty\left([-1,1]\right)}\\[0.5em]
\lesssim &\ 
\epsilon \left(
1  + 4 G_1\right)
\frac{|W|}{G_0} \min\left\{1, \frac{1}{|W|}\right\} \|f\|_{L^\infty\left([-1,1]\right)}
\end{aligned}
\label{levnum:6}
\end{equation}
holds regardless of the values of $|W|$ and  $G_1$; that is, $k$ can be chosen independent of the magnitude of $g'$.
To see this, we first observe that (\ref{levnum:5}) together with Formulas~(\ref{preliminaries:chebound1}) and (\ref{preliminaries:chebound2}) 
in Subsection~\ref{section:preliminaries:chebyshev}
imply that  we can choose $k$ independently of $G_1$ and $|W|$ such that 
\begin{equation}
\left\| \mathsf{P}_k\left[p_b\right] - p_b \right\|_{L^\infty\left([-1,1]\right)} \leq \epsilon
\left\| p_b \right\|_{L^\infty\left([-1,1]\right)}
\label{levnum:7}
\end{equation}
and
\begin{equation}
\left\| \mathsf{P}_k\left[p_b\right]' -  p_b' \right\|_{L^\infty\left([-1,1]\right)} 
\leq  \epsilon \left\| p_b' \right\|_{L^\infty\left([-1,1]\right)}.
\label{levnum:8}
\end{equation}
Clearly, we can also ensure that 
\begin{equation}
\left\| \mathsf{P}_k\left[g'\right] - g' \right\|_{L^\infty\left([-1,1]\right)} \leq \epsilon
\left\| g' \right\|_{L^\infty\left([-1,1]\right)}
\label{levnum:9}
\end{equation}
and
\begin{equation}
\left\| \mathsf{P}_k\left[f\right] -f \right\|_{L^\infty\left([-1,1]\right)} \leq \epsilon
\left\| f \right\|_{L^\infty\left([-1,1]\right)}
\label{levnum:9.5}
\end{equation}
with $k$ still chosen independently of $G_1$ and $|W|$.
Now  (\ref{levnum:7}) and (\ref{levnum:9}) together with the assumption that
$0 < \epsilon < 1$ imply
\begin{equation}
\left\| \mathsf{P}_k[p_b] \right\|_{L^\infty\left([-1,1]\right)}
\leq 2 \left\| p_b \right\|_{L^\infty\left([-1,1]\right)}
\end{equation}
and
\begin{equation}
\left\| \mathsf{P}_k[g'] \right\|_{L^\infty\left([-1,1]\right)}
\leq 2 \left\| g' \right\|_{L^\infty\left([-1,1]\right)}.
\label{levnum:10}
\end{equation}
It follows readily that 
\begin{equation}
\begin{aligned}
&\left\| \mathsf{P}_k\left[g'\right]  \mathsf{P}_k\left[p_b\right] - g' p_b\right\|_{L^\infty\left([-1,1]\right)}\\
\leq 
&\left\| \mathsf{P}_k\left[p_b\right]\right\|_{L^\infty\left([-1,1]\right)}
\left\| \mathsf{P}_k\left[g'\right]-g'\right\|_{L^\infty\left([-1,1]\right)}\\ 
&\ \ \ \ \ + \left\| \mathsf{P}_k\left[g'\right]\right\|_{L^\infty\left([-1,1]\right)} 
\left\| \mathsf{P}_k\left[p_b\right]-p_b\right\|_{L^\infty\left([-1,1]\right)}\\
\leq
&4\epsilon  \|g'\|_{L^\infty\left([-1,1]\right)} \|p_b\|_{L^\infty\left([-1,1]\right)}.
\end{aligned}
\label{levnum:11}
\end{equation}
Now (\ref{levnum:8}) and (\ref{levnum:11})  imply the first inequality in (\ref{levnum:6}),
and the second follows from (\ref{levnum:3}) and (\ref{levnum:4}).

Combining (\ref{levnum:2}) with (\ref{levnum:6}) shows that 
\begin{equation}
\begin{aligned}
&\left| \mathsf{P}_k\left[p_b\right]'(x) 
+\mathsf{P}_k\left[g'\right](x)  \mathsf{P}_k\left[p_b\right](x) - f(x)\right|\\
\lesssim
&
\epsilon
\left(1+\left(1  + 4 G_1 \right)\min\left\{1, \frac{1}{|W|}\right\}\right)
\frac{|W|}{G_0} \|f\|_{L^\infty\left([-1,1]\right)}
\end{aligned}
\label{levnum:12}
\end{equation}
for all $x\in [-1,1]$.    If we now let
\begin{equation}
\overline{p_b} = 
\left(\begin{array}{c}
p_b\left(x^{\mbox{\tiny cheb}}_{1,k}\right)\\
p_b\left(x^{\mbox{\tiny cheb}}_{2,k}\right)\\
\vdots\\
p_b\left(x^{\mbox{\tiny cheb}}_{k,k}\right)\\
\end{array}\right)
\ \ \ \mbox{and} \ \ \ 
\overline{f}=
\left(\begin{array}{c}
f\left(x^{\mbox{\tiny cheb}}_{1,k}\right)\\
f\left(x^{\mbox{\tiny cheb}}_{2,k}\right)\\
\vdots\\
f\left(x^{\mbox{\tiny cheb}}_{k,k}\right)\\
\end{array}\right)
\label{levnum:13}
\end{equation}
be the vectors of values of $p_b$ and $f$ at the $k$-point extremal Chebyshev grid,
and define the $k\times k$ matrix $\mathscr{G}$ via
\begin{equation}
\mathscr{G} = \left(
\begin{array}{cccccc}
g'\left(x^{\mbox{\tiny cheb}}_{1,k}\right) & & & \\
& g'\left(x^{\mbox{\tiny cheb}}_{2,k}\right) & & & \\
& & \ddots \\
& & & &g'\left(x^{\mbox{\tiny cheb}}_{k,k}\right)
\end{array}
\right),
\label{levnum:14}
\end{equation}
 then (\ref{levnum:12}) implies that 
\begin{equation}
\left( \mathscr{D}_k + i \mathscr{G}\right) \overline{p_b} = \overline{f} + \overline{\delta },
\label{levnum:15}
\end{equation}
where
\begin{equation}
\left\|\overline{p_b}\right\|
\lesssim 
\frac{|W|}{G_0} \min\left\{1,\frac{1}{|W|}\right\}
\|f\|_{L^\infty\left([-1,1]\right)}
\label{levnum:16}
\end{equation}
and
\begin{equation}
\left\| 
\overline{\delta}
\right\|
\lesssim
\epsilon
\left(1+\left(1  + 4 G_1 \right)\min\left\{1, \frac{1}{|W|}\right\}\right)
\frac{|W|}{G_0} \|f\|_{L^\infty\left([-1,1]\right)}.
\label{levnum:17}
\end{equation}
From (\ref{levnum:16}), (\ref{levnum:17})  and the fact that 
\begin{equation}
\|\mathscr{D}_k + i \mathscr{G}  \| \lesssim \max\{G_1, k^2\},
\label{levnum:18}
\end{equation}
we see that 
\begin{equation}
\left\|\overline{\delta}\right\|
\lesssim \epsilon \left\|\mathscr{D}_k + i \mathscr{G}  \right\| \left\|\overline{p_b}\right\|
\label{levnum:19}
\end{equation}
regardless of the magnitude of $g'$.  In particular, the left-hand side
of (\ref{levnum:19}) is on the order
of $(G_1/G_0 + |W|/G_0) \sim G_1/G_0$ when $G_1$ is large and $G_1/G_0|W|$
when $G_1$ is small, while
the right-hand side is on the order of $G_1/G_0$ when $G_1$ is large
and $G_1/G_0 k^2 > G_1/G_0|W|$ when $G_1$ is small.
So by  Corollary~\ref{preliminaries:numerical:svd}, when we solve the linear system
\begin{equation}
\left(\mathscr{D}_k + i \mathscr{G} \right)\overline{x} = \overline{f}
\label{levnum:20}
\end{equation}
via a singular value decomposition which as been truncated at precision on the order of $\epsilon \left\|\mathscr{D}_k + \mathscr{G} \right\|$,
we obtain a solution  $\overline{p_1}$ such that 
\begin{equation}
\left\| \overline{p_1} \right\| 
\lesssim \left\|\overline{p_b} \right\|
\lesssim
\frac{|W|}{G_0} \min\left\{1,\frac{1}{|W|}\right\}
\|f\|_{L^\infty\left([-1,1]\right)}
\label{levnum:21}
\end{equation}
and
\begin{equation}
\begin{aligned}
\left\| \left(\mathscr{D}_k + \mathscr{G} \right)\overline{p_1} - \overline{f}\right\| 
&\lesssim
\epsilon \left\|\mathscr{D}_k + \mathscr{G} \right\| \left\|\overline{p_b}\right\|\\
&\lesssim
\epsilon 
\frac{|W|}{G_0}  \max\{G_1, k^2\}  \min\left\{1,\frac{1}{|W|}\right\}
\|f\|_{L^\infty\left([-1,1]\right)}.
\end{aligned}
\label{levnum:22}
\end{equation}

Now we let $p_1$ and $\delta_1$ be the $k$-term Chebyshev expansions
whose values at the extremal Chebyshev nodes agree with the values of the vector
$\overline{p_1}$ and those of $\left(\mathscr{D}_k + \mathscr{G} \right)\overline{p_1}-\overline{f}$, respectively.
From (\ref{levnum:21}) and (\ref{levnum:22}) and the fact that the Chebyshev polynomials are bounded in $L^\infty\left([-1,1]\right)$
norm by $1$, we have that
\begin{equation}
\left\|
p_1
\right\|_{L^\infty\left([-1,1]\right)}
\lesssim
\frac{|W|}{G_0} \min\left\{1,\frac{1}{|W|}\right\}
\|f\|_{L^\infty\left([-1,1]\right)}
\label{levnum:23}
\end{equation}
and
\begin{equation}
\left\|
\delta_1 
\right\|_{L^\infty\left([-1,1]\right)}
\lesssim
\epsilon  
\frac{|W|}{G_0}  \max\{G_1, k^2\}  \min\left\{1,\frac{1}{|W|}\right\}
\|f\|_{L^\infty\left([-1,1]\right)}.
\label{levnum:24}
\end{equation}
Moreover, it follows from  (\ref{levnum:20}) that 
%
\begin{equation}
p_1'\left(x_{j,k}^{\mbox{\tiny cheb}}\right) 
+ i g'\left(x_{j,k}^{\mbox{\tiny cheb}}\right) 
p_1\left(x_{j,k}^{\mbox{\tiny cheb}}\right) = 
f\left(x_{j,k}^{\mbox{\tiny cheb}}\right)  + \delta_1\left(x_{j,k}^{\mbox{\tiny cheb}}\right)
\label{levnum:25}
\end{equation}
for all $j=1,\ldots,k$; that is,
\begin{equation}
\mathsf{P}_k
\left[ p_1' + i g'p_1 \right](x) 
= \mathsf{P}_k\left[f  +\delta_1\right](x).
\label{levnum:26}
\end{equation}
But using (\ref{levnum:9}), (\ref{levnum:9.5}), (\ref{levnum:23})  and the fact that $p_1$ and $\delta_1$ are  $k$-term Chebyshev expansions, 
we see  that 
\begin{equation}
\begin{aligned}
&\left\|
\mathsf{P}_k\left[p_1' + ig' p_1\right] - \left(p_1'+ig' p_1\right)
\right\|_{L^\infty\left([-1,1]\right)}\\
&\leq 
\epsilon G_1 \|p_1\|_{L^\infty\left([-1,1]\right)}\\
&\lesssim
\epsilon\frac{G_1}{G_0} |W| \min\left\{1,\frac{1}{|W|}\right\}
\|f\|_{L^\infty\left([-1,1]\right)}
\label{levnum:27}
\end{aligned}
\end{equation}
and
\begin{equation}
\left\|
\mathsf{P}_k\left[f+\delta_1\right] - \left(f+\delta_1\right)
\right\|_{L^\infty\left([-1,1]\right)}
\leq 
\epsilon \|f\|_{L^\infty\left([-1,1]\right)}.
\label{levnum:28}
\end{equation}
It follows from (\ref{levnum:24}), (\ref{levnum:26}), (\ref{levnum:27}) and (\ref{levnum:28}) that 
\begin{equation}
\begin{aligned}
&\left| p_1'(x) + i g'(x)p_1(x) - f(x) \right| \\
\lesssim
&\epsilon
\left(1+
\frac{G_1}{G_0} +
\frac{1}{G_0}  \max\{G_1, k^2\}  
 \right)
 |W| \min\left\{1,\frac{1}{|W|}\right\} \left\|f\right\|_{L^\infty\left([-1,1]\right)}
\end{aligned}
\label{levnum:29}
\end{equation}
for all $x\in [-1,1]$.  

To complete our analysis for the case in which $G_0 > 0$, we let 
\begin{equation}
I = \int_{-1}^1 f(t) \exp(ig(t))\ dt
\label{levnum:30}
\end{equation}
be the true value of the oscillatory integral we hope to compute via the Levin method,
and take $I_1$ to be the estimate
\begin{equation}
I_1 = p_1(1) \exp(i g(1)) - p_1(-1)\exp(ig(-1)).
\label{levnum:31}
\end{equation}
%
Now
\begin{equation}
\begin{aligned}
I_1- I      &= \int_{-1}^1 \frac{d}{dt} \left( p_1(t) \exp(ig(t))\right)\,dt - \int_{-1}^1 f(t) \exp(ig(t))\,dt \\
        &=\int_{-1}^1 \left(p_1'(t) + i g'(t) p_1(t)-f(t)\right)   \exp(i g(t)) \,dt,
\end{aligned}
\label{levnum:32}
\end{equation}
and it follows from this and (\ref{levnum:29}) that 
\begin{equation}
\begin{aligned}
&\left|I_1 -I  \right|   \\
\lesssim
&\epsilon
\left(1+
\frac{G_1}{G_0}  +\frac{1}{G_0}  \max\{G_1, k^2\}
 \right)
|W|\min\left\{1,\frac{1}{|W|}\right\} \left\|f\right\|_{L^\infty\left([-1,1]\right)}.
\end{aligned}
\label{levnum:33}
\end{equation}
The quantity
\begin{equation}
|W|  \min\left\{1,\frac{1}{|W|}\right\}
\end{equation}
is clearly  bounded independent of the magnitude of $g'$.
Moreover, $G_1/G_0$ is small  when $g'$ does not vary in magnitude too much over the interval (a reasonable assumption
when analyzing an adaptive scheme).  So 
(\ref{levnum:33}) shows that the absolute
error in the computed integral is bounded independent of the magnitude of $g'$,
assuming $g'$ does not vary too much over the interval.

We now consider the case in which $G_1 < 1/4$.  By invoking Theorem~\ref{leveq:thm2}, we 
can find a bandlimited function $p_b$  such that
\begin{equation}
\left|p_b'(x) +i g'(x) p_b(x) - f(x) \right| \leq 
2\epsilon
\left(1+\frac{G_1}{1-2G_1}\right)
  \left\|f\right\|_{L^\infty\left([-1,1]\right)}
\label{levnum:101}
\end{equation}
for all $x \in [-1,1]$,
\begin{equation}
\|p_b\|_{L^\infty\left([-1,1]\right)} \leq
 \frac{2}{1-2G_1}
\|f\|_{L^\infty\left([-1,1]\right)}
\label{levnum:102}
\end{equation}
and
\begin{equation}
\|p_b'\|_{L^\infty\left([-1,1]\right)} \leq
4 \left(1+ 
\frac{G_1}{1-2G_1}
\right)
\|f\|_{L^\infty\left([-1,1]\right)}.
\label{levnum:103}
\end{equation}
Since the bandlimit of $p_b$ is bounded under our assumption on $G_1$, the Chebyshev
coefficients of $p_b$ are bounded by a rapidly decaying function which is
independent of $G_1$.  Moreover, we can once again choose $k$ independent of $G_1$
such that  (\ref{levnum:7}) through (\ref{levnum:11}) hold.
Then, proceeding as we did before shows that 
\begin{equation}
\begin{aligned}
&\left\|
\mathsf{P}_k\left[p_b\right]' + i \mathsf{P}_k \left[g'\right] \mathsf{P}_k \left[p_b\right]
- \left( p_b' + i g' p_b\right)
\right\|_{L^\infty\left([-1,1]\right)}\\
\lesssim 
&\epsilon 
\left(1 
 + \frac{4G_1}{1-2G_1}
\right)
\|f\|_{L^\infty\left([-1,1]\right)}
\end{aligned}
\label{levnum:104}
\end{equation}
for all $x \in [-1,1]$.   If we define $\overline{p_b}$, $\overline{\delta}$ and $\mathscr{G}$ as before,
then we see that 
\begin{equation}
\left\|\overline{\delta}\right\| \lesssim
\epsilon 
\left(1 
 + \frac{4G_1}{1-2G_1}
\right)
\|f\|_{L^\infty\left([-1,1]\right)}
\label{levnum:105}
\end{equation}
while
\begin{equation}
\left\|\left(\mathscr{D}_k + \mathscr{G}\right) \overline{p_b}\right\|  \lesssim
 \frac{2}{1-2G_1}
\max\{G_1,k^2\}
\|f\|_{L^\infty\left([-1,1]\right)},
\label{levnum:106}
\end{equation}
so there is no difficulty in applying Corollary~\ref{preliminaries:numerical:svd} to
see that solving (\ref{levnum:20}) via a truncated singular value decompostion yields  a solution  $\overline{p_1}$
such that 
\begin{equation}
\left\| \overline{p_1} \right\| 
\lesssim
 \frac{2}{1-2G_1}
\|f\|_{L^\infty\left([-1,1]\right)}
\label{levnum:107}
\end{equation}
and
\begin{equation}
\begin{aligned}
\left\| \left(\mathscr{D}_k + \mathscr{G} \right)\overline{p_1} - \overline{f}\right\| 
&\lesssim
\epsilon \left\|\mathscr{D}_k + \mathscr{G} \right\| \left\|\overline{p_b}\right\|\\
&\lesssim
\epsilon 
  \max\{G_1, k^2\}  \frac{2}{1-2G_1}
\|f\|_{L^\infty\left([-1,1]\right)}.
\end{aligned}
\label{levnum:108}
\end{equation}
Defining $p_1$ and $\delta_1$ as before, we obtain the bound
\begin{equation}
\begin{aligned}
&\left| p_1' (x) + i g'(x) p_1(x) - f(x) \right|\\
\lesssim
&\epsilon 
\left(
  \max\{G_1, k^2\}  \frac{2}{1-2G_1}
+
 \frac{1}{1-2G_1}
\right)\|f\|_{L^\infty\left([-1,1]\right)},
\end{aligned}
\label{levnum:109}
\end{equation}
which holds for all $x \in [-1,1]$.
Arguing as before, we see that the absolute error in the computed value of the integral
(\ref{levnum:31})  is bounded by a constant multiple of 
\begin{equation}
\epsilon 
\left(
  \max\{G_1, k^2\}  \frac{2}{1-2G_1}
+
 \frac{1}{1-2G_1}
\right)\|f\|_{L^\infty\left([-1,1]\right)}.
\label{levnum:110}
\end{equation}
\end{section}
Since we are considering an interval on which $G_1$ is small, the constant 
in (\ref{levnum:110}) is small as well.

\begin{section}{The Adaptive Levin Method}
\label{section:algorithm}

In this section, we describe the adaptive Levin method
for the numerical calculation of the integral
\begin{equation}
\int_a^b \exp(i g(x)) f(x)\,dx.
\label{algorithm:int}
\end{equation}
Trivial modifictions allow for the evaluation of 
\begin{equation}
\int_a^b \cos(g(x)) f(x)\,dx\ \ \ \mbox{or} \ \ \
\int_a^b \sin(g(x)) f(x)\,dx
\end{equation}
instead.

Before detailing the algorithm proper, we describe a subroutine which estimates
the value of 
\begin{equation}
\int_{a_0}^{b_0} \exp(ig(x)) f(x)\,dx,
\label{algorithm:subint}
\end{equation}
where $[a_0,b_0]$ is a subinterval of $[a,b]$. It takes as input the interval $[a_0,b_0]$, 
an integer $k$ which controls the number of Chebyshev nodes used to discretize
the Levin equation  and an external subroutine for evaluating the functions $f$ and $g$.
The subroutine proceeds as follows:
\begin{enumerate}

\item Use the external subroutine supplied by the user to evaluate the 
functions $f$ and $g$ at the extremal Chebyshev nodes
$x_{1,k}^{\mbox{\tiny cheb}},\ x_{2,k}^{\mbox{\tiny cheb}},\ \ldots,\ x_{k,k}^{\mbox{\tiny cheb}}$.

\item 
Calculate approximate values
\begin{equation}
\widetilde{g'\left(x_{1,k}^{\mbox{\tiny cheb}}\right)},\ \ldots,\ \widetilde{g'\left(x_{k,k}^{\mbox{\tiny cheb}}\right)}, 
\end{equation}
of the derivatives of the function $g$ at the extremal Chebyshev nodes by applying
the spectral differential matrix $\mathscr{D}_k$ to the vector
of values of $g$; that is, via the formula
\begin{equation}
\left(\begin{array}{c}
\widetilde{g'\left(x_{1,k}^{\mbox{\tiny cheb}}\right)}\\[.5em]
\widetilde{g'\left(x_{2,k}^{\mbox{\tiny cheb}}\right)}\\[.5em]
\vdots\\
\widetilde{g'\left(x_{k,k}^{\mbox{\tiny cheb}}\right)}\\[.5em]
\end{array}
\right)
=
\mathscr{D}_k 
\left(\begin{array}{c}
g\left(x_{1,k}^{\mbox{\tiny cheb}}\right)\\[0.5em]
g\left(x_{2,k}^{\mbox{\tiny cheb}}\right)\\[0.5em]
\vdots\\[0.5em]
g\left(x_{k,k}^{\mbox{\tiny cheb}}\right)
\end{array}
\right).
\end{equation}

\item Form the matrix
\begin{equation}
\mathscr{A} =  \mathscr{D}_k + 
i \left(\begin{array}{ccccccc}
\widetilde{g'\left(x_{1,k}^{\mbox{\tiny cheb}}\right)} &       &        &\\
        & \widetilde{g'\left(x_{2,k}^{\mbox{\tiny cheb}}\right)} &        &\\
        &         & \ddots \\&
        &         &          & \widetilde{g'\left(x_{k,k}^{\mbox{\tiny cheb}}\right)}
\end{array}
\right).
\end{equation}

\item Construct a singular value decomposition
\begin{equation}
\mathscr{A} = 
\left(
\begin{array}{cccc}
\overline{u_1} & \overline{u_2} & \cdots & \overline{u_k}
\end{array}  
\right)
\left(
\begin{array}{ccccc}
\sigma_1 &           &         & \\
         & \sigma_2  &         & \\
         &           & \ddots  & \\
         &           &         &\sigma_k \\
\end{array}
\right) 
\left(
\begin{array}{cccc}
\overline{v_1} & \overline{v_2} & \cdots & \overline{v_k}
\end{array}  
\right)^*
\end{equation}
of the matrix $\mathscr{A}$.  

\item
Find the least integer $1\leq l \leq k$ such that $\sigma_l \ge \|\mathscr{A}\| \epsilon_0$, where $\epsilon_0$
is machine zero.  If no such integer exists, return the estimate $0$ for (\ref{algorithm:subint}).

\item Let
\begin{equation*}
\left(
\begin{array}{c}
\widetilde{p\left(x_{1,k}^{\mbox{\tiny cheb}}\right)}\\[1.1em]
\widetilde{p\left(x_{2,k}^{\mbox{\tiny cheb}}\right)}\\[1.1em]
\vdots \\
\widetilde{p\left(x_{k,k}^{\mbox{\tiny cheb}}\right)}\\[1.1em]
\end{array}
\right)
  = 
\left(
\begin{array}{ccccc}
\overline{v_1} & \cdots  & \overline{v_l}
\end{array}
\right)
\left(
\begin{array}{ccccc}
\frac{1}{\sigma_1} &           &         & \\
         & \frac{1}{\sigma_2}  &         & \\
         &           & \ddots  & \\
         &           &         &\frac{1}{\sigma_l} \\
\end{array}
\right) 
\left(\begin{array}{cccccc}
\overline{u_1} & \cdots & \overline{u_l }
\end{array}\right)^*
\left(
\begin{array}{c}
f(x_1) \\
f(x_2) \\
\vdots\\
f(x_k)
\end{array}
\right).
\end{equation*}
%
The entries of this vector
approximate the values of a function $p$ such that 
\begin{equation}
\frac{d}{dx} \left( p(x) \exp(i g(x)) \right) = f(x) \exp(i g(x))
\end{equation}
at the extermal Chebyshev nodes on $[a_0,b_0]$. 

\item  
Return the estimate
\begin{equation}
\widetilde{p\left(x_{k,k}^{\mbox{\tiny cheb}}\right)}\exp(i g(x_k)) -  \widetilde{p\left(x_{1,k}^{\mbox{\tiny cheb}}\right)}\exp(i g(x_1))
\label{algorithm:est0}
\end{equation}
for the value of (\ref{algorithm:subint}).
\end{enumerate}

The algorithm proper takes as input a tolerance parameter $\epsilon >0$, 
the endpoints $a < b$ of the integration domain, an integer $k$ specifying the number
of Chebyshev nodes used to discretize the Levin equation on each subinterval considered,
 and an external subroutine which returns the values
of the functions $f$ and $g$ at a specified collection of points.
It maintains an estimated value $val$ for (\ref{algorithm:int}) and a list of intervals.
Initially, the list of intervals contains only $[a,b]$ and the value of the estimate is set to $0$. 
As long as  the list of intervals is nonempty the following steps are repeated:
\begin{enumerate}
\item 
Remove an entry $[a_0,b_0]$ from the list of intervals.
\item
Calculate an estimate $val_0$ of 
\begin{equation}
 \int_{a_0}^{b_0} \exp(ig(x)) f(x)\,dx
\end{equation}
using the subprocedure described above.

\item
Calculate estimates $val_{\mbox{\tiny L}}$ and $val_{\mbox{\tiny R}}$
of 
\begin{equation}
 \int_{a_0}^{c_0} \exp(ig(x)) f(x)\,dx\ \ \mbox{and} \ \ \
 \int_{c_0}^{b_0} \exp(ig(x)) f(x)\,dx,
\label{algorithm:phase}
\end{equation}
where $c_0 = (a_0+b_0)/2$,  using the subprocedure described above.

\item
If $\left| val_0 - val_{\mbox{\tiny L}} - val_{\mbox{\tiny R}}\right| < \epsilon$, then
update the estimate $val$ by letting  $val = val + val_0$.  Otherwise, add the intervals
$[a_0,c_0]$ and $[c_0,b_0]$ to the list of intervals.

\end{enumerate}
In the end, the procedure returns the estimate $val$ for (\ref{algorithm:int}).
Because the condition number of the oscillatory integral (\ref{algorithm:int})
increases with the magnitude of $g$, some loss of accuracy is expected
when calculating its value numerically.  In the case of the adaptive Levin method, 
the principal loss of accuracy occurs when exponentials of large magnitude are evaluated 
in  (\ref{algorithm:est0}).  Of course, the magnitudes of many 
integrals of the form (\ref{algorithm:int}) decrease with the magnitude $g$,
with the consequence that the absolute error in the calculated
value of (\ref{algorithm:int}) often remains constant or even decays as the  magntiude
of $g$ increases.

\vskip 1em
\begin{remark}
The truncated singular value decomposition is quite expensive.  In our implementation
of the adaptive Levin method, we used a rank-revealing QR decomposition in lieu
of the truncated singular value decomposition to 
solve the linear system which results from discretizing the ordinary differential equation.
This was found to be about 5 times faster and leads to no apparent loss in accuracy.
\end{remark}



\end{section}

\begin{section}{Numerical experiments}
\label{section:experiments}

In this section, we present the results of numerical experiments
which were conducted to illustrate the properties of the
algorithm of this article.  The code for these experiments
was written in Fortran and compiled with version 12.10
of the GNU Fortran compiler.  
They were performed on a desktop computer equipped 
with an AMD Ryzen 3900X processor and 32MB of memory.  No attempt 
was made to parallelize our code.    We used a 12-point
Chebyshev spectral method in our implementation of the 
adaptive Levin scheme (i.e., the parameter $k$ was taken to be $12$).   
As discussed in Section~\ref{section:algorithm}, the condition
number of evaluation of most integrals of the form 
\begin{equation}
\int_a^b \exp(i g(x)) f(x)\,dx
\label{experiments:expint}
\end{equation}
increases with the magnitude of $g'$.  This is often offset by a commensurate decrease
in the magnitude of the integral with the consequence that, in most cases,
it is reasonable to expect absolute errors in the calculated
values of (\ref{experiments:expint}) to be largely independent of the magnitude of $g'$.

In some of these experiments, we compared the results of the algorithm of this paper 
with those of an adaptive Gaussian quadrature code 
for evaluating integrals of the from
\begin{equation}
\int_a^b f(x) \,dx, \ \ -\infty < a < b <\infty.
\label{experiments:gauss0}
\end{equation}
Our implementation of this algorithm is written is Fortran and is quite standard.  
It maintains a list of intervals, which is initialized with the single interval
$[a,b]$, and a running tally  of the value of the integral, which is initialized to zero.
As long as the list of intervals is nonempty, the algorithm removes one interval $[a_0,b_0]$ from
the list, compares the value of 
\begin{equation}
\int_{a_0}^{b_0} f(x)\,dx
\label{experiments:gauss1}
\end{equation}
as computed by a $30$-point Gauss-Legendre quadrature rule to the value of the sum
\begin{equation}
\int_{a_0}^{(a_0+b_0)/2} f(x)\,dx +  \int_{(a_0+b_0)/2}^{b_0} f(x)\,dx,
\label{experiments:gauss2}
\end{equation}
where each integral is separately estimated with a 30-point Gauss-Legendre rule.
If the difference is larger than $\epsilon$, where $\epsilon$ is a tolerance parameter
supplied  by the user, then the intervals $\left[a_0,(a_0+b_0)/2\right]$
and  $\left[(a_0+b_0)/2,b_0\right]$ are placed in the list of intervals.
Otherwise, the value of (\ref{experiments:gauss1}) is added to the running
tally of the value of the integral (\ref{experiments:gauss0}).
In all of the experiments described here, the tolerance parameter for the adaptive Gaussian
quadrature code was taken to be $\epsilon = 10^{-15}$.  We found it necessary to set
the  tolerance parameter for
adaptive Gaussian quadrature  to be somewhat smaller than that 
for the adaptive Levin method in order to obtain accurate results from the former.
We used a 30-point Gauss-Legendre rule because we found it to be more efficient
than rules of other orders in most cases.

The experiments of Subsection~\ref{section:experiments:bes} concern Bessel functions, those
of Subsection~\ref{section:experiments:alf} involve the associated Legendre functions
and the experiments of Subsection~\ref{section:experiments:herm} concern
Hermite polynomials. In order to apply the adaptive Levin method  to integrals involving these functions,
we constructed phase function representations of them
via the method of \cite{BremerPhase2}.   That algorithm  applies to
second order linear ordinary differential equations of the form
\begin{equation}
y''(x) + q(x) y(x) = 0,\ \ \ a <x < b,
\label{experiments:second}
\end{equation}
where the coefficient $q$ is real-valued and slowly-varying.  We note that essentially  any second order
differential equation, including the differential equation defining the  associated Legendre functions
and Bessel's differential equation, can be put into the form (\ref{experiments:second})
via a simple transformation.   The method of \cite{BremerPhase2} constructs a piecewise Chebyshev expansion
representing a slowly-varying phase function $\psi$ such that
\begin{equation}
\frac{\exp\left(i\psi(x)\right)}{\sqrt{\psi'(x)}}
\ \ \ \mbox{and} \ \ \ 
\frac{\exp\left(-i\psi(x)\right)}{\sqrt{\psi'(x)}}
\end{equation}
constitute a basis in the space of solutions of (\ref{experiments:second}).
The derivative of the phase function is uniquely determined, but the phase function itself is only defined
up to a constant.  We typically use this degree of freedom to ensure that either
\begin{equation}
\frac{\sin\left(\psi(x)\right)}{\sqrt{\psi'(x)}}
\ \ \ \mbox{or}\ \ \
\frac{\cos\left(\psi(x)\right)}{\sqrt{\psi'(x)}}
\end{equation}
represent the special function we wish to integrate. 

\begin{subsection}{Certain integrals involving elementary functions for which explicit formulas are available}
\label{section:experiments:elem}

In the experiments described in this subsection, we considered the integrals
\begin{equation*}
\begin{aligned}
I_1(\lambda) &= \int_{-1}^1 \cos\left(\lambda \arctan(x) \right) \frac{1}{1+x^2}\,dx = \frac{2}{\lambda} \sin\left(\frac{\pi \lambda}{4}\right)
,\\
I_2(\lambda) &=\int_0^\infty  \frac{\exp\left(i \lambda x^2 \right)}{\sqrt{x}} \,dx =\exp\left(\frac{\pi i}{8}\right)
\frac{2 \Gamma\left(\frac{5}{4}\right)}{\lambda^{\frac{1}{4}}} ,\\
I_3(\lambda) &= \int_0^1 \exp\left(\frac{i \lambda}{\sqrt{x}}\right) \frac{1}{x}\,dx
= 2\Gamma(0,-i\lambda)
 \ \ \mbox{and} \\
I_4(\lambda) &= \int_0^{10} \exp(i\lambda \exp(x)) \exp(x)\,dx = 
\frac{i}{\lambda} \left(\exp(i\lambda) - \exp\left(i\exp\left(10\right)\lambda\right)\right).
\end{aligned}
\label{experiment:elem:ints}
\end{equation*}
All of the above formulas can be found in \cite{Gradshteyn}.

We proceeded by first sampling $l=200$ equispaced points $x_1,\ldots,x_{l}$  in the interval $[1,7]$.
Then, we used the adaptive Levin method to
evaluate the integrals appearing above  for each $\lambda = 10^{x_1}, 10^{x_2}, \ldots, 10^{x_l}$.
Figures~\ref{figure:elemplots1} and \ref{figure:elemplots2}
give the time require to evaluate these integrals and the absolute errors
in the obtained values.

\begin{figure}[h!]
\hfil
\includegraphics[width=.40\textwidth]{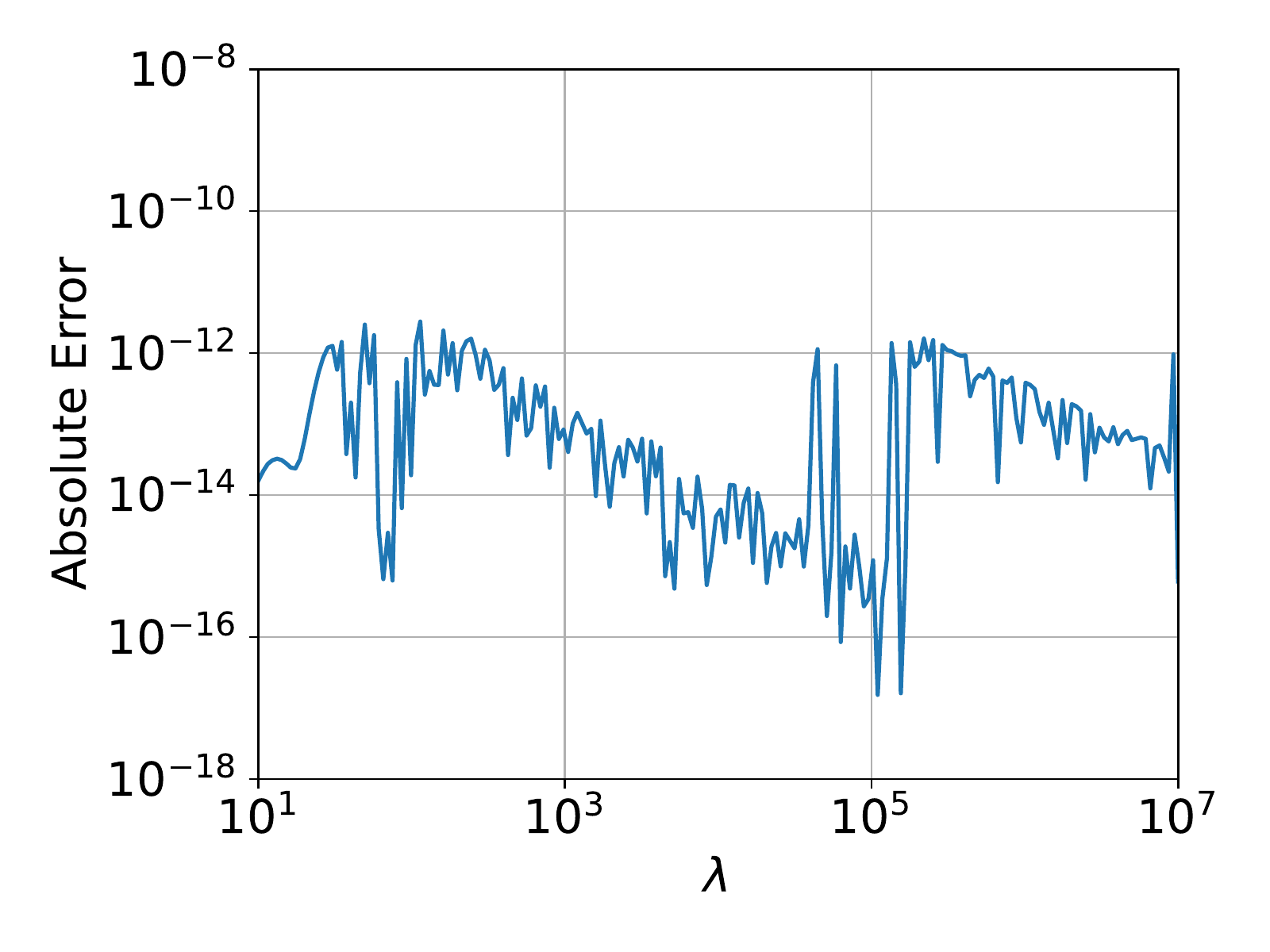}
\hfil
\includegraphics[width=.40\textwidth]{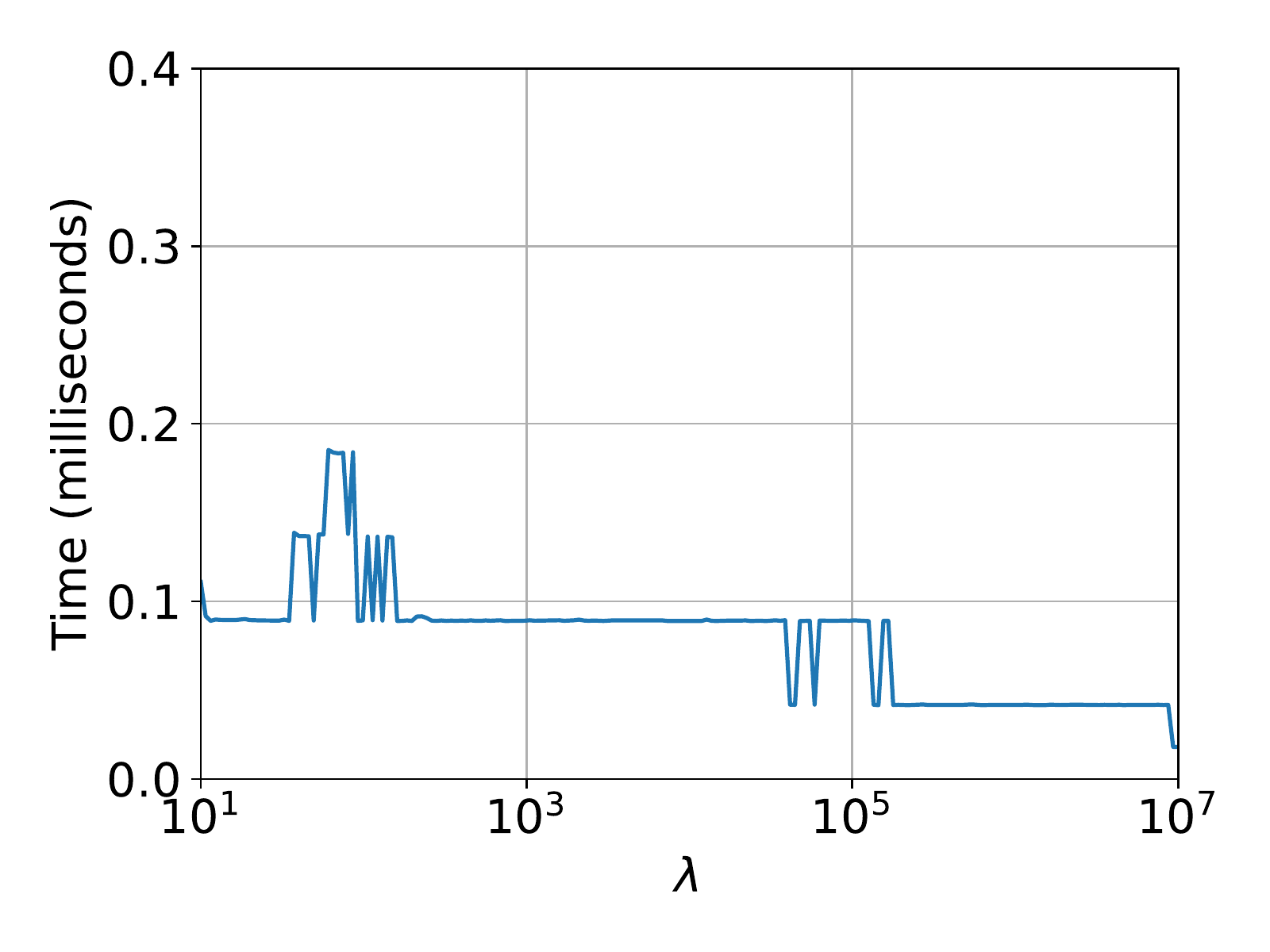}
\hfil

\hfil
\includegraphics[width=.40\textwidth]{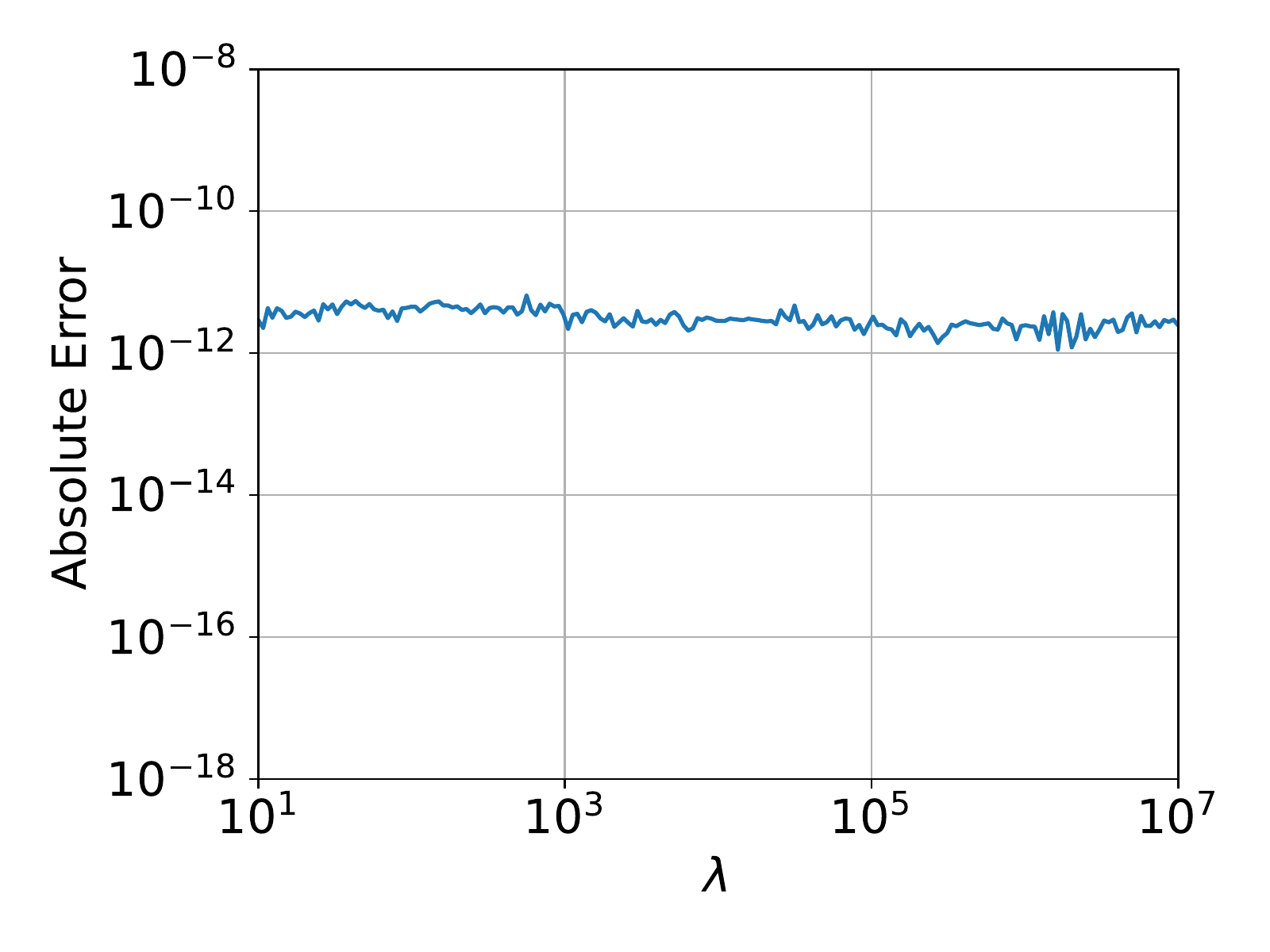}
\hfil
\includegraphics[width=.40\textwidth]{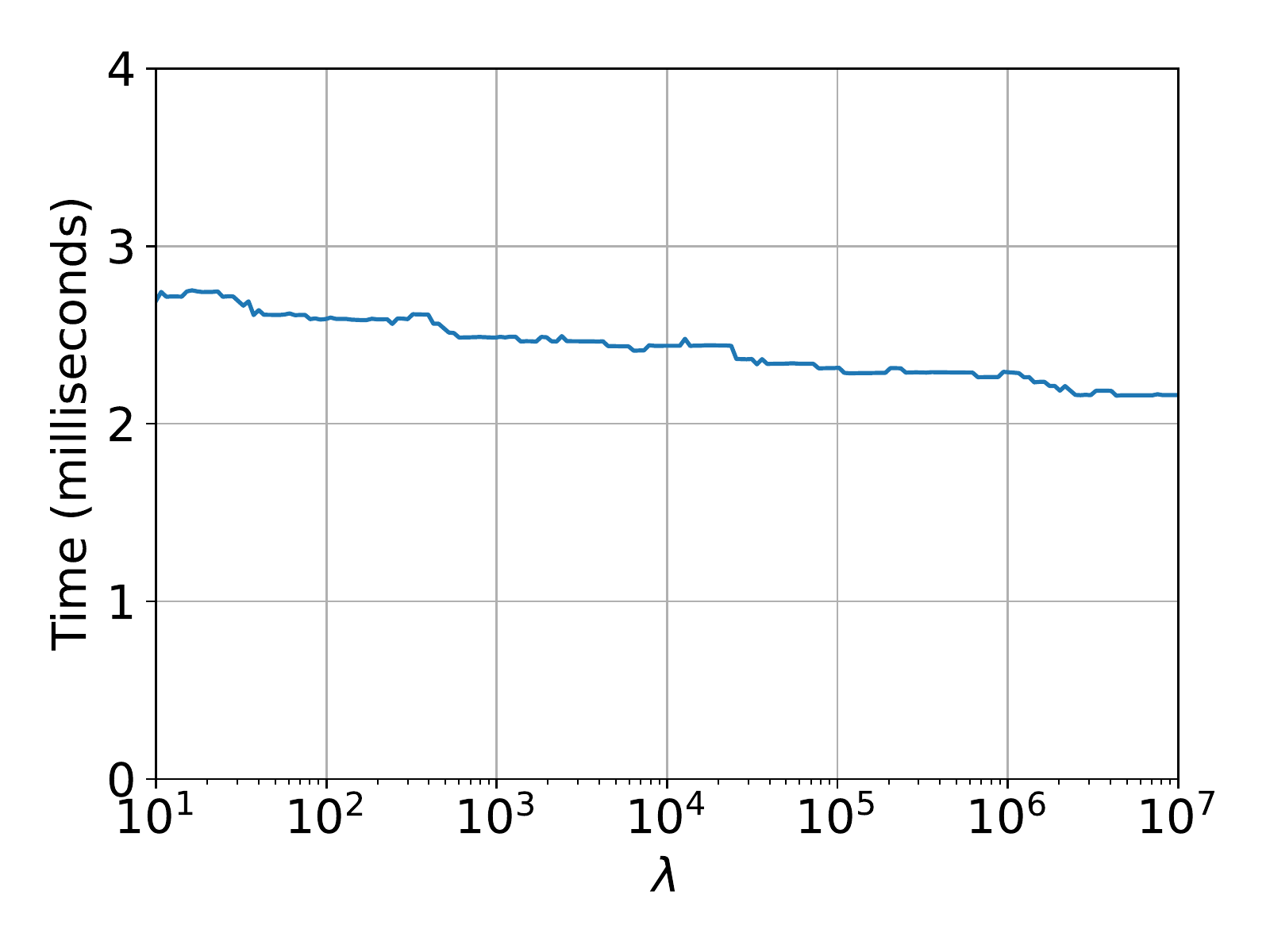}
\hfil

\caption{\small
The results of the first two experiments of Section~\ref{section:experiments:elem}.
The first row of plots pertains to the integral $I_1$ and  the second to $I_2$.
 Each plot on the left gives the error in the value of the integral
computed via the adaptive Levin method as a function of $\lambda$, while
each plot on the right gives the running time in milliseconds as a function
of $\lambda$.
}
\label{figure:elemplots1}
\end{figure}

\begin{figure}[h!]
\hfil
\includegraphics[width=.40\textwidth]{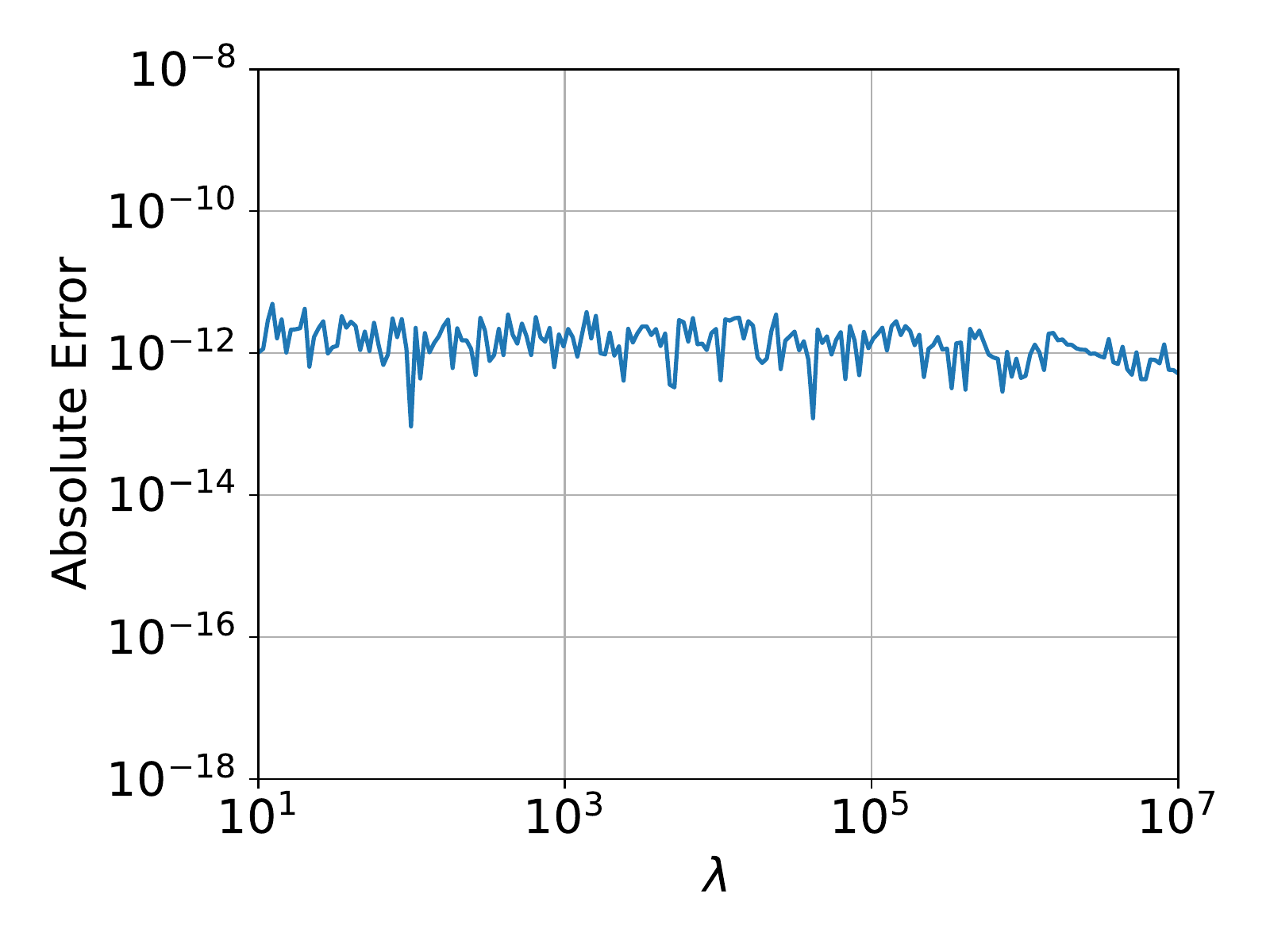}
\hfil
\includegraphics[width=.40\textwidth]{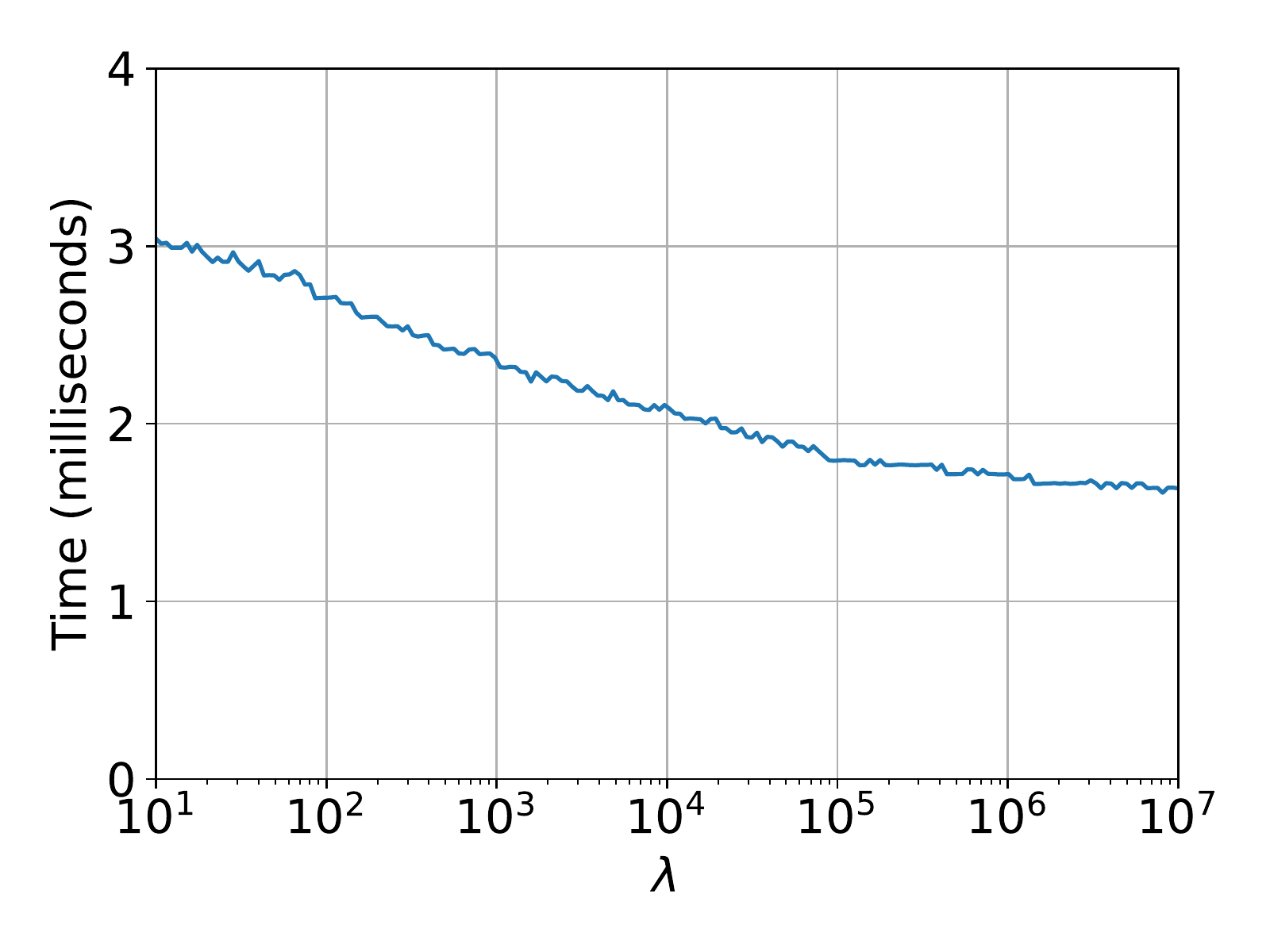}
\hfil

\hfil
\includegraphics[width=.40\textwidth]{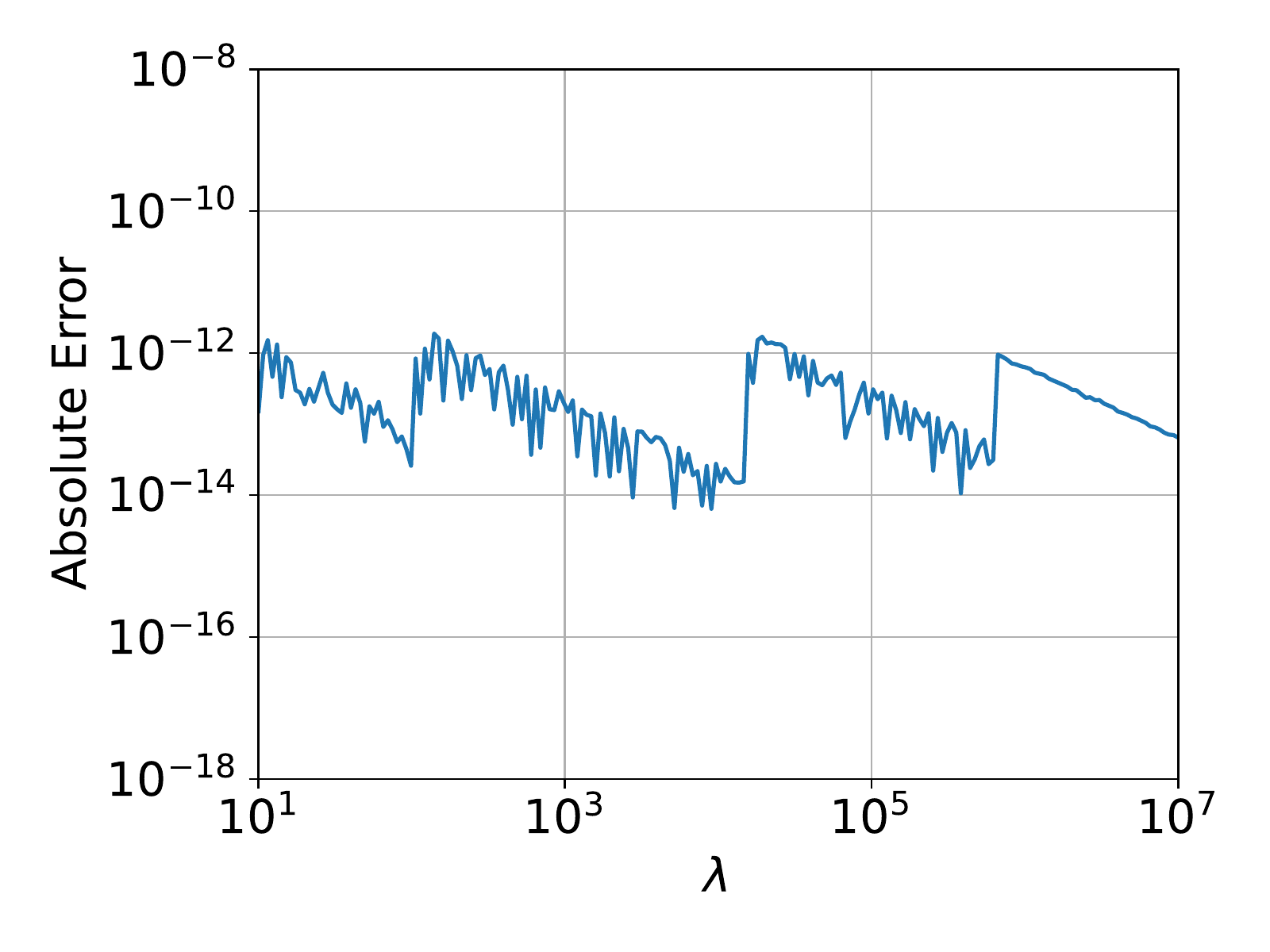}
\hfil
\includegraphics[width=.40\textwidth]{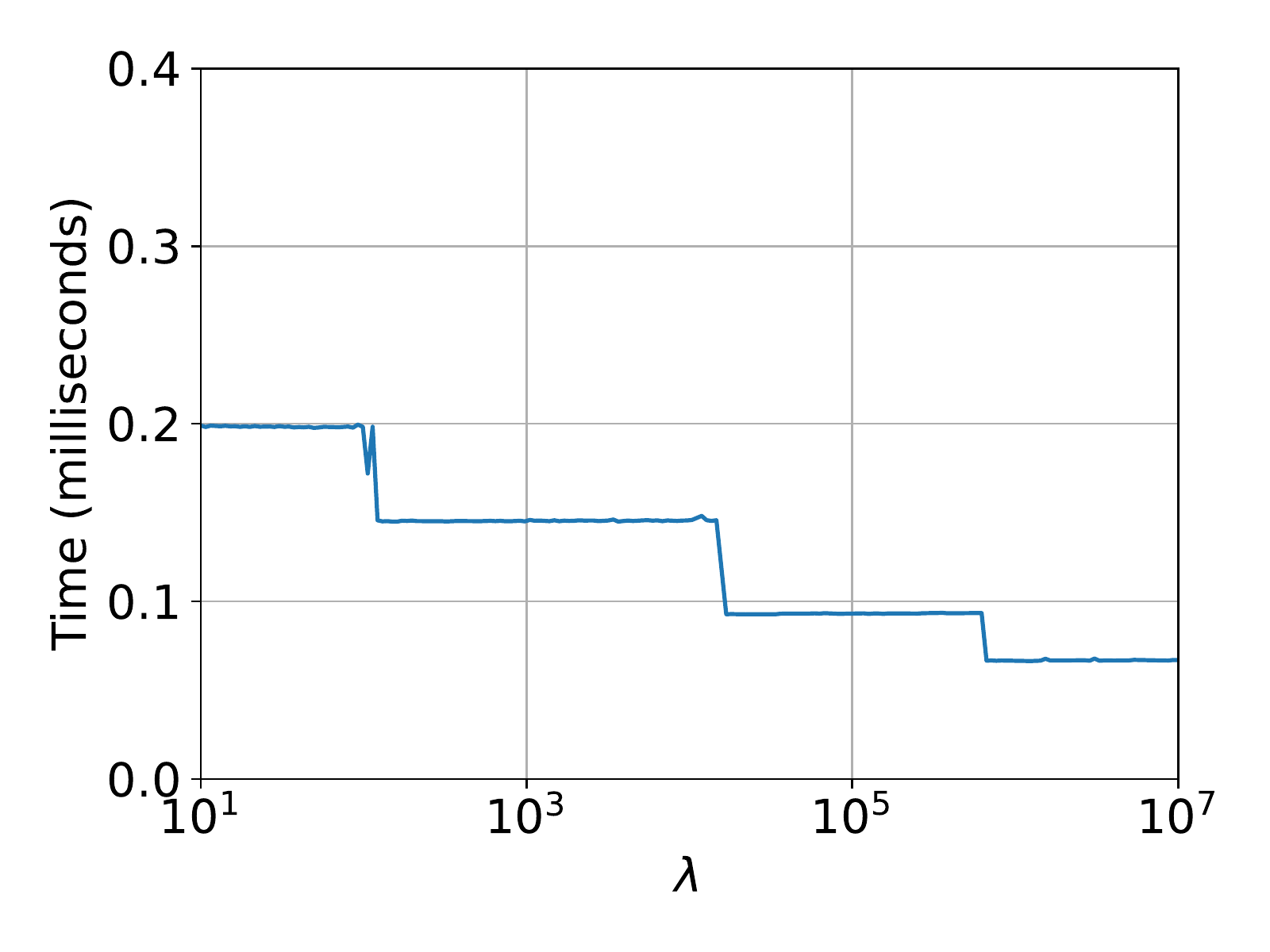}
\hfil

\caption{\small
The results of the last two experiments of Section~\ref{section:experiments:elem}.
The first row of plots pertains to the integral $I_3$ and the second to $I_4$.
 Each plot on the left gives the error in the value of the integral
computed via the adaptive Levin method as a function of $\lambda$, while
each plot on the right gives the running time in milliseconds as a function
of $\lambda$.
}
\label{figure:elemplots2}
\end{figure}

\end{subsection}

\begin{subsection}{Certain other integrals involving elementary functions}
\label{section:experiments:gauss}

In the experiments of this subsection, we evaluated the integrals
\begin{equation}
\begin{aligned}
I_5(\lambda) &= \int_0^1 \exp\left(i \lambda  x^2\right) \exp\left(-x\right) x\,dx,
\\
I_6(\lambda) &= \int_{-1}^{1} \exp\left( i\lambda x^2\right)  \left(1+x^2\right) \,dx,
\\
I_7(\lambda) &= \int_{-4}^4 \exp\left( i\lambda  x^2\right)\, dx
 \ \ \mbox{and} \\
I_8(\lambda) &= \int_{-1}^{1} \exp\left(i \lambda x^4\right) \frac{1}{0.01+x^4}\,dx.
\end{aligned}
\label{experiments:elem}
\end{equation}
We considered various ranges of values of $\lambda$.  We randomly sampled $200$ values of  $\lambda$ 
in each range and, then, for each such value of $\lambda$, we calculated
 $I_5$, $I_6$, $I_7$ and $I_8$ using both
the adaptive Levin algorithm and our adaptive Gaussian quadrature code.
The time taken by each method was measured, and the absolute differences
between the values of the integrals calculated with each method recorded.
Table~\ref{table:gauss:table1}  gives the result.  
Each row gives the results for one integral and range of values of $\lambda$.

\begin{table}[h!]
\center
\begin{tabular}{c@{\hspace{  2.00em}}c@{\hspace{  1.00em}}c@{\hspace{  1.00em}}c@{\hspace{  1.00em}}r@{\hspace{  1.00em}}c@{\hspace{  0.25em}}}
\toprule
\multirow[c]{2}{*}{Integral } & \multirow[c]{2}{*}{Range of $\lambda$ } & Avg Time & Avg Time & \multirow[c]{2}{*}{Ratio } & Max Observed \\ 
 &  & Adap Levin & Adap Gauss &  & Difference \\ 
\midrule 
\addlinespace[  0.50em]
$I_5$ & $10^{ 0}  - 10^{ 1}$ & $  4.83\times 10^{-05}$ & $  2.23\times 10^{-06}$ &       0.05 & $  9.94\times 10^{-13}$\\
 & $10^{ 1}  - 10^{ 2}$ & $  9.16\times 10^{-05}$ & $  5.93\times 10^{-06}$ &       0.06 & $  1.32\times 10^{-12}$\\
 & $10^{ 2}  - 10^{ 3}$ & $  1.23\times 10^{-04}$ & $  4.71\times 10^{-05}$ &       0.38 & $  1.01\times 10^{-12}$\\
 & $10^{ 3}  - 10^{ 4}$ & $  1.58\times 10^{-04}$ & $  4.15\times 10^{-04}$ &       2.64 & $  7.53\times 10^{-13}$\\
 & $10^{ 4}  - 10^{ 5}$ & $  2.02\times 10^{-04}$ & $  4.26\times 10^{-03}$ &      21.12 & $  9.99\times 10^{-13}$\\
 & $10^{ 5}  - 10^{ 6}$ & $  2.29\times 10^{-04}$ & $  3.99\times 10^{-02}$ &     173.87 & $  1.00\times 10^{-12}$\\
 & $10^{ 6}  - 10^{ 7}$ & $  2.51\times 10^{-04}$ & $  3.71\times 10^{-01}$ &    1476.94 & $  4.00\times 10^{-13}$\\
\addlinespace[  0.50em]
$I_6$ & $10^{ 0}  - 10^{ 1}$ & $  1.38\times 10^{-04}$ & $  2.29\times 10^{-06}$ &       0.02 & $  1.94\times 10^{-12}$\\
 & $10^{ 1}  - 10^{ 2}$ & $  2.88\times 10^{-04}$ & $  1.38\times 10^{-05}$ &       0.05 & $  1.97\times 10^{-12}$\\
 & $10^{ 2}  - 10^{ 3}$ & $  4.17\times 10^{-04}$ & $  1.07\times 10^{-04}$ &       0.26 & $  3.58\times 10^{-12}$\\
 & $10^{ 3}  - 10^{ 4}$ & $  4.74\times 10^{-04}$ & $  9.00\times 10^{-04}$ &       1.90 & $  3.32\times 10^{-12}$\\
 & $10^{ 4}  - 10^{ 5}$ & $  5.25\times 10^{-04}$ & $  8.84\times 10^{-03}$ &      16.85 & $  2.20\times 10^{-12}$\\
 & $10^{ 5}  - 10^{ 6}$ & $  5.77\times 10^{-04}$ & $  8.05\times 10^{-02}$ &     139.52 & $  3.53\times 10^{-12}$\\
 & $10^{ 6}  - 10^{ 7}$ & $  6.44\times 10^{-04}$ & $  8.01\times 10^{-01}$ &    1243.52 & $  2.57\times 10^{-12}$\\
\addlinespace[  0.50em]
$I_7$ & $10^{ 0}  - 10^{ 1}$ & $  3.61\times 10^{-04}$ & $  3.03\times 10^{-05}$ &       0.08 & $  2.57\times 10^{-12}$\\
 & $10^{ 1}  - 10^{ 2}$ & $  4.84\times 10^{-04}$ & $  2.44\times 10^{-04}$ &       0.51 & $  2.97\times 10^{-12}$\\
 & $10^{ 2}  - 10^{ 3}$ & $  5.32\times 10^{-04}$ & $  2.31\times 10^{-03}$ &       4.34 & $  3.67\times 10^{-12}$\\
 & $10^{ 3}  - 10^{ 4}$ & $  5.86\times 10^{-04}$ & $  2.24\times 10^{-02}$ &      38.15 & $  3.41\times 10^{-12}$\\
 & $10^{ 4}  - 10^{ 5}$ & $  6.35\times 10^{-04}$ & $  2.21\times 10^{-01}$ &     347.36 & $  2.52\times 10^{-12}$\\
 & $10^{ 5}  - 10^{ 6}$ & $  6.93\times 10^{-04}$ & $  2.31\times 10^{-00}$ &    3337.39 & $  3.29\times 10^{-12}$\\
 & $10^{ 6}  - 10^{ 7}$ & $  7.61\times 10^{-04}$ & $  2.66\times 10^{-01}$ &   34967.12 & $  5.68\times 10^{-12}$\\
\addlinespace[  0.50em]
$I_8$ & $10^{ 0}  - 10^{ 1}$ & $  8.66\times 10^{-04}$ & $  4.26\times 10^{-05}$ &       0.05 & $  3.48\times 10^{-12}$\\
 & $10^{ 1}  - 10^{ 2}$ & $  7.32\times 10^{-04}$ & $  1.70\times 10^{-04}$ &       0.23 & $  6.57\times 10^{-12}$\\
 & $10^{ 2}  - 10^{ 3}$ & $  7.58\times 10^{-04}$ & $  1.43\times 10^{-03}$ &       1.88 & $  4.17\times 10^{-12}$\\
 & $10^{ 3}  - 10^{ 4}$ & $  7.49\times 10^{-04}$ & $  1.36\times 10^{-02}$ &      18.09 & $  7.30\times 10^{-12}$\\
 & $10^{ 4}  - 10^{ 5}$ & $  7.41\times 10^{-04}$ & $  1.31\times 10^{-01}$ &     177.18 & $  6.40\times 10^{-12}$\\
 & $10^{ 5}  - 10^{ 6}$ & $  7.57\times 10^{-04}$ & $  1.19\times 10^{-00}$ &    1570.01 & $  3.62\times 10^{-12}$\\
 & $10^{ 6}  - 10^{ 7}$ & $  8.23\times 10^{-04}$ & $  1.37\times 10^{-01}$ &   16697.05 & $  3.76\times 10^{-12}$\\
\addlinespace[  0.25em]
\bottomrule
\end{tabular}

\caption{\small
The results of the experiments of Section~\ref{section:experiments:gauss}
in which the performance of the adaptive Levin method was compared with the 
performance of an adaptive Gaussian quadrature scheme.    Each row corresponds
to one of the integrals $I_5(\lambda)$, $I_6(\lambda)$, $I_7(\lambda)$ or $I_8(\lambda)$
and to one range of values of $\lambda$.
The average time taken by the adaptive Levin method and by an adaptive Gaussian quadrature
scheme, the ratio
of the average time taken by the adaptive Gaussian quadrature algorithm to the average time taken by the
adaptive Levin method, and the
maximum observed difference in the values of the integrals computed using each method  are reported.
}
\label{table:gauss:table1}
\end{table}

\end{subsection}

\begin{subsection}{Behavior in the presence of a stationary point}
\label{section:experiments:saddle}

In the experiments of this subsection, we considered the integral
\begin{equation}
\begin{aligned}
I_9(\lambda, m)  &= \int_{-1}^{1} \exp\left(i\lambda x^m\right) \frac{\cos(x)}{1+x^2}\,dx
\end{aligned}
\label{experiments:saddle:ints}
\end{equation}
in order to understand the behavior of the adaptive Levin method in the presence
of a stationary point.

In the first experiment, we sampled $l=200$ equispaced points
$x_1,\ldots,x_l$ in the interval $[1,7]$ and, for each $\lambda=10^{x_1},10^{x_2},\ldots,
10^{x_l}$ and $m=2,3,\ldots,9$, we evaluated  $I_9(\lambda,m)$ 
twice using the adaptive Levin method.  The first time the scheme was executed, the 
tolerance parameter taken to be $\epsilon = 10^{-12}$. The second time, we set $\epsilon = 10^{-7}$.
The results are shown in the first two rows of Figure~\ref{figure:saddleplot1}.
The first column gives the results for $\epsilon=10^{-7}$ and 
the second for $\epsilon=10^{-12}$.
The plots in the first row show the number of subintervals in 
the adaptively determined subdivision of $[-1,1]$ used to compute the integral
as a function of $\lambda$ for $m=2,3,4,5$, while those in the second
give the absolute error in the computed value of the integral
as a function of $\lambda$ for $m=2,3,4,5$.

\begin{figure}[b!]

\hfil
\includegraphics[width=.40\textwidth]{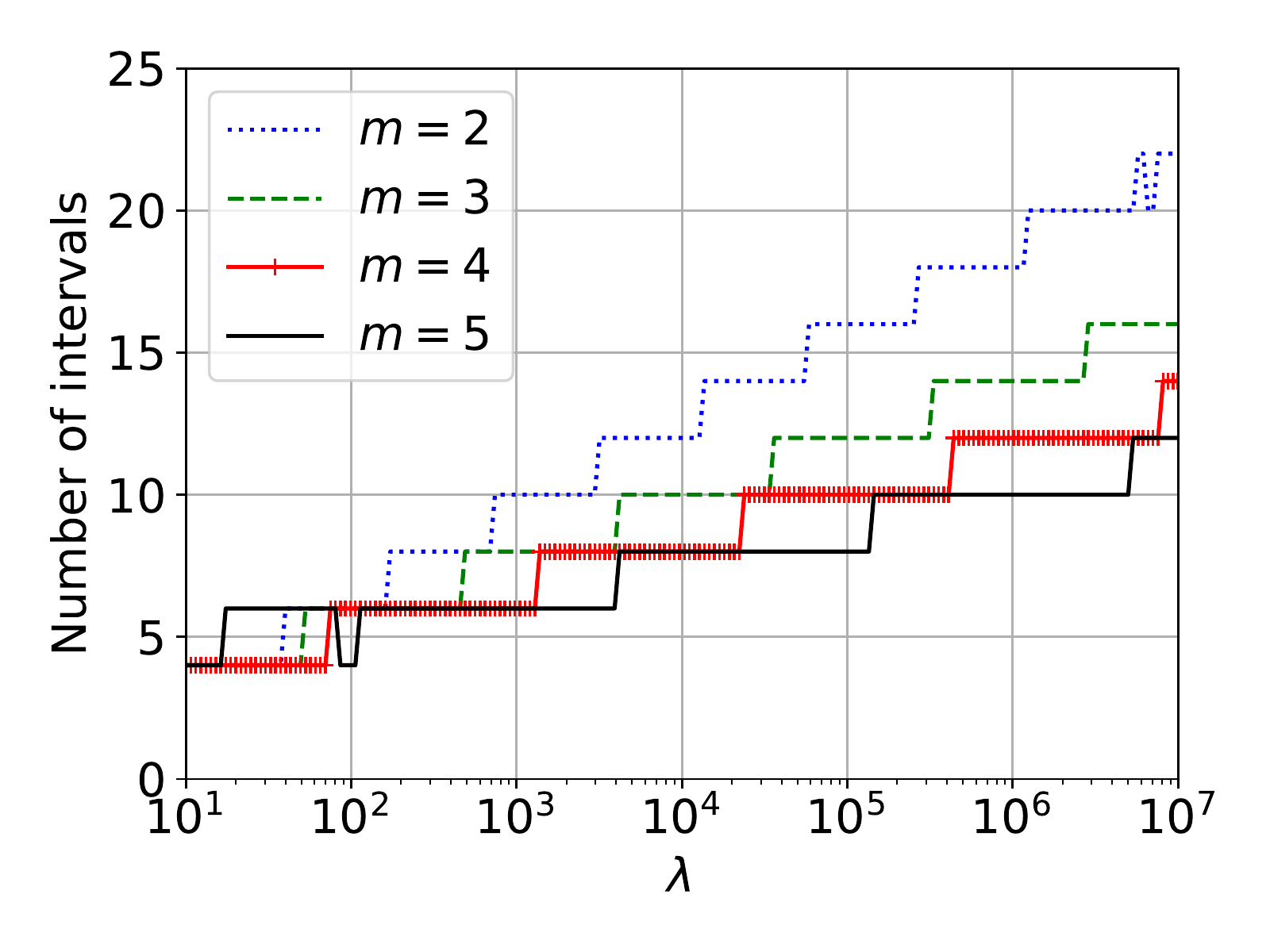}
\hfil
\includegraphics[width=.40\textwidth]{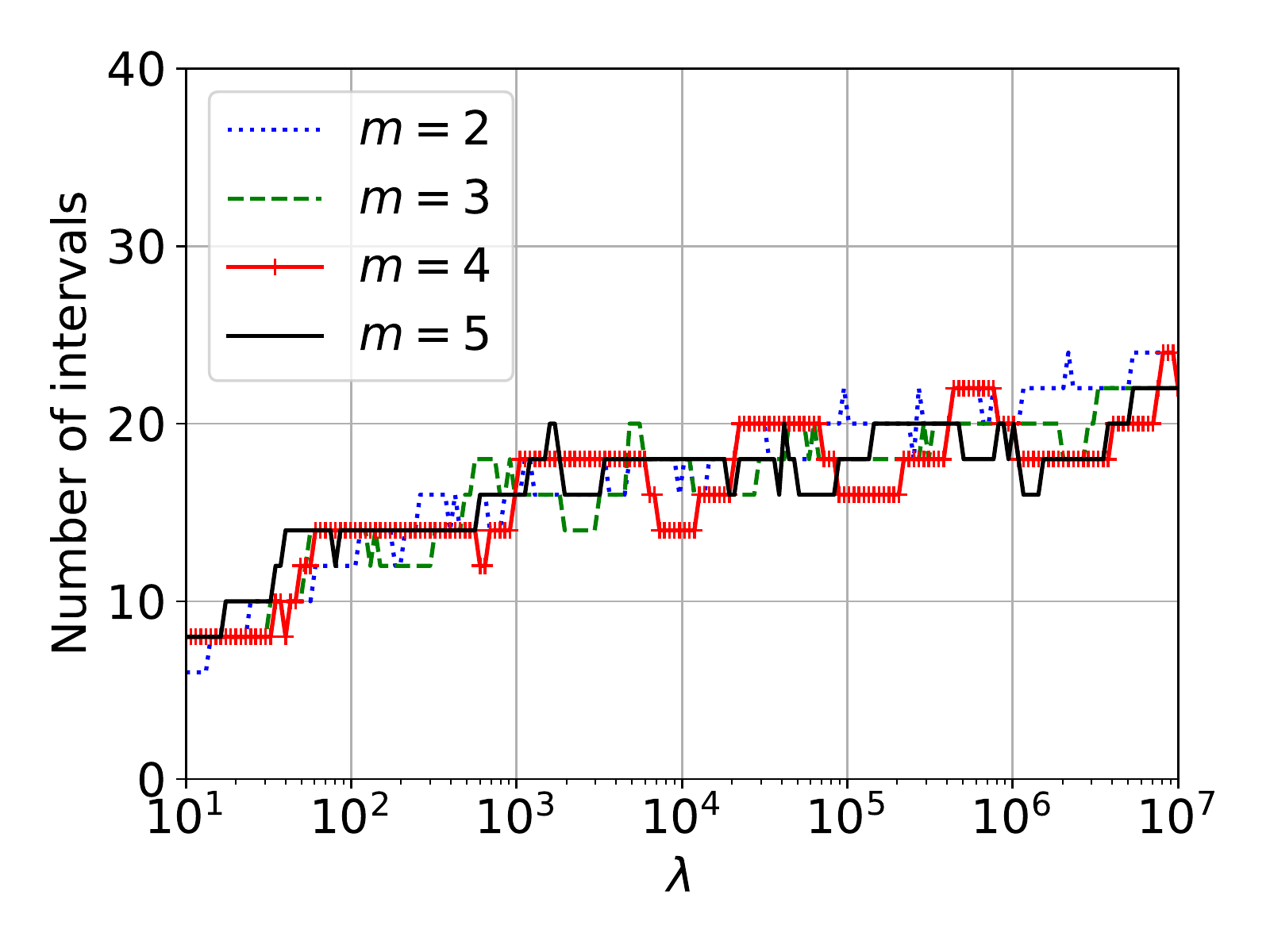}
\hfil

\hfil
\includegraphics[width=.40\textwidth]{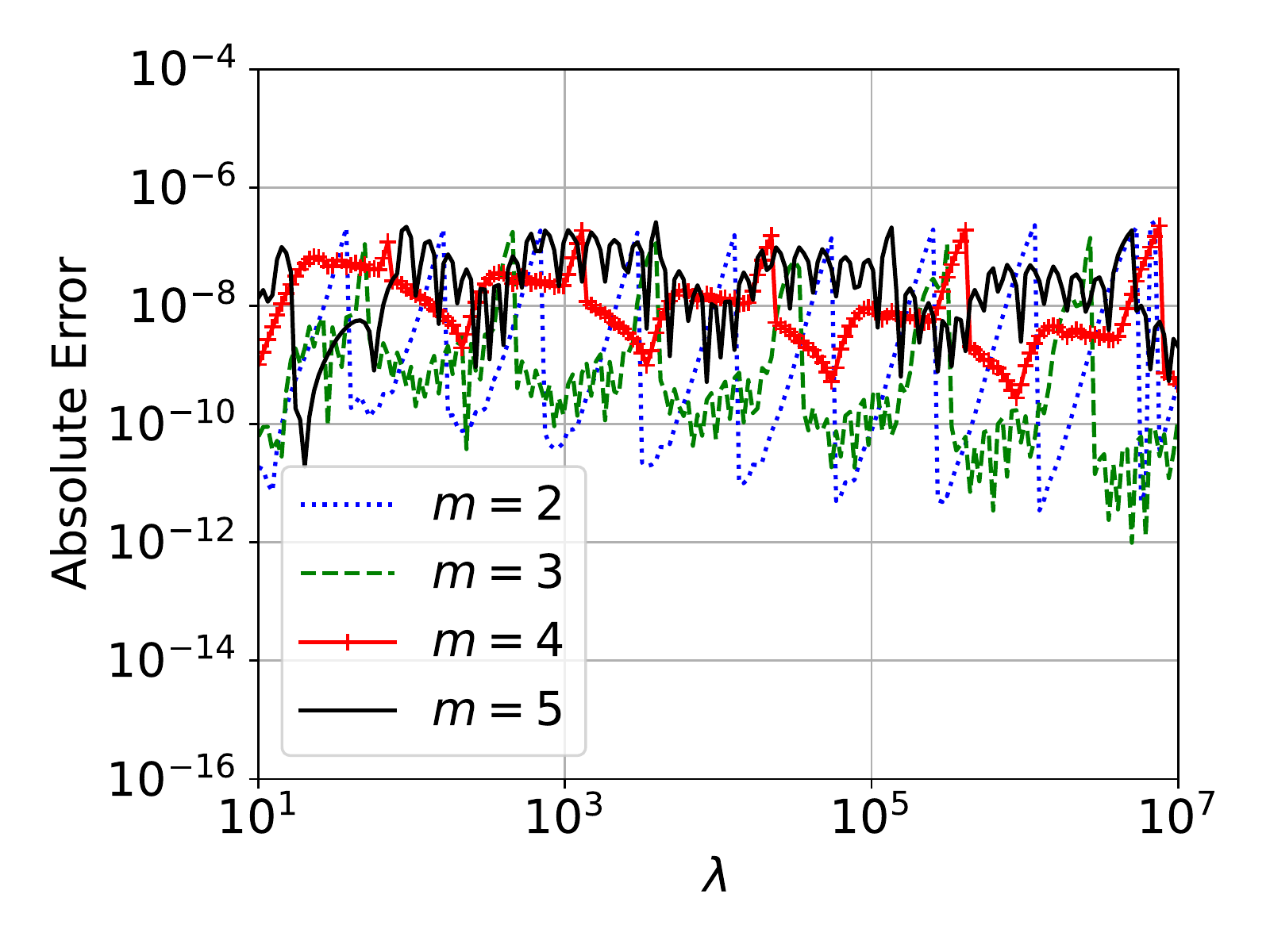}
\hfil
\includegraphics[width=.40\textwidth]{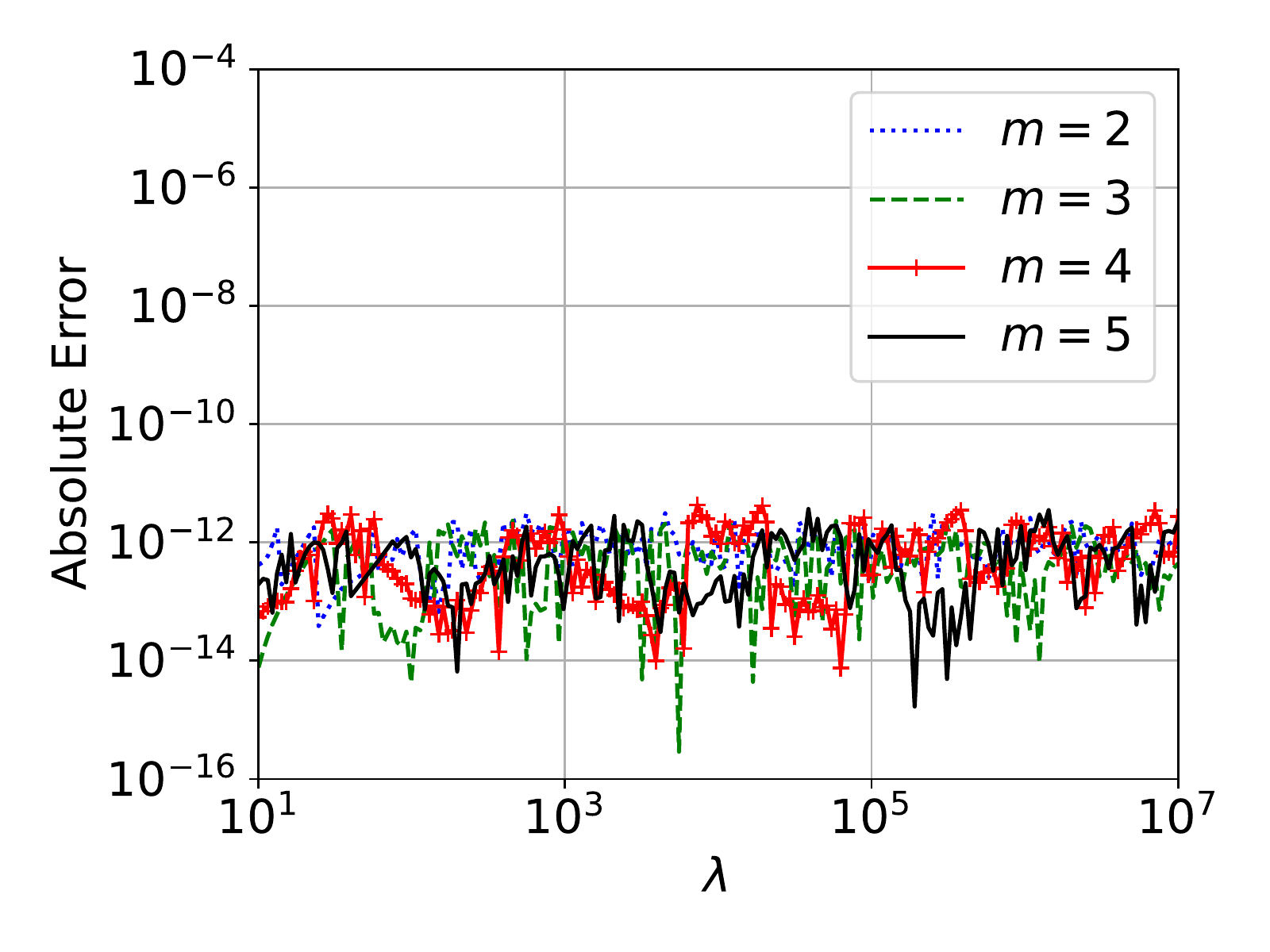}
\hfil

\hfil
\includegraphics[width=.40\textwidth]{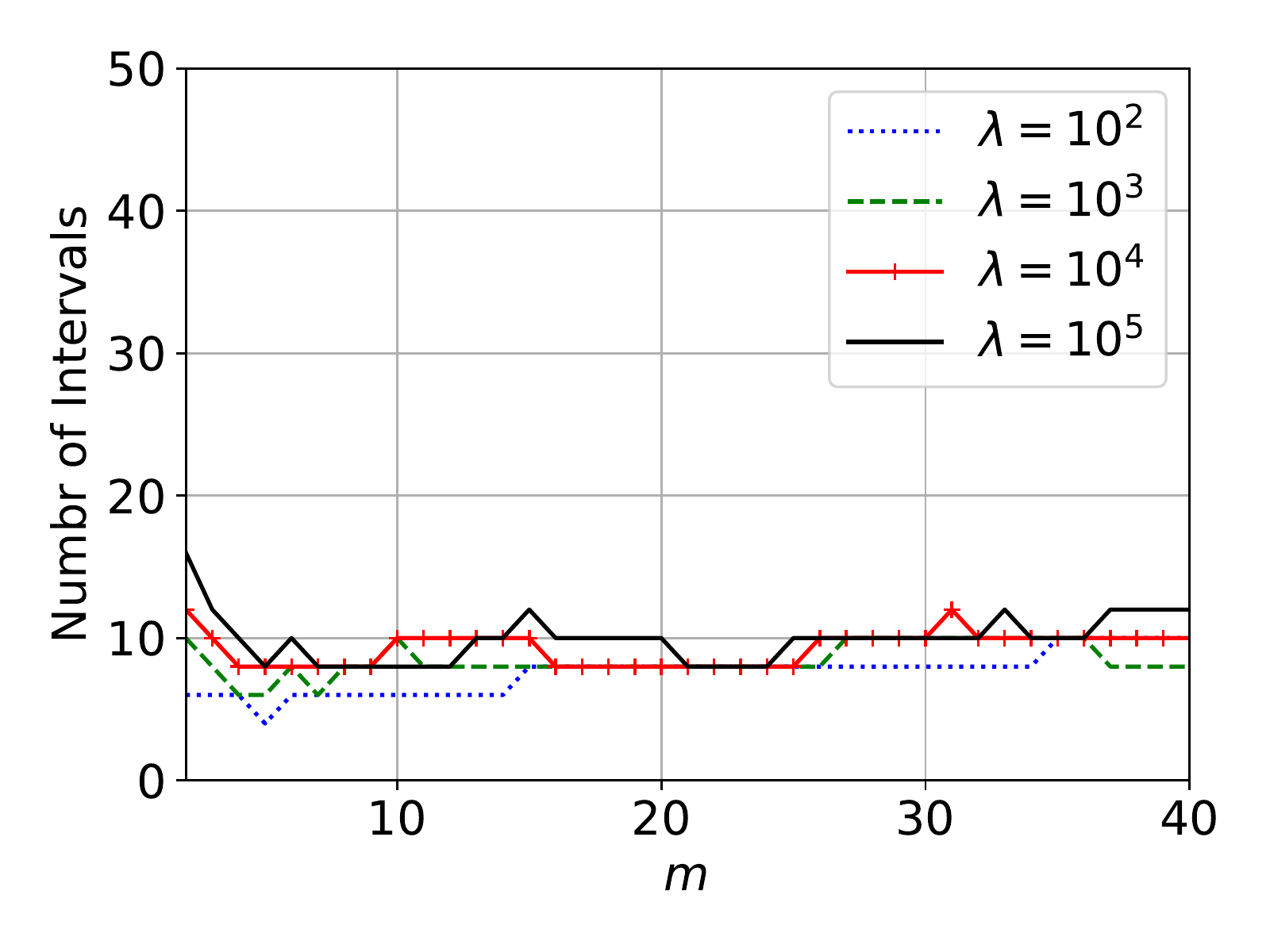}
\hfil
\includegraphics[width=.40\textwidth]{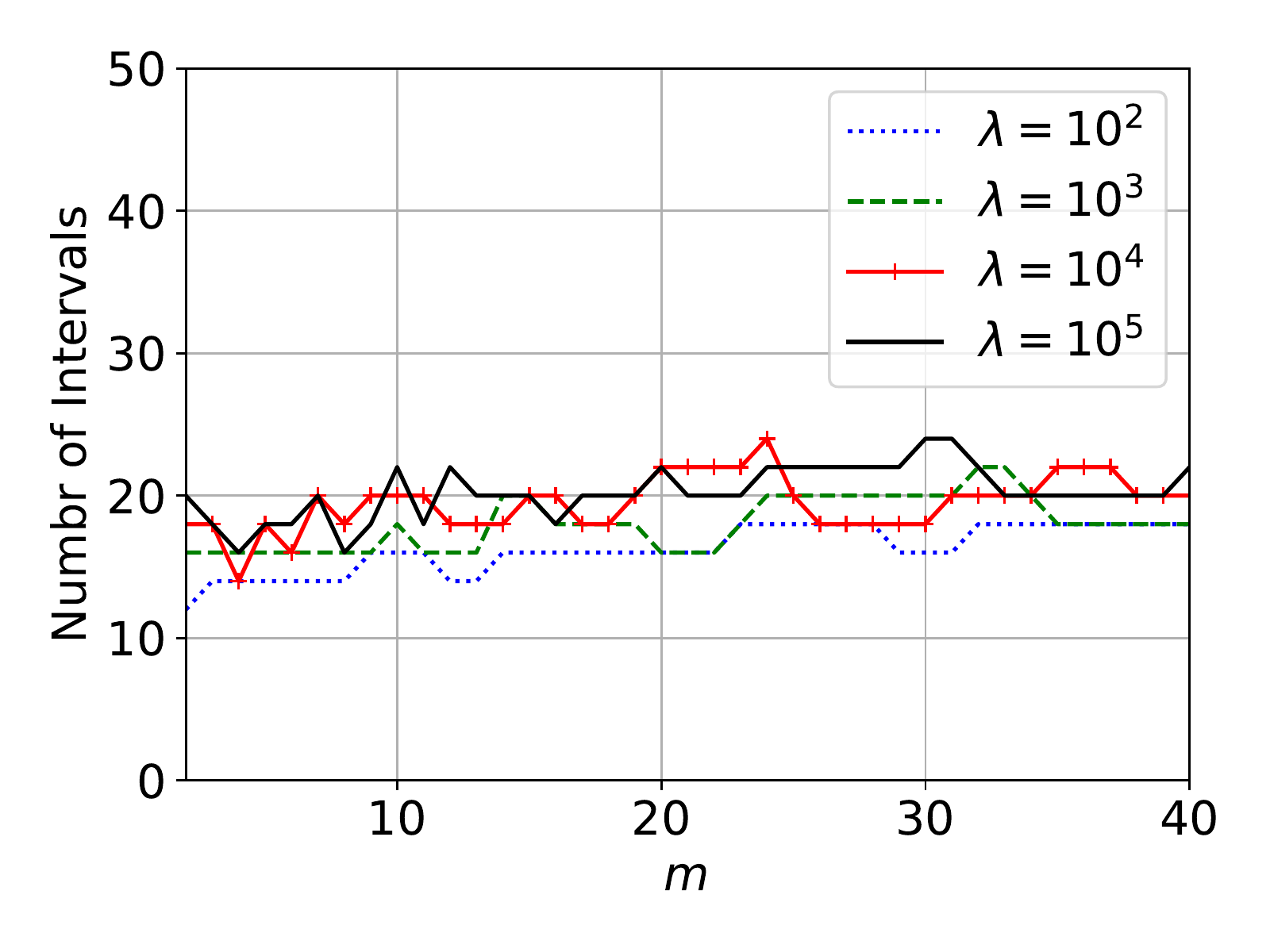}
\hfil

\caption{\small
The results of the experiments of Section~\ref{section:experiments:saddle}.
The plots in the first row give the number of subintervals in the adaptive discretization of $[-1,1]$ formed
in the course of evaluating $I_9(\lambda,m)$ as a function of $\lambda$ for $m=2,3,4,5$; those in the second 
give the absolute error in the calculated error in the 
value  of $I_9(\lambda,m)$ as a  function of $\lambda$ for $m=2,3,4,5$;
and the plots in the third row show  the number of subintervals
in the adaptive discretization of $[-1,1]$ as a function of $m$ for
$\lambda=10^2, 10^3, 10^4, 10^5$.
The plots in the column on the left concern experiments executed with
the precision parameter  $\epsilon$ taken to be  $10^{-7}$,
while those on the right correspond to $\epsilon = 10^{-12}$.
}
\label{figure:saddleplot1}
\end{figure}

In a second experiment, 
for each $m=2,3,4,\ldots,40$ and $\lambda=10^2, 10^3, 10^4, 10^5$ we evaluated $I_9(\lambda,m)$ twice using the adaptive
Levin method.  In the first run, the tolerance parameter was taken to be $\epsilon=10^{-12}$
and, in the second, it was $\epsilon=10^{-7}$.  
The results appear in the third row of Figure~\ref{figure:saddleplot1}.  Each plot there gives
the number of subintervals in the adaptive subdivision of $[-1,1]$ formed by our algorithm
as a  function of $m$ for each of the values of $\lambda$ considered.  The plot on the left
corresponds to $\epsilon=10^{-12}$ and the plot on the right to $\epsilon=10^{-7}$.

We see  that, when $\epsilon = 10^{-12}$,
the number of subintervals in the adaptive subdivision of $[-1,1]$
formed by the adaptive Levin method varies only slightly with $\lambda$ and 
is largely independent of $m$.  
However, when $\epsilon = 10^{-7}$,
the number of subintervals grows roughly logarithmically with $\lambda$ and
decreases with $m$ for small values of  $m$.
It is essentially independent of $m$ for moderate to large values of $m$.

This behavior can be understood in light of the the analysis presented in 
Section~\ref{section:leveq}.
Theorem~\ref{leveq:thm2} implies that the Levin method will yield an accurate
result on subintervals of the form $[-\delta,\delta]$ provided
$g'$ is sufficiently small there.
Theorem~\ref{leveq:thm1} indicates that the Levin method will yield an accurate
result on any subinterval of $[-1,1]$ which is bounded away from $0$
provided the ratio of the maximum to minimum value of $g'$ is small
and
\begin{equation}
h(z) = f\left(u^{-1}(z)\right) \frac{d u^{-1}}{dz}(z)
= \frac{f\left(\mbox{sign}(z)\left|z\right|^{\frac{1}{m}}\right)}
{m |z|^{1-\frac{1}{m}}}
\label{experiments:saddle:h}
\end{equation}
can be represented via a well-behaved function with a small bandlimit there.  
That $h$ can be represented by a function with small bandlimit 
is more-or-less equivalent to the requirement that $h$ be represented via a Chebyshev expansion of small
fixed order.  
We expect the adaptive Levin method to subdivide $[-1,1]$ 
until one or both of Theorems~\ref{leveq:thm1} and \ref{leveq:thm2} apply to each of the resulting subintervals.

Since $g'(x) = \lambda m x^{m-1}$, in order for $|g'(x)|$ to be bounded by a constant $C$ on $[-\delta,\delta]$, we must have
\begin{equation}
\delta <  \left(\frac{C}{m \lambda}\right)^{\frac{1}{m-1}}.
\end{equation}
So we  expect the adaptive Levin method to divide $[-1,1]$
into at least
\begin{equation}
\log(\delta) \sim \frac{\log(m \lambda)}{m-1}
\end{equation}
subintervals.   The algorithm will further divide
the subintervals contained in $[-1,-\delta]$ and  $[\delta,1]$ 
until, on each of the resulting subintervals,
the ratio of the maximum to minimum value of $g'$ is small
and $h$ is well-represented via a Chebyshev expansion of small fixed order.
 The number of subintervals required for these conditions to be met
depends only weakly on $m$ and $\lambda$.   The dependence
on $\lambda$ arises because the
greater the distance $\delta$ between these intervals and $0$,
the fewer subdivisons required, and the distance $\delta$ depends on $\lambda$.

When $\epsilon = 10^{-12}$,  the the cost of our algorithm
is dominated by the need to represent $h(z)$.  This depends
only weakly on $m$ and $\lambda$, and this behavior is reflected in 
 the results shown in Figure~\ref{figure:saddleplot1}.
They indicate that the cost of the adaptive Levin method 
grows only very slowly with $\lambda$ in this case, and that 
it is essentially independent of $m$.

When $\epsilon=10^{-7}$, the difficulty
of representing $h(x)$ is lower,  and the running time  our algorithm 
is dominated by the number of subdivisions needed to ensure
$\delta$ is sufficiently small, at least when $m$ is small.
 This is on the order of $\log(m \lambda)/(m-1)$,
and we see this behavior in the plot appearing on the left-hand side of the first row
of Figure~\ref{figure:saddleplot1}.
As $m$ increases, though, the cost of representing $h(z)$
begins to dominate the running time of the algorithm and it ceases
to depends strongly on $\lambda$.
This is reflected in plot appearing on the lower-left corner of 
 Figure~\ref{figure:saddleplot1}.

The behavior of the adaptive Levin method in the  presence of stationary points
is somewhat complicated, but from this analysis and the experiments
described above, it is safe to conclude that at worst
the algorithm grows logarithmically with the frequency of $g'$.

\end{subsection}

\begin{subsection}{Certain integrals involving the Bessel functions}
\label{section:experiments:bes}
In the experiments described in this subsection, we used the adaptive Levin method to
evaluate the integrals
\begin{equation}
\begin{aligned}
I_{10}(\nu) &= \int_0^\infty\frac{ J_\nu(x) }{\sqrt{x}} \,dx
=
\frac{\Gamma\left(\frac{\nu}{2}+\frac{1}{4}\right)}{\sqrt{2}\, \Gamma\left(\frac{\nu}{2}+\frac{3}{4}\right)}
,\\
I_{11}(\lambda) &= \int_0^{\infty} Y_{\frac{1}{2}}(\lambda x) \ \exp(-x)\,dx = 
-\frac{\sqrt{\sqrt{\lambda ^2+1}+1}}{\sqrt{\lambda ^3+\lambda}}
,\\
I_{12}(\lambda) &=  \int_0^{\frac{\pi}{2}} 
J_{2\lambda}\left(2\lambda \cos(x))\right)\, dx = \frac{\pi}{2} J_\lambda^2\left(\lambda\right)
\ \ \mbox{and}\\
I_{13}(\lambda) &= 
\int_0^\infty J_0(\lambda x) J_{\frac{1}{2}}(\lambda x) \exp(-x)\sqrt{x}\, dx = 
\sqrt{
\frac
{-1 + \sqrt{1+4\lambda^2}}{\pi\lambda + 4\pi\lambda^3}
}.
\end{aligned}
\label{experiments:bessel:int}
\end{equation}
Each of the above formulas can  be found either in  \cite{Gradshteyn} or \cite{TABLESII}.

In the first experiment, we sampled  $l=200$ equispaced points $x_1,x_2,\ldots,x_l$ in the interval $[1,7]$.  
Then, for each $\nu = 10^{x_1},\ldots,10^{x_l}$, we constructed a phase function for the normal form 
\begin{equation}
y''(x) + \left(1+\frac{\frac{1}{4}-\nu^2}{x^2}\right) y(x) = 0
\label{experiments:1:besseleq}
\end{equation}
of Bessel's differential equation using the algorithm of \cite{BremerPhase2}.  
This gave us the following representations of the Bessel function of the first and second kinds of order 
$\lambda$:
\begin{equation}
\begin{aligned}
J_\nu(x) = \sqrt{\frac{2}{\pi x}} \frac{\sin\left(\psi^{\mbox{\tiny bes}}_\nu(x)\right)}{\sqrt{\frac{d}{dx}\psi^{\mbox{\tiny bes}}_\nu(x)}}
\ \ \ \mbox{and} \ \ \
Y_\nu(x) = \sqrt{\frac{2}{\pi x}} \frac{\cos\left(\psi^{\mbox{\tiny bes}}_\nu(x)\right)}{\sqrt{\frac{d}{dx}\psi^{\mbox{\tiny bes}}_\nu(x)}}.
\end{aligned}
\label{experiments:bessel:phase}
\end{equation}
The first of these representations was used, together with the adaptive Levin method, to evaluate $I_{10}(\nu)$.  
More explicitly, we took the input functions for the adaptive Levin method to be
\begin{equation}
\begin{aligned}
g(x) =  \psi^{\mbox{\tiny bes}}_\nu(x)\ \ \ \ \ \mbox{and} \ \ \ \mbox{and}\ \ \ 
f(x) =  \frac{1}{x} \sqrt{\frac{2}{\pi\, \frac{d}{dx}\psi^{\mbox{\tiny bes}}_\nu(x)}}.
\end{aligned}
\end{equation}
The results are given in the first row of Figure~\ref{figure:besplots1}.  Note that for this experiment,
we report the time taken by the adaptive Levin method and the time required to construct the phase
function separately.

\begin{figure}[t!]

\hfil
\includegraphics[width=.40\textwidth]{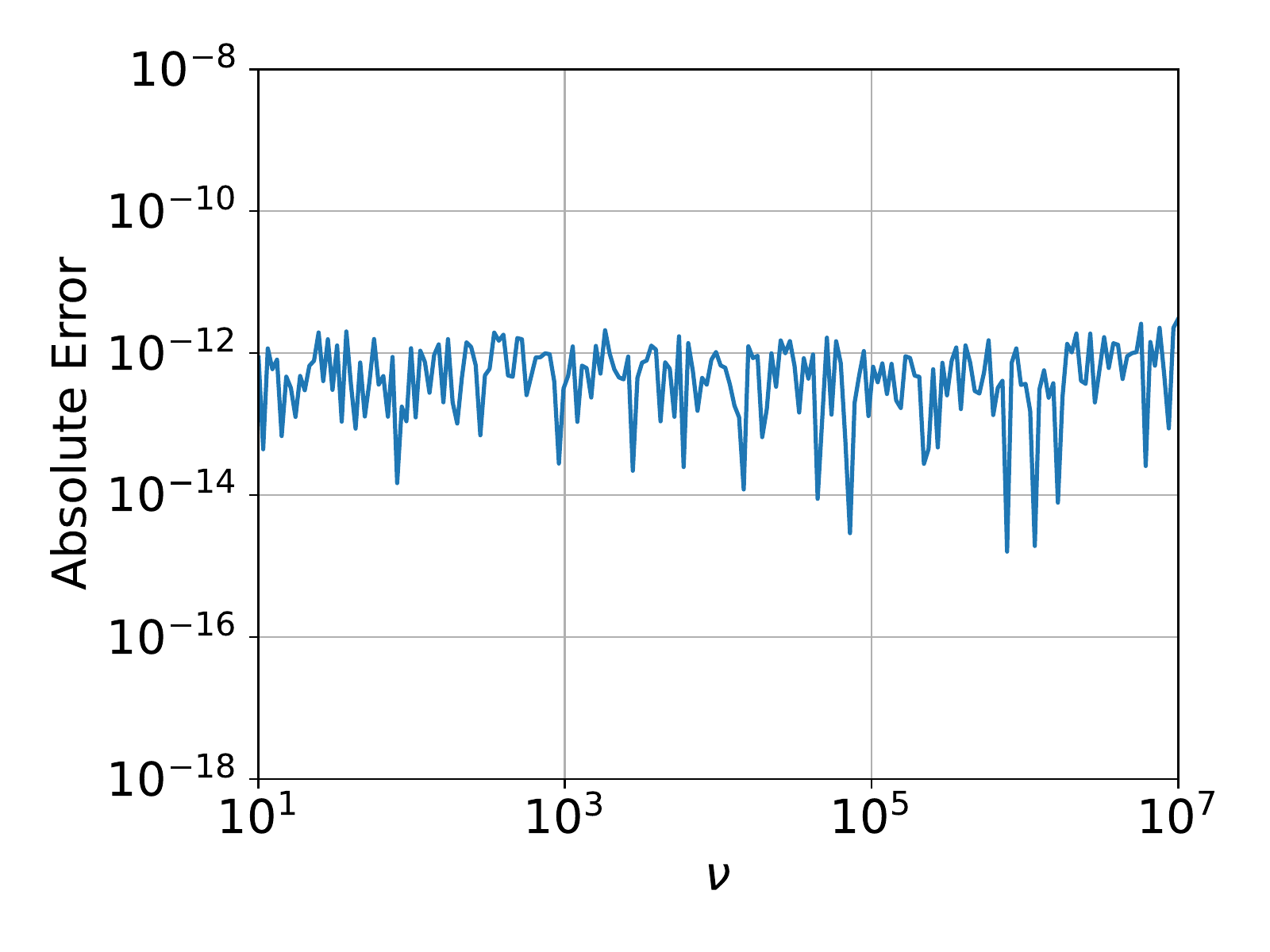}
\hfil
\includegraphics[width=.40\textwidth]{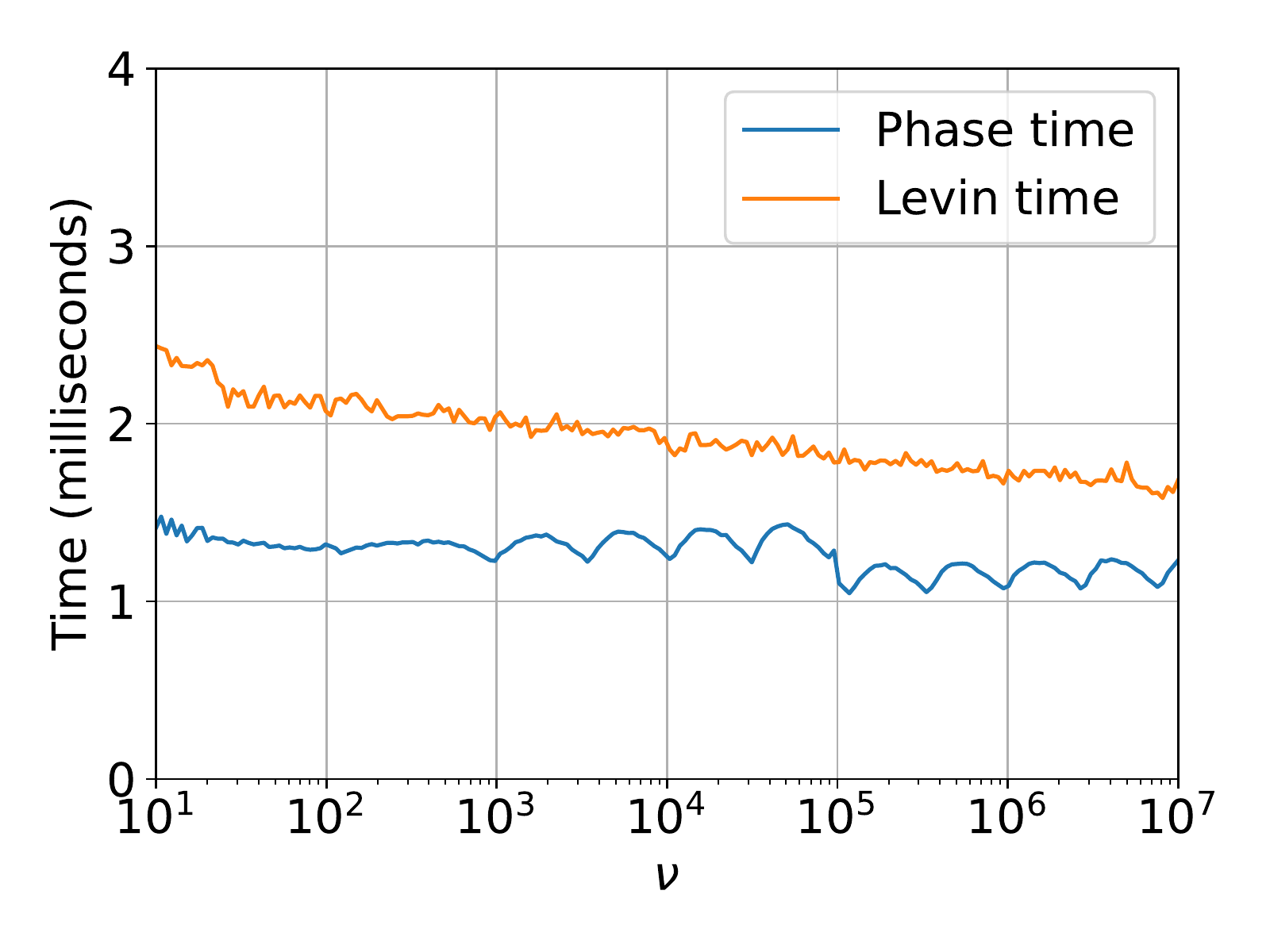}
\hfil

\hfil
\includegraphics[width=.40\textwidth]{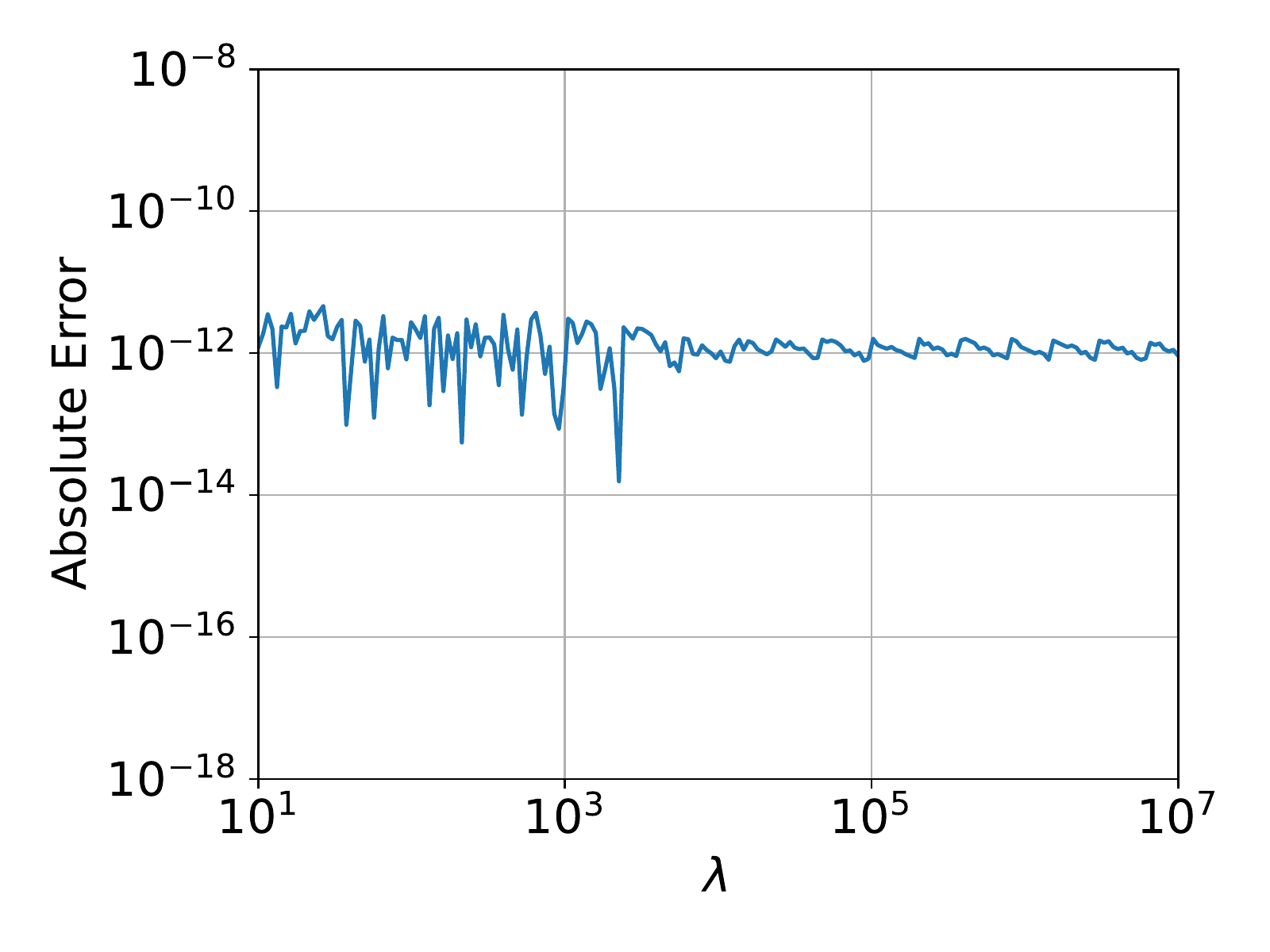}
\hfil
\includegraphics[width=.40\textwidth]{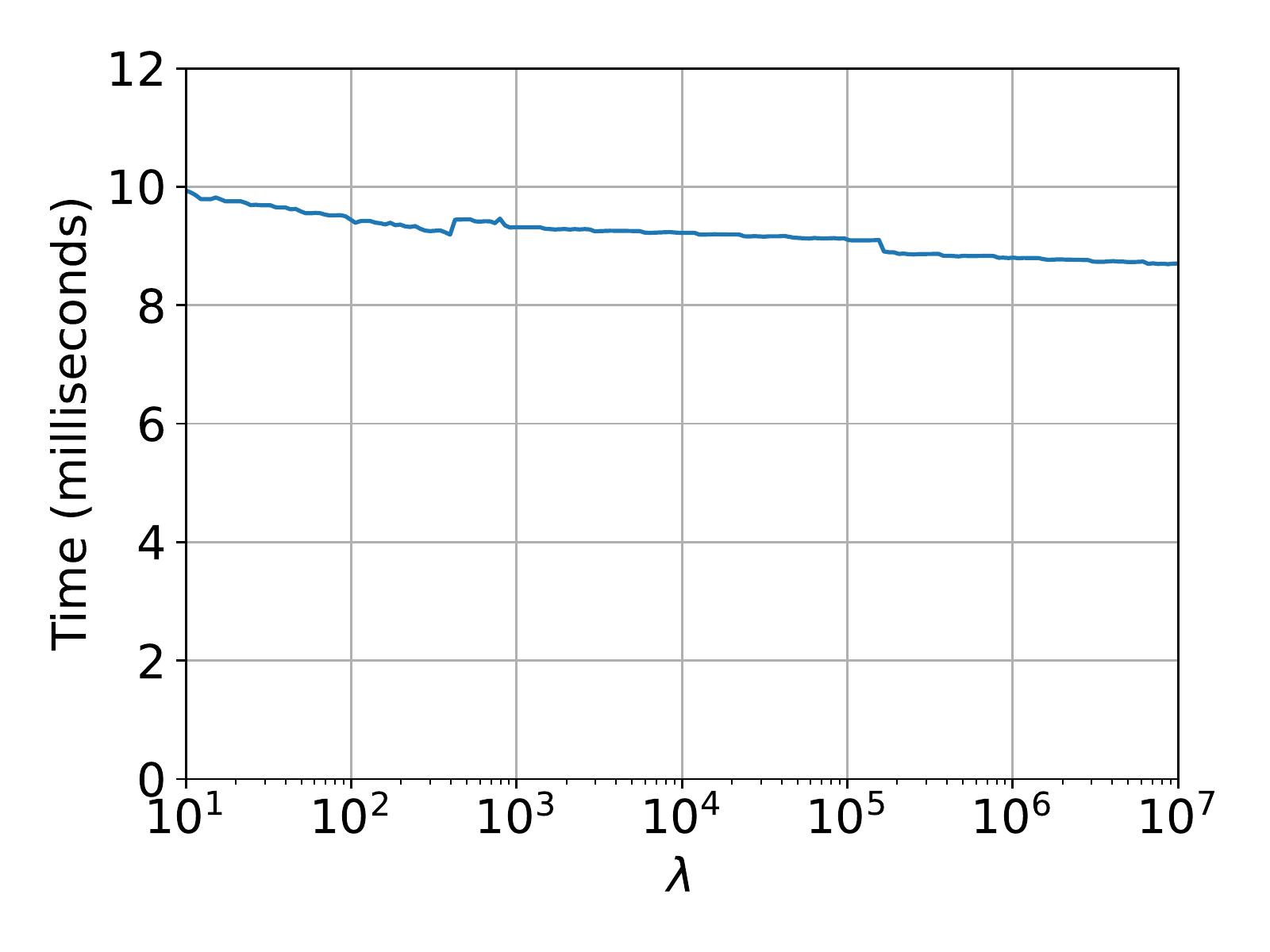}
\hfil

\caption{\small
The results of the first two experiments of Section~\ref{section:experiments:bes}.
In the first row, the plot on the left gives the absolute error in the calculation
of the integral $I_{10}(\nu)$ as a function $\nu$  and the plot on the right shows
the time take by the adaptive Levin method and the time taken to constuct
the phase function as functions of $\nu$.
The plot at bottom left gives the absolute error in the obtained
value of $I_{11}(\lambda)$ as a function of $\lambda$, and 
the bottom-right plot gives the time required to compute $I_{11}(\lambda)$,
including the time required to construct any necessary phase functions, as a function of $\lambda$.
}
\label{figure:besplots1}
\end{figure}

We began our second experiment by sampling $l=200$ equispaced points $x_1,\ldots,x_l$ in the interval $[1,7]$.
Then, for each $\lambda = 10^{x_1},10^{x_2},\ldots,10^{x_l}$, we constructed
the phase function  $\psi^{\mbox{\tiny bes}}_{1/2}$ and
executed the adaptive Levin method with the input functions taken to be 
\begin{equation}
\begin{aligned}
g(x) =  \psi^{\mbox{\tiny bes}}_{1/2}(\lambda x)\ \ \ \mbox{and}\ \ \ 
f(x) =  \sqrt{\frac{2}{\pi}} \exp(-x)
\sqrt{\frac{1}{\lambda x\, \frac{d}{dx}\psi^{\mbox{\tiny bes}}_{1/2}(\lambda x)}}
\end{aligned}
\end{equation}
in order to evaluate $I_{11}(\lambda)$.  The results are shown in the second
row of Figure~\ref{figure:besplots1}.  The reported times include the cost of constructing the phase
function as well as the time required by the adaptive Levin method.

We began our third experiment regarding Bessel functions by sampling $l=200$ equispaced points $x_1,\ldots,x_l$ in the interval $[1,7]$.
Then, for each $\lambda = 10^{x_1},10^{x_2},\ldots,10^{x_l}$, we constructed
a phase function  $\psi^{\mbox{\tiny bes}}_{2\lambda}$ representing
the Bessel functions of order $2\lambda$ and used the  adaptive Levin method to evaluate $I_{12}(\lambda)$;
the input functions were taken to be 
\begin{equation*}
g(x) =  \psi^{\mbox{\tiny bes}}_{2\lambda}(2 \lambda \cos(x))\ \ \ \mbox{and}\ \ \ 
f(x) = \sqrt{\frac{2}{\pi}} \sqrt{\frac{1}{2\lambda\cos(x) \, \frac{d}{dx}\psi^{\mbox{\tiny bes}}_{2\lambda}(2 \lambda \cos(x))}}.
\end{equation*}
%
The results are shown in the first
row of Figure~\ref{figure:besplots2}.  Once  again, the reported times include the cost of constructing 
the necessary phase functions as well as the time required by the adaptive Levin method.

\begin{figure}[t!]

\hfil
\includegraphics[width=.40\textwidth]{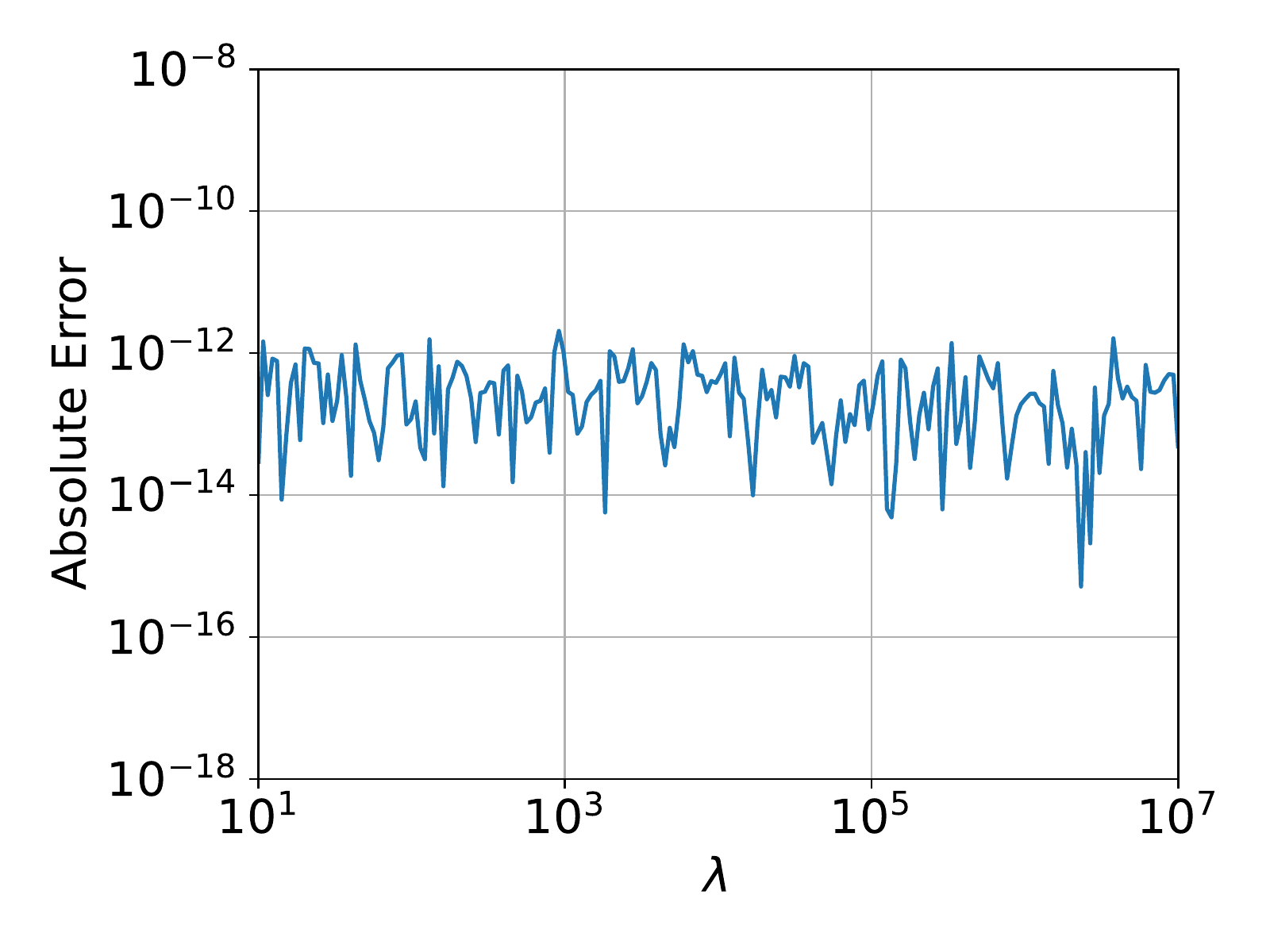}
\hfil
\includegraphics[width=.40\textwidth]{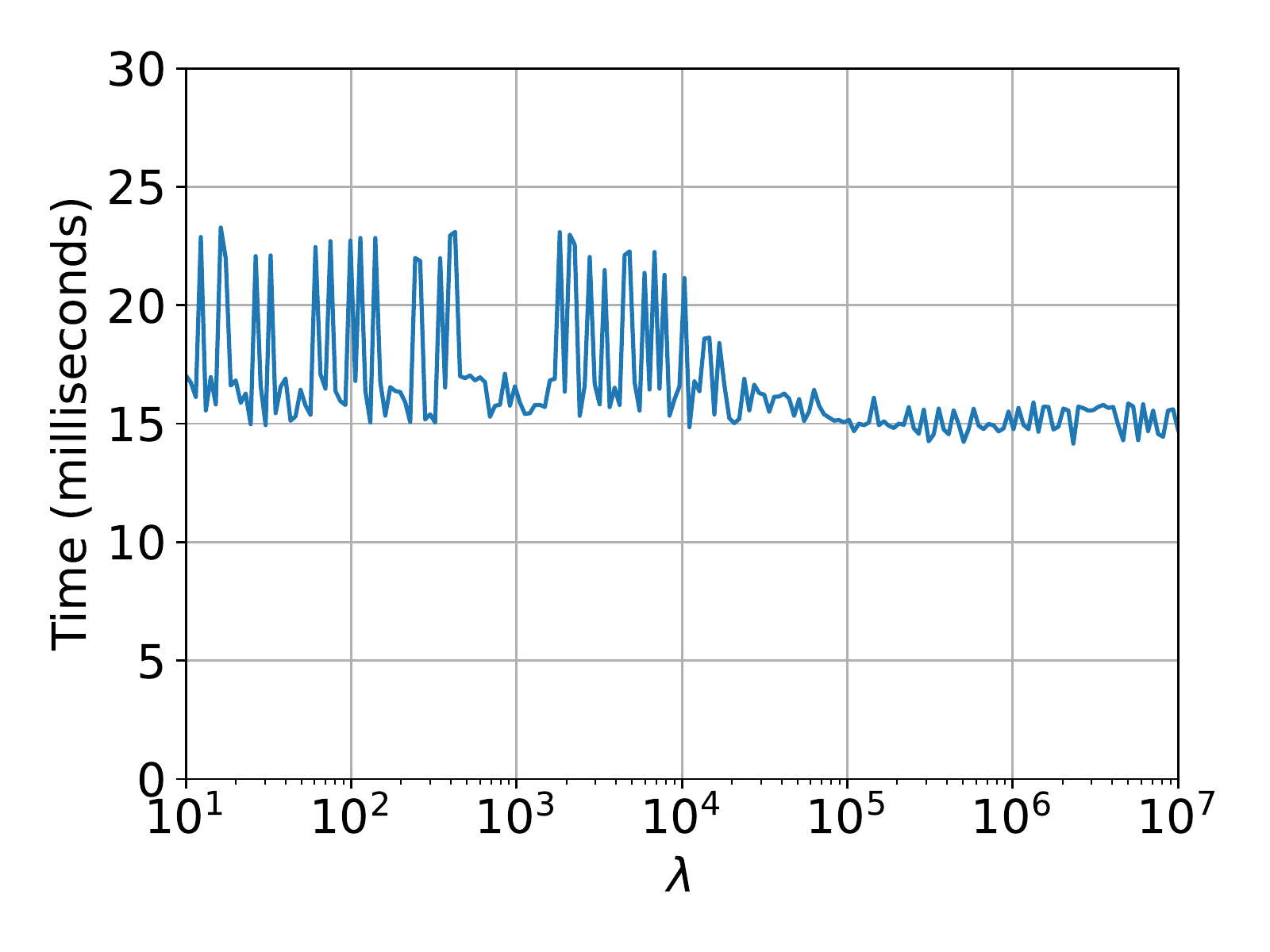}
\hfil

\hfil
\includegraphics[width=.40\textwidth]{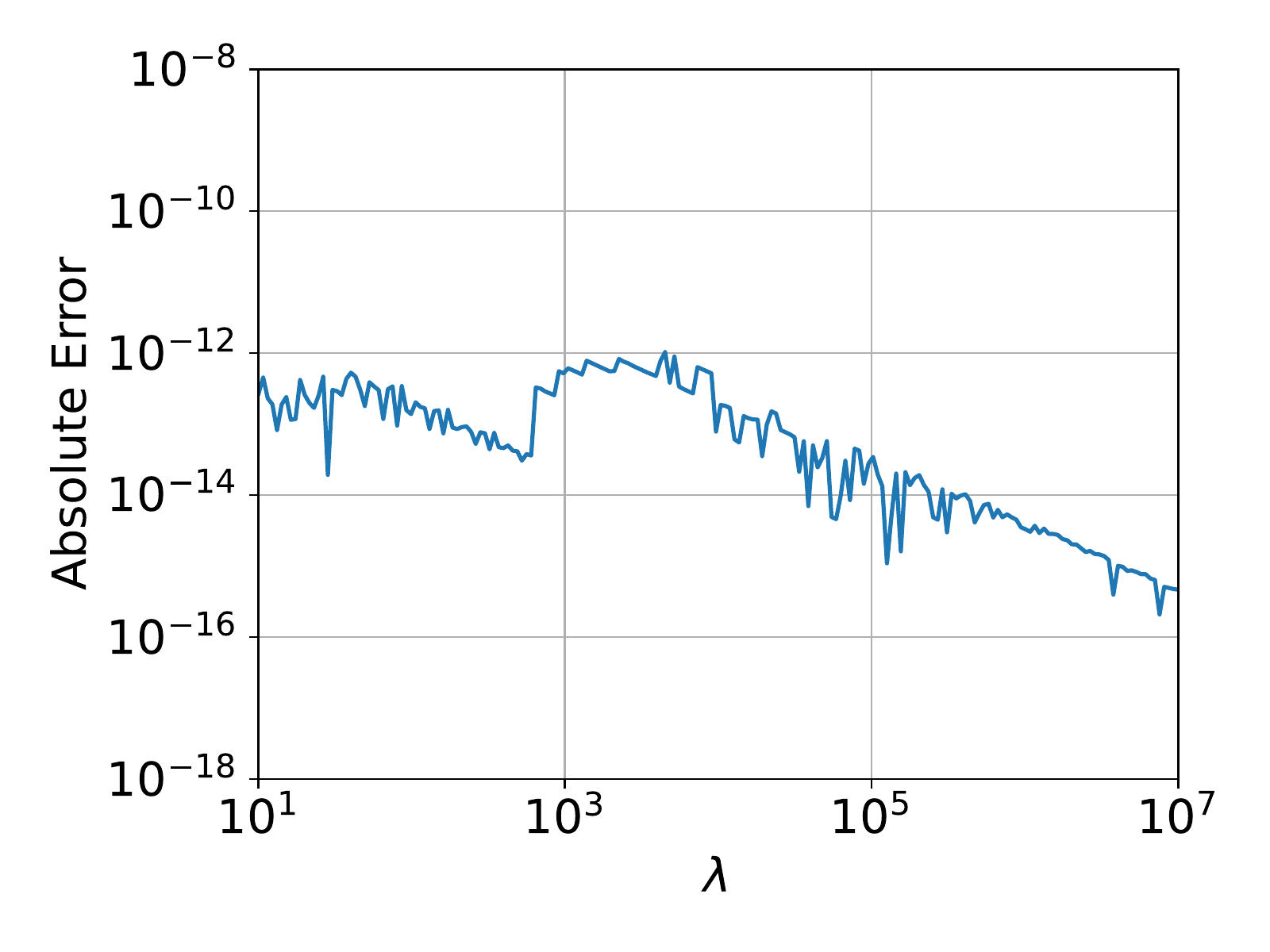}
\hfil
\includegraphics[width=.40\textwidth]{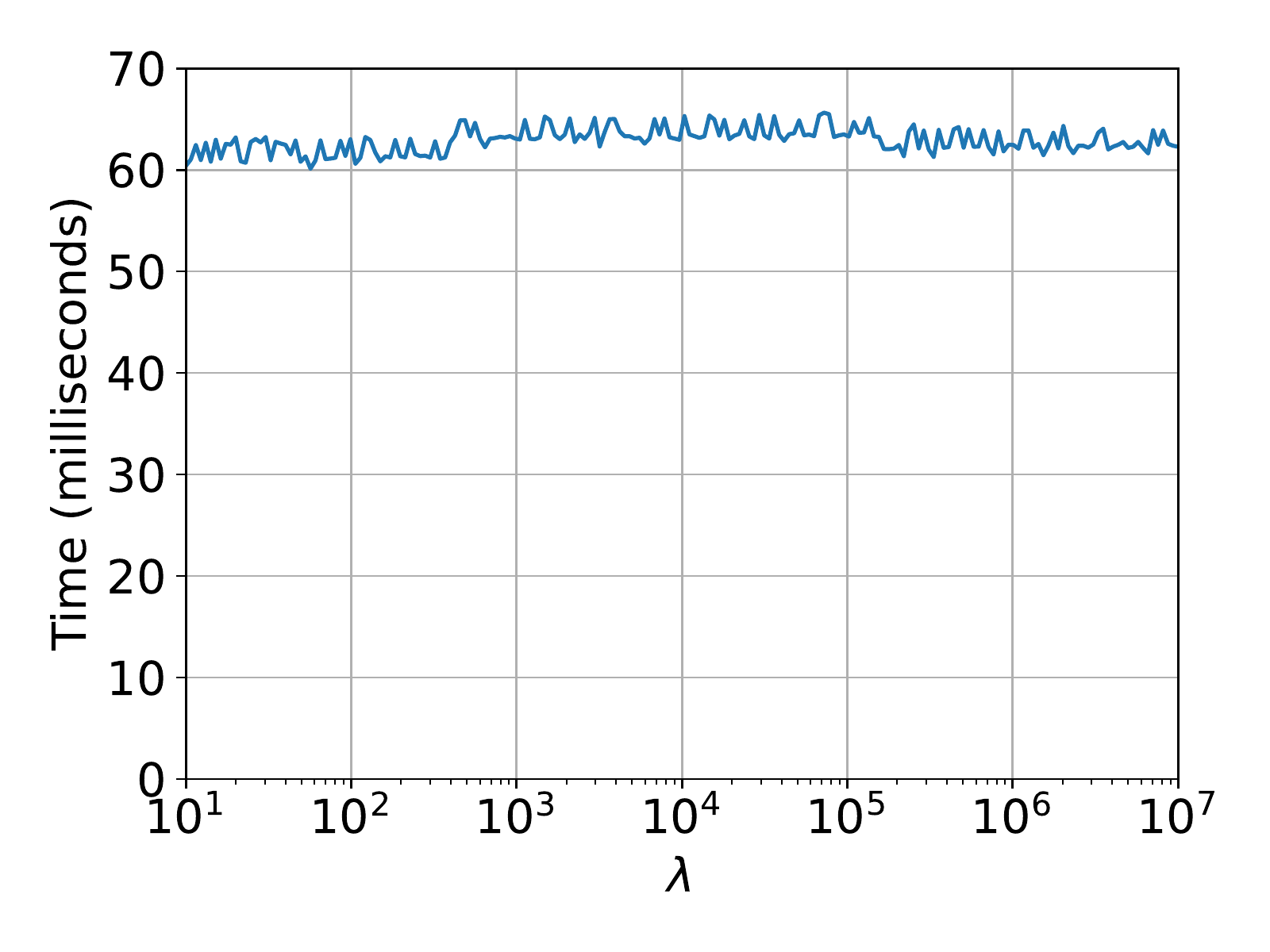}
\hfil

\caption{\small
The results of the last two experiments of Section~\ref{section:experiments:bes}.
The first row pertains to   $I_{12}(\lambda)$ while the second concerns $I_{13}(\lambda)$.
In each row, the plot on the left gives the absolute error in the calculation
of the integral as a function of $\lambda$ and the plot on the right shows
the total time required to compute the integral via the adaptive Levin method
and to construct any necessary phase functions as a function of $\lambda$.
}
\label{figure:besplots2}
\end{figure}

In our fourth experiment concerning Bessel functions, we first sampled $l=200$ equispaced points $x_1,\ldots,x_l$ in the interval $[1,7]$.
Then, for each $\lambda = 10^{x_1},10^{x_2},\ldots,10^{x_l}$, we constructed the phase functions
$\psi^{\mbox{\tiny bes}}_{0}$ and  $\psi^{\mbox{\tiny bes}}_{1/2}$, one representing the Bessel functions
of order $0$ and the other the Bessel functions of order $1/2$.  Since
\begin{equation}
I_{13}(\lambda) = 
\frac{2}{\pi } \int_0^\infty 
\frac{\sin\left(\psi^{\mbox{\tiny bes}}_{0}\left(\lambda x\right)\right)}
{\sqrt{\lambda x\,  \frac{d}{dx} \psi^{\mbox{\tiny bes}}_{0}\left(\lambda x\right)}}
\frac{\sin\left(\psi^{\mbox{\tiny bes}}_{1/2}\left(\lambda x\right)\right)}
{\sqrt{\lambda x\,  \frac{d}{dx} \psi^{\mbox{\tiny bes}}_{1/2}\left(\lambda x\right)}}
  \exp(-x)\sqrt{x}\, dx
\end{equation}
and
\begin{equation}
\sin(x)\sin(y) = \frac{\cos(x-y)-\cos(x+y)}{2},
\end{equation}
we were able to compute $I_{13}$ via  the formula
\begin{equation}
I_{13} (\lambda)= \frac{2}{\lambda \pi } \frac{I_{13a} (\lambda)-I_{13b} (\lambda)}{2},
\end{equation}
where
\begin{equation}
I_{13a} (\lambda) =   \int_0^\infty\frac{\cos\left(\psi^{\mbox{\tiny bes}}_{1/2}(x) -\psi^{\mbox{\tiny bes}}_{0}(x) \right)}
{\sqrt{x\  \frac{d}{dx} \psi^{\mbox{\tiny bes}}_{0}(x)\, \frac{d}{dx} \psi^{\mbox{\tiny bes}}_{1/2}(x)}}
 \exp(-x)\, dx
\end{equation}
and
\begin{equation}
I_{13b} (\lambda) =  \int_0^\infty \frac{\cos\left(\psi^{\mbox{\tiny bes}}_{0}(x) + \psi^{\mbox{\tiny bes}}_{1/2}(x) \right)}
{\sqrt{x\  \frac{d}{dx} \psi^{\mbox{\tiny bes}}_{0}(x)\, \frac{d}{dx} \psi^{\mbox{\tiny bes}}_{1/2}(x)}}
 \exp(-x)\, dx.
\end{equation}
The integrals $I_{13a}$ and $I_{13b}$ were, of course, evaluated using the adaptive Levin method.
 The results are shown in the second
row of Figure~\ref{figure:besplots2}.  The reported times include the cost of constructing the two necessary
phase functions as well as the time required by the adaptive Levin method.


\end{subsection}

\begin{subsection}{Certain integrals involving the associated Legendre functions}
\label{section:experiments:alf}

In the experiments described in this subsection, we used the adaptive Levin method to
evaluate the integrals
\begin{equation*}
\begin{aligned}
I_{15}(\nu,\mu) &= \int_0^1 \widetilde{P}_{\nu}^\mu (x) (1-x^2)^{\frac{\mu}{2}}\,dx\\
&= 
\sqrt{\left(\nu+\frac{1}{2}\right)\frac{\Gamma(\nu+\mu+1)}{\Gamma(\nu-\mu+1)}}
\frac{(-1)\mu 2^{-\mu-1} \sqrt{\pi}} {\Gamma\left(1+\frac{\mu}{2}-\frac{\nu}{2}\right)
\Gamma\left(\frac{3}{2}+\frac{\mu}{2}+\frac{\nu}{2}\right)},
\\
I_{16}(\lambda) &= \int_0^1 \widetilde{P}^{\frac{1}{2}}_\lambda(x) \widetilde{Q}^{\frac{1}{2}}_\lambda(x)\, dx \ \ \ \mbox{and}
\\
I_{17}(\lambda) &= \int_0^{\frac{\pi}{2}} \widetilde{P}^{1}_\lambda(\cos(x))\, dx.
\end{aligned}
\end{equation*}
Here, we use $\widetilde{P}_\nu^\mu$ and $\widetilde{Q}_\nu^\mu$  to denote normalized version of the Ferrer's functions
of the first and second kinds of degree $\nu$ and order $\mu$.  
For $\nu \geq \mu$, the usual Ferrer's function
${P}_\nu^\mu$ is the the unique solution of the associated Legendre differential equation
\begin{equation}
(1-x^2) y''(x) - 2x y'(x) + \left(\nu(\nu+1) - \frac{\mu^2}{1-x^2}\right)y(x) = 0
\label{experiments:alf:ode}
\end{equation}
which is regular at the singular point $x=1$ and such that $P_\nu^\mu(1)=1$. 
Because, for most values of the parameters $\nu$ and $\mu$,  $P_\nu^\mu$  is exponentially decaying on some part of the interval $[0,1]$,   
it can take on extremely large values.  Accordingly, we prefer to work with the normalized Ferrer's function
\begin{equation}
\begin{aligned}
\widetilde{P}_\nu^\mu(x) &= \sqrt{\left(\nu+\frac{1}{2}\right)\frac{\Gamma\left(\nu+\mu+1\right)}{\Gamma\left(\nu-\mu+1\right)}} P_\nu^\mu(x)
\end{aligned}
\end{equation}
whose $L^2(-1,1)$ is $1$ when $n \geq m$ are integers.  The Ferrer's function $Q_\nu^\mu$ of the second kind 
is (essentially) $\pi/2$ times  the Hilbert transform of $P_\nu^\mu$ and we define its normalized version 
via
\begin{equation}
\begin{aligned}
\widetilde{Q}_\nu^\mu(x) &= \frac{2}{\pi} \sqrt{\left(\nu+\frac{1}{2}\right)\frac{\Gamma\left(\nu+\mu+1\right)}{\Gamma\left(\nu-\mu+1\right)}} Q_\nu^\mu(x).
\end{aligned}
\end{equation}
We refer the reader to Section~5.15 of \cite{Olver} for a thorough discussion of of the Ferrer's functions.

\begin{figure}[t!]

\hfil
\includegraphics[width=.40\textwidth]{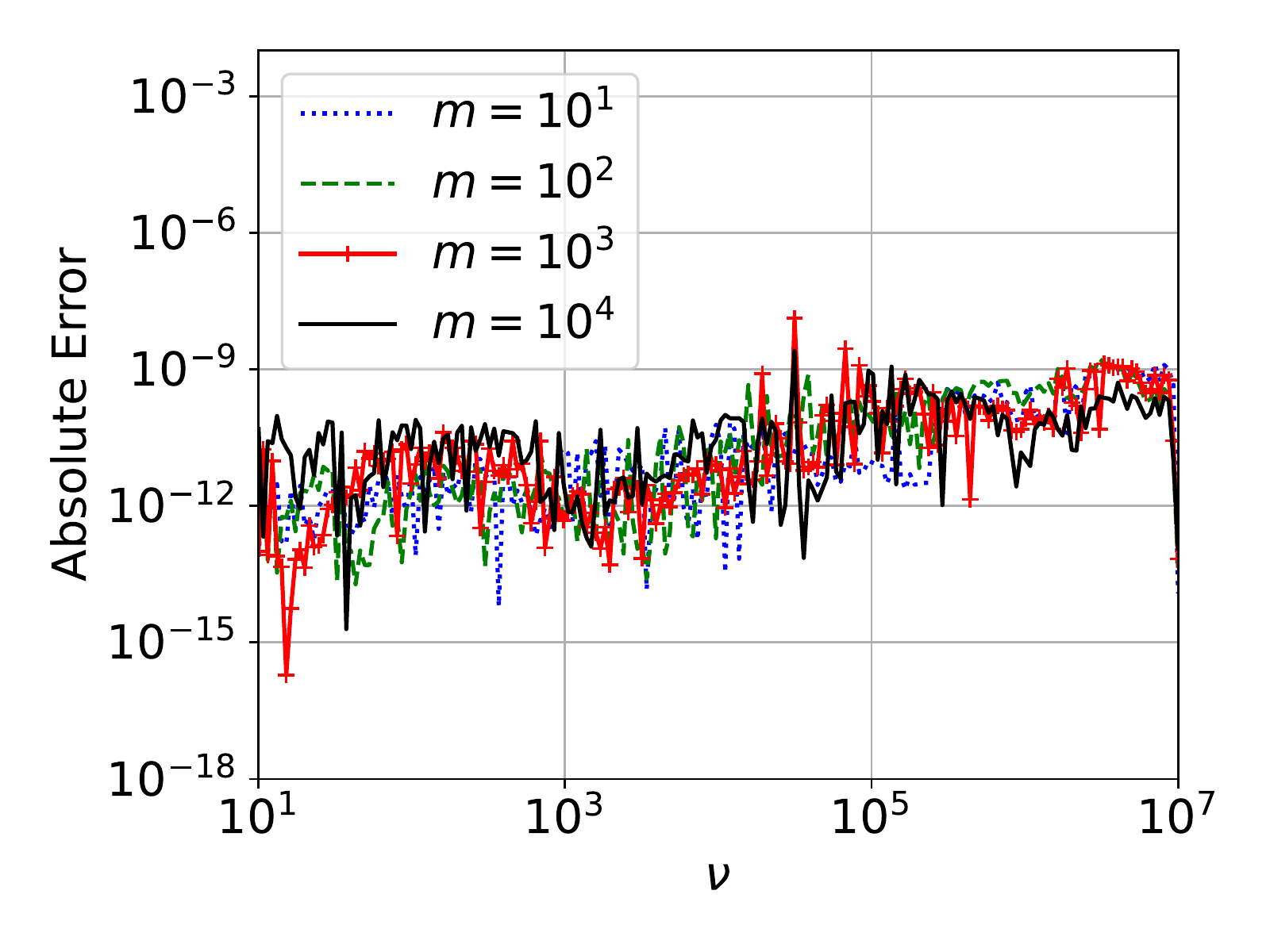}
\hfil
\includegraphics[width=.40\textwidth]{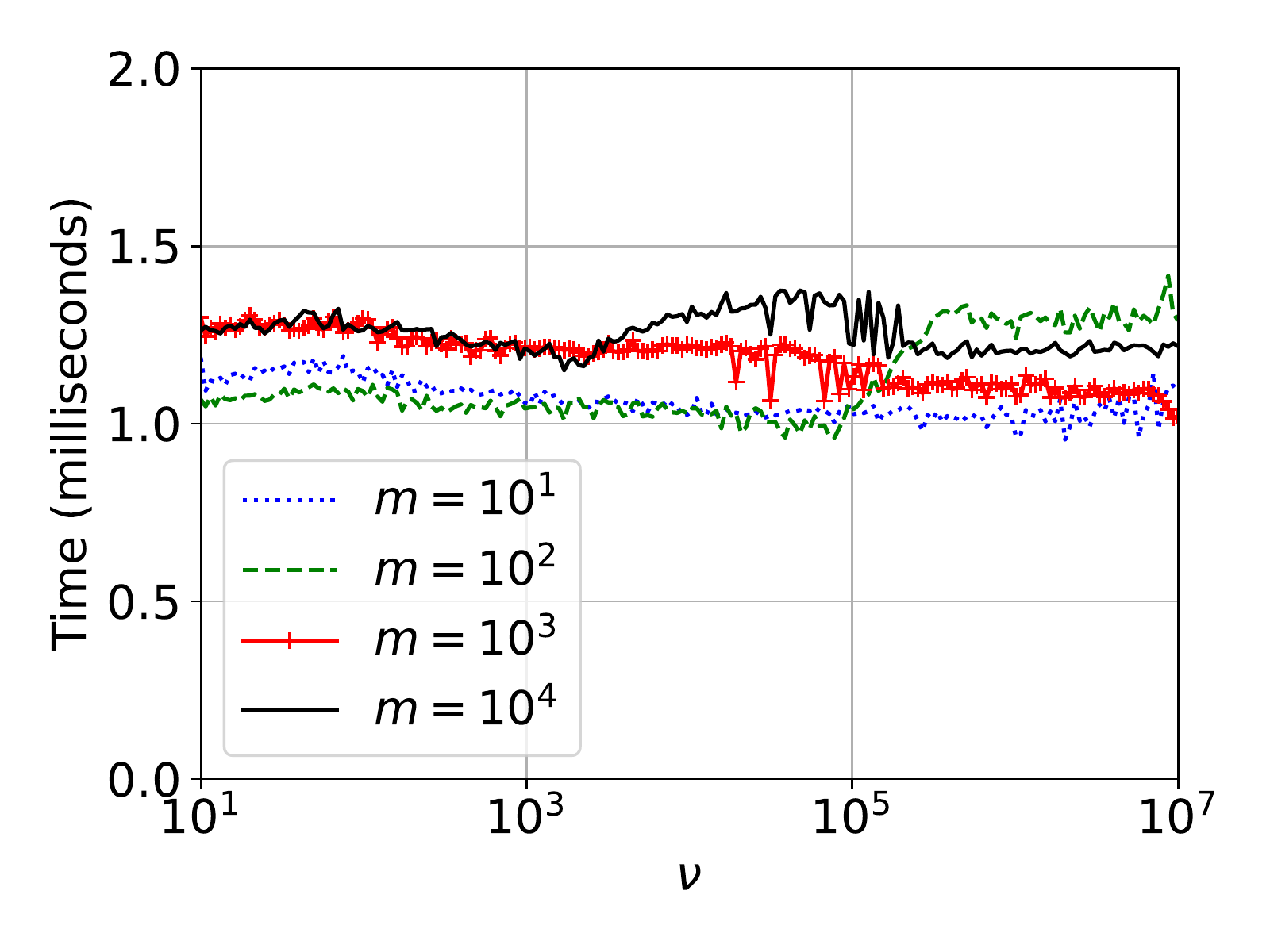}
\hfil

\caption{\small
The results of the first experiment of Section~\ref{section:experiments:alf}.
The plots give the error in the calculated value of  $I_{15}(m+\nu,m)$
and the time required to calculate it via the adaptive Levin method
(incuding the time spent constructing the phase function)
as functions of $\nu$ for $m=10, 10^2, 10^3, 10^4$.
}
\label{figure:alfplots1}
\end{figure}

Because (\ref{experiments:alf:ode}) has singular points at $\pm 1$, it is convenient to introduction the change of variables
\begin{equation}
z(w) = y(\tanh(w)),
\end{equation}
which yields the new differential equation
\begin{equation}
z''(w) + \left(\nu(\nu+1)\sech^2(w) - \mu^2\right)z(w) = 0.
\label{experiments:alf:ode2}
\end{equation}
Applying the algorithm of \cite{BremerPhase2} to (\ref{experiments:alf:ode2})
gives us the representations
\begin{equation}
\begin{aligned}
\frac{\sin\left(\psi_{\nu,\mu}^{\mbox{\tiny alf}} (w) \right)}
{\sqrt{\frac{d}{dw} \psi_{\nu,\mu}^{\mbox{\tiny alf}}(w)}}
= 
\sqrt{\frac{\pi}{2\nu+1}}
\cos(\pi (\mu+1)) \widetilde{P}_\nu^\mu \left(\tanh(w)\right)  \\
-\sqrt{\frac{\pi}{2\nu+1}}
\sin(\pi(\mu+1)) \widetilde{Q}_\nu^\mu \left(\tanh(w)\right)
\end{aligned}
\label{experiments:alf:sine}
\end{equation}
and
\begin{equation}
\begin{aligned}
\frac{\cos\left(\psi_{\nu,\mu}^{\mbox{\tiny alf}} (w) \right)}
{\sqrt{\frac{d}{dw} \psi_{\nu,\mu}^{\mbox{\tiny alf}}(w)}}
= \sqrt{\frac{\pi}{2\nu+1}} \cos(\pi (\mu+1)) \widetilde{P}_\nu^\mu \left(\tanh(w)\right) \\
+ \sqrt{\frac{\pi}{2\nu+1}}\sin(\pi(\mu+1)) \widetilde{Q}_\nu^\mu \left(\tanh(w)\right).
\end{aligned}
\label{experiments:alf:cos}
\end{equation}

\begin{figure}[t!]

\hfil
\includegraphics[width=.40\textwidth]{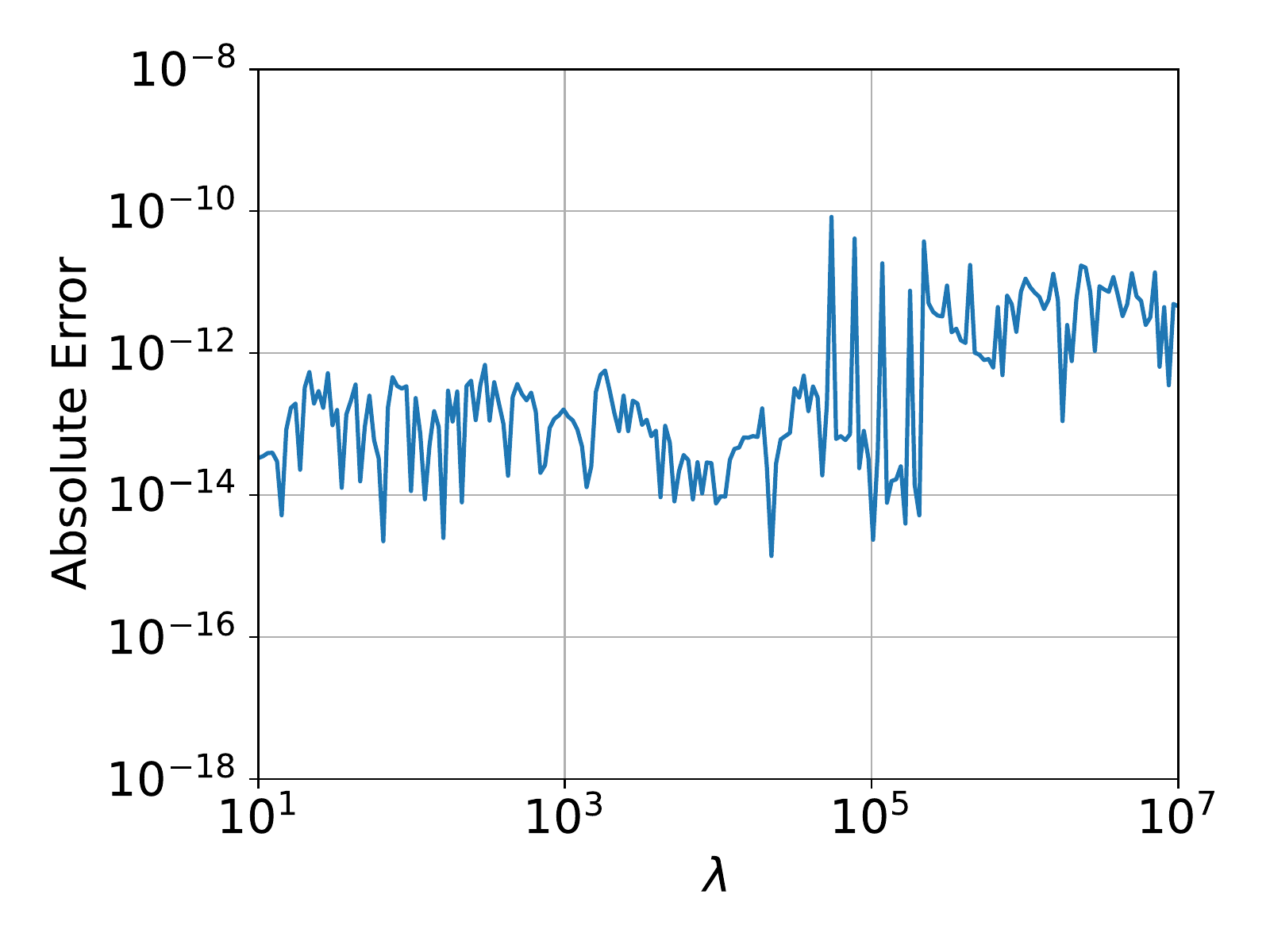}
\hfil
\includegraphics[width=.40\textwidth]{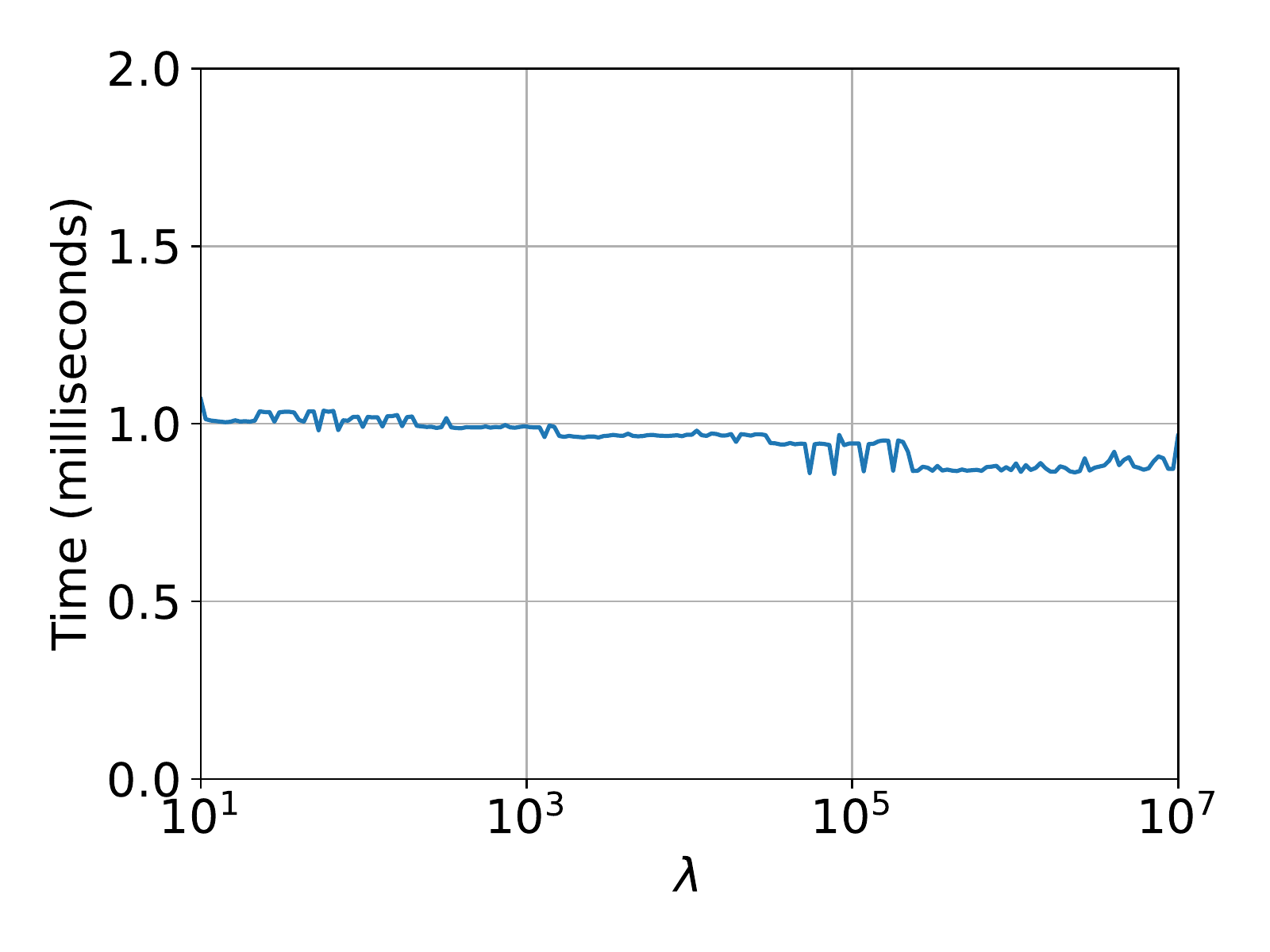}
\hfil

\hfil
\includegraphics[width=.40\textwidth]{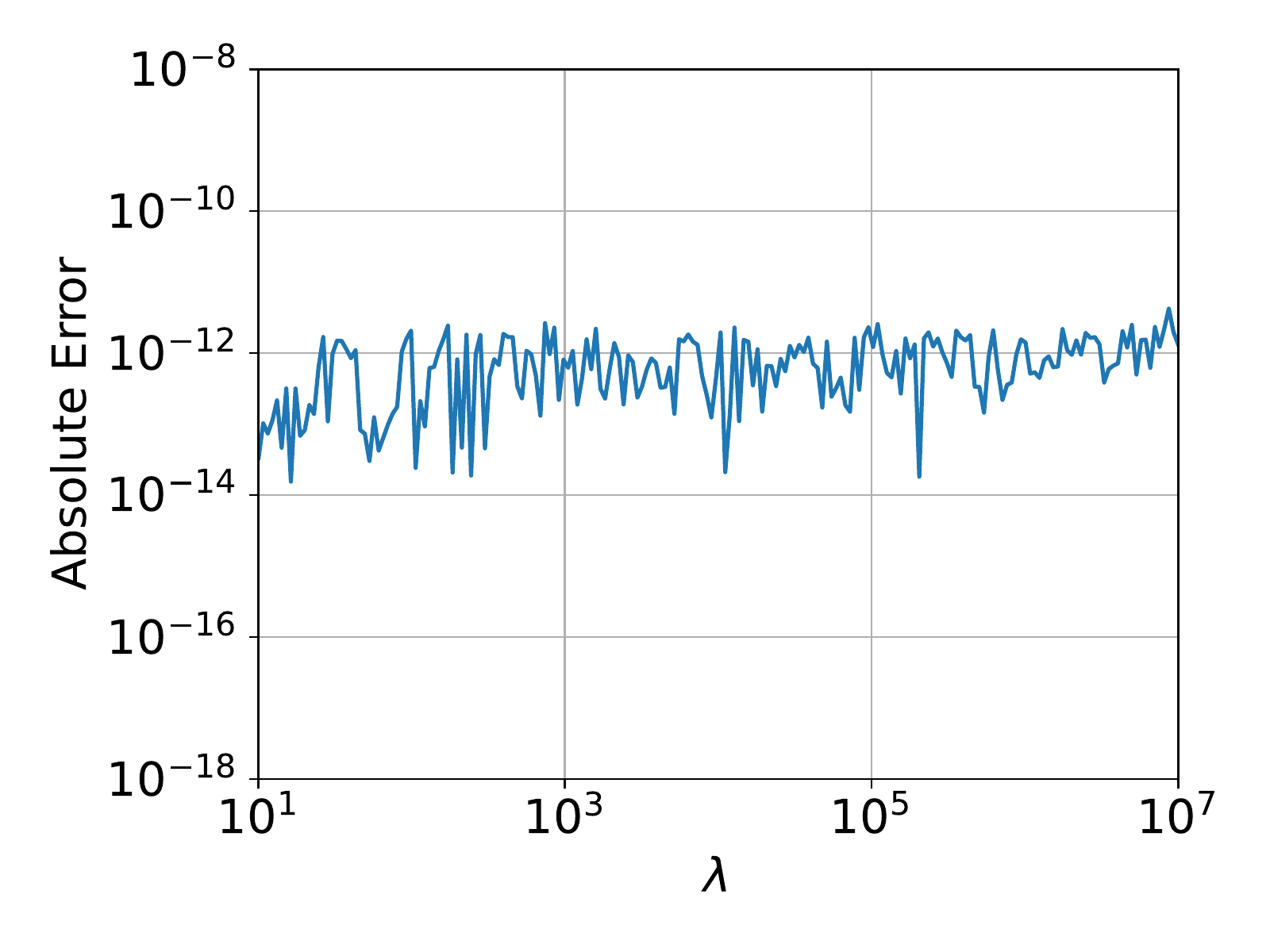}
\hfil
\includegraphics[width=.40\textwidth]{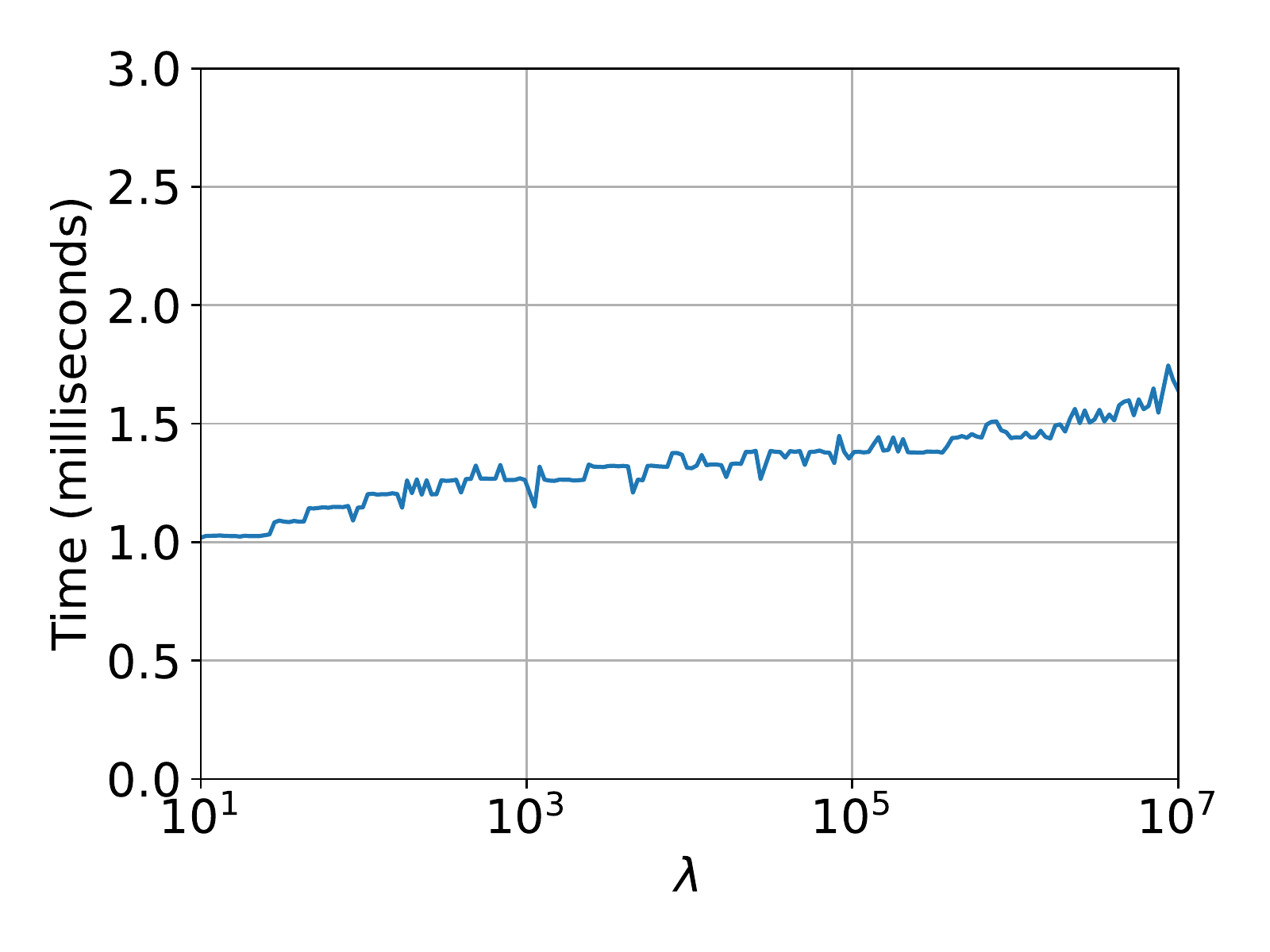}
\hfil

\caption{\small
The results of the last two experiments of Section~\ref{section:experiments:alf}.
The first row pertains to   $I_{16}(\lambda)$ while the second concerns $I_{17}(\lambda)$.
In each row, the plot on the left gives the absolute error in the calculation
of the integral as a function of $\lambda$ and the plot on the right shows
the total time required to compute the integral via the adaptive Levin method
and to construct any necessary phase functions as a function of $\lambda$.
}
\label{figure:alfplots2}
\end{figure}

In our first experiment, we sampled $l=200$ equispaced  points $x_1,\ldots,x_l$ 
in the interval $[1,7]$.  Then, for each $\nu = 10^{x_1},\ldots,10^{x_l}$
and  each $m=10, 10^2, 10^3, 10^4$,  we  constructed the
phase function $\psi^{\mbox{\tiny alf}}_{m+\nu,m}$ 
and applied the adaptive Levin method to the integral
\begin{equation}
(-1)^{m+1} \sqrt{\frac{2\nu+2m+1}{\pi}} \int_0^\infty 
\frac{\sin\left(\psi_{m+\nu,m}^{\mbox{\tiny alf}} (w) \right)}
{\sqrt{\frac{d}{dw} \psi_{m+\nu,m}^{\mbox{\tiny alf}}(w)}
} 
\sech(w)^{2+m}
\, dw
\end{equation}
in order to evaluate
$I_{15}(m+\nu,m)$.  We note that the associated Legendre functions are generally
only defined when the degree is greater than or equal to the order, hence our decision to write the
degree in the form $m + \nu$ with $m$ the order of the associated Legendre function being integrated.
Moreover, by choosing integer values of $m$ we ensured that only
the term involving the function of the first kind was nonzero in (\ref{experiments:alf:sine}).
The results are given in Figure~\ref{figure:alfplots1}.  The timings reported there
include both the cost to apply the adaptive Levin method and the time required to construct
the phase function.

In our second experiment concerning the associated Legendre functions,
we sampled $l=200$ equispaced  points $x_1,\ldots,x_l$ 
in the interval $[1,7]$.  Then, for each $\lambda = 10^{x_1},\ldots,10^{x_l}$,  we  constructed the
phase function $\psi^{\mbox{\tiny alf}}_{\lambda,\frac{1}{2}}$ 
and applied the adaptive Levin method to the integral
\begin{equation}
 \frac{2\lambda+1}{2\pi}
 \int_0^\infty 
\frac{\sin\left(2\psi_{\lambda,\frac{1}{2}}^{\mbox{\tiny alf}} (w) \right)}
{\frac{d}{dw} \psi_{\lambda,\frac{1}{2}}^{\mbox{\tiny alf}}(w)} 
\sech(w)^2
\, dw,
\end{equation}
which is equal to $I_{16}(\lambda)$ by 
virute of the fact that $\sin(2x) = 2 \sin(x)\cos(x)$.
We took the input functions for the adaptive Levin method to be
\begin{equation}
g(x) =  2\psi^{\mbox{\tiny alf}}_{\lambda,\frac{1}{2}}\left(x\right)
\end{equation}
and
\begin{equation}
f(x) =
 \frac{2\lambda+1}{2\pi}
\frac{1}
{\psi^{\mbox{\tiny alf}}_{\lambda,\frac{1}{2}}\left(\mbox{arctanh}\left(\cos(x)\right)\right)}.
\end{equation}
The results are given in the first row of Figure~\ref{figure:alfplots2}.  The timings reported there
include both the cost to apply the adaptive Levin method and the time required to construct
the required phase function.

In our third and final experiment regarding the associated Legendre functions, 
we sampled $l=200$ equispaced  points $x_1,\ldots,x_l$ 
in the interval $[1,7]$.  Then, for each $\lambda = 10^{x_1},\ldots,10^{x_l}$,  we  constructed the
phase function $\psi^{\mbox{\tiny alf}}_{\lambda,1}$
via the algorithm of \cite{BremerPhase2} and evaluated the integral $I_{17}(\lambda)$
via the adaptive Levin method.  The input functions were taken to be 
\begin{equation}
g(x) =  \psi^{\mbox{\tiny alf}}_{\lambda,1}\left(\mbox{arctanh}\left(\cos(x)\right)\right)
\end{equation}
and
\begin{equation}
f(x) = \sqrt{ \frac{2\lambda+1}{2\pi}} 
\frac{1}
{ \sqrt{\psi^{\mbox{\tiny alf}}_{\lambda,1}\left(\mbox{arctanh}\left(\cos(x)\right)\right)}}.
\end{equation}
The results are shown in the second row of Figure~\ref{figure:alfplots2}.  The timings reported there
include both the cost to apply the adaptive Levin method and the time required to construct
the required phase function.

\end{subsection}

\begin{subsection}{Certain integrals involving the Hermite polynomials}
\label{section:experiments:herm}

In the experiments of this subsection, we used the adaptive Levin method to evaluate the integrals
\begin{equation*}
\begin{aligned}
I_{18}(n) &= \int_0^\infty   \widetilde{H}_n(x)\exp\left(\frac{-x^2}{2}\right)\, dx
,
\\
I_{19}(n) &= \int_0^\infty  \widetilde{H}_n(x)(x) \cos\left(n x \right) \exp(-x) \, dx 
\ \ \ \mbox{and}
\\
I_{20}(n) &= \int_0^\infty \widetilde{H}_n(\exp(x))\, dx ,
\end{aligned}
\end{equation*}
where $\widetilde{H_n}$ denotes a normalized version of the Hermite polynomial of 
degree $n$.    The Hermite polynomial $H_n$ is the unique solution of 
\begin{equation}
y''(x)  - 2 x y'(x) + 2 n y(x) = 0
\end{equation}
which decays to $0$ as $x \to \pm \infty$ and such that
\begin{equation}
H_n(0) = \frac{2^n \sqrt{\pi}} { \Gamma\left(\frac{1-n}{2}\right)}.
\end{equation}
We define $\widetilde{H}_n$ via
\begin{equation}
\widetilde{H}_n(x) = \sqrt{\frac{\Gamma(n+1)}{2^n \sqrt{\pi}} }\exp\left(-\frac{x^2}{2}\right) H_n(x);
\end{equation}
it is a solution of 
\begin{equation}
y''(x) + \left(1 + 2n - x^2\right) y(x)=0
\label{experiments:herm:ode}
\end{equation}
whose $L^2(-\infty,\infty)$ norm is $1$.
Applying the algorithm of \cite{BremerPhase2} to (\ref{experiments:herm:ode})
gives us the representation 
\begin{equation}
  \widetilde{H}_n \left(x\right)  = C_n^{\mbox{\tiny herm}} \frac{\sin\left(\psi_{n}^{\mbox{\tiny herm}} (x) \right)}
{\sqrt{\frac{d}{dx} \psi_{n}^{\mbox{\tiny herm}}(x)}}
\label{experiments:herm:sine}
\end{equation}
of    the normalized  Hermite polynomial of degree $n$.  The authors are not aware of a convenient
formula for the constant $C_n^{\mbox{\tiny herm}}$; we determine it numerically after the phase function
has been constructed.

In our first experiment, we sampled $l=200$ equispaced  points $x_1,\ldots,x_l$ 
in the interval $[1,7]$.  Then, for each $n = \left\lfloor10^{x_1}\right\rfloor,\ldots,\left\lfloor10^{x_l}\right\rfloor$,  we  constructed the
phase function $\psi^{\mbox{\tiny herm}}_{n}$ 
and used it to evaluate $I_{18}(n)$.
The results are given in the first row of Figure~\ref{figure:hermplots1}.  The timings reported there
include both the cost to apply the adaptive Levin method and the time required to construct
the phase function.

In a second experiment concerning Hermite polynomials, we 
 we sampled $l=200$ equispaced  points $x_1,\ldots,x_l$ 
in the interval $[1,7]$. 
 Then, for each $n = \left\lfloor10^{x_1}\right\rfloor,\ldots,\left\lfloor10^{x_l}\right\rfloor$,  we  constructed the
phase function $\psi^{\mbox{\tiny herm}}_{n}$  and used it to evaluate the integrals
\begin{equation}
I_{19a} (n) =  \int_0^\infty
\frac{\sin\left(\psi^{\mbox{\tiny herm}}_{n}(x) -n  x \right)}
{\sqrt{\  \frac{d}{dx} \psi^{\mbox{\tiny herm}}_{n}(x)}}
 \exp(-x)\, dx
\end{equation}
and
\begin{equation}
I_{19b} (n) =   \int_0^\infty
\frac{\sin\left(\psi^{\mbox{\tiny herm}}_{n}(x) + n x \right)}
{\sqrt{\  \frac{d}{dx} \psi^{\mbox{\tiny herm}}_{n}(x)}}
 \exp(-x)\, dx
\end{equation}
The value of $I_{19}(n)$ is then given by 
\begin{equation}
I_{19}(n) =  C_n^{\mbox{\tiny herm}} \frac{I_{19a}\left(n\right) + I_{19b}\left(n\right)}{2}
\end{equation}
since
\begin{equation}
\sin(x)\cos(y) = \frac{\sin(x-y) + \sin(x+y)}{2}.
\end{equation}
The results of this experiment are given in the second row of Figure~\ref{figure:hermplots1}.  Again, the reported timings 
include both the cost to apply the adaptive Levin method and the time required to construct
the phase function.

We began our third and final experiment regarding the Hermite polynomials by
sampling $l=200$ equispaced  points $x_1,\ldots,x_l$  in the interval $[1,7]$. 
 Then, for each $n = \left\lfloor10^{x_1}\right\rfloor,\ldots,\left\lfloor10^{x_l}\right\rfloor$,  we  constructed the
phase function $\psi^{\mbox{\tiny herm}}_{n}$  and used it to evaluate $I_{20}(n)$.  We took the input
functions for the adaptive Levin method to be 
\begin{equation}
g(x) =\psi^{\mbox{\tiny herm}}_{n}(\exp(x))
\ \ \ \mbox{and} \ \ \ 
f(x) =  C_n^{\mbox{\tiny herm}} \frac{1}{\sqrt{\frac{d}{dx}\psi^{\mbox{\tiny herm}}_{n}(\exp(x)) }}.
\end{equation}
The results are shown in the third row of Figure~\ref{figure:hermplots1}.  The timings reported there
include both the cost to apply the adaptive Levin method and the time required to construct
the phase function.

\begin{figure}[h!]

\hfil
\includegraphics[width=.40\textwidth]{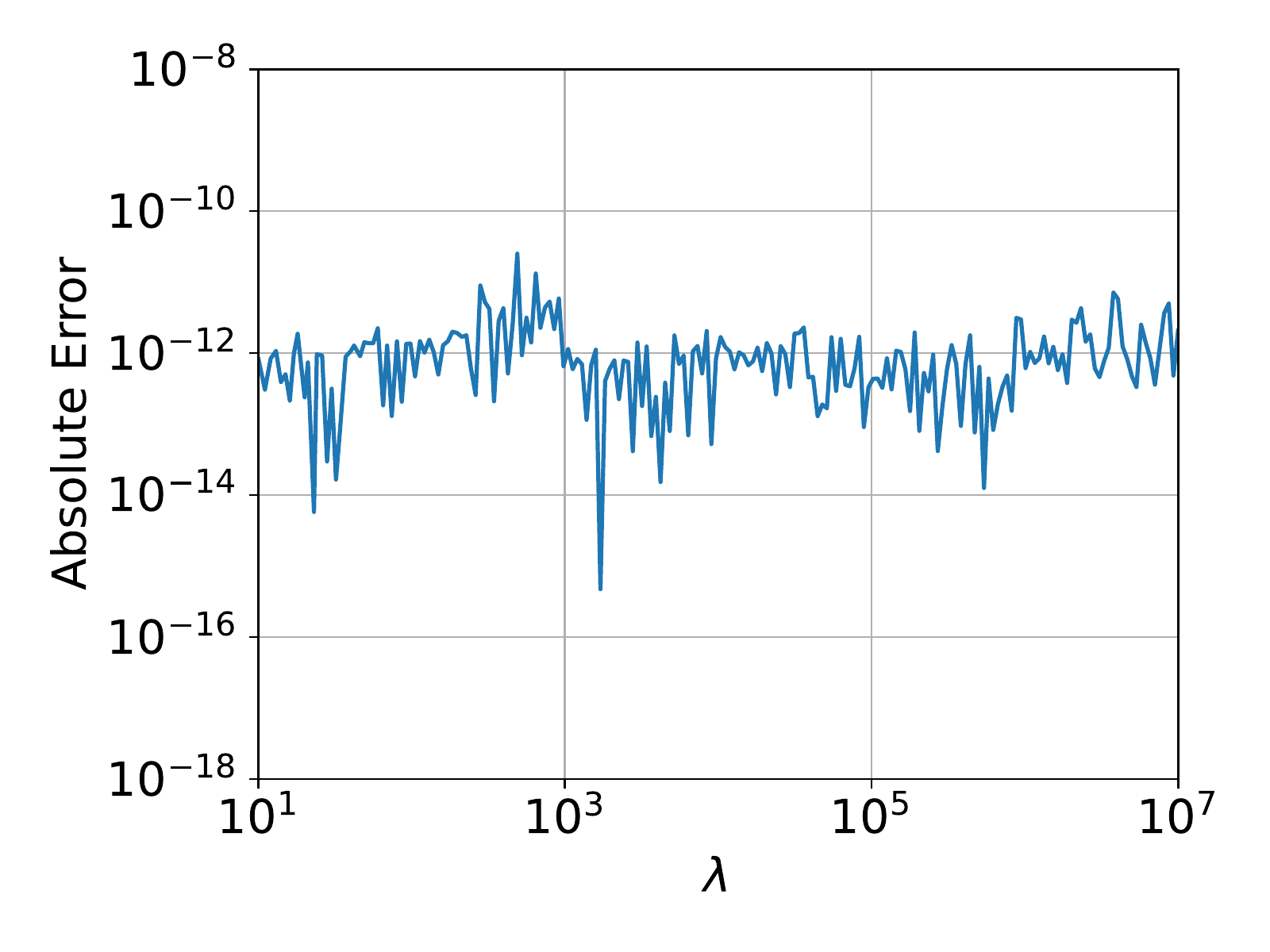}
\hfil
\includegraphics[width=.40\textwidth]{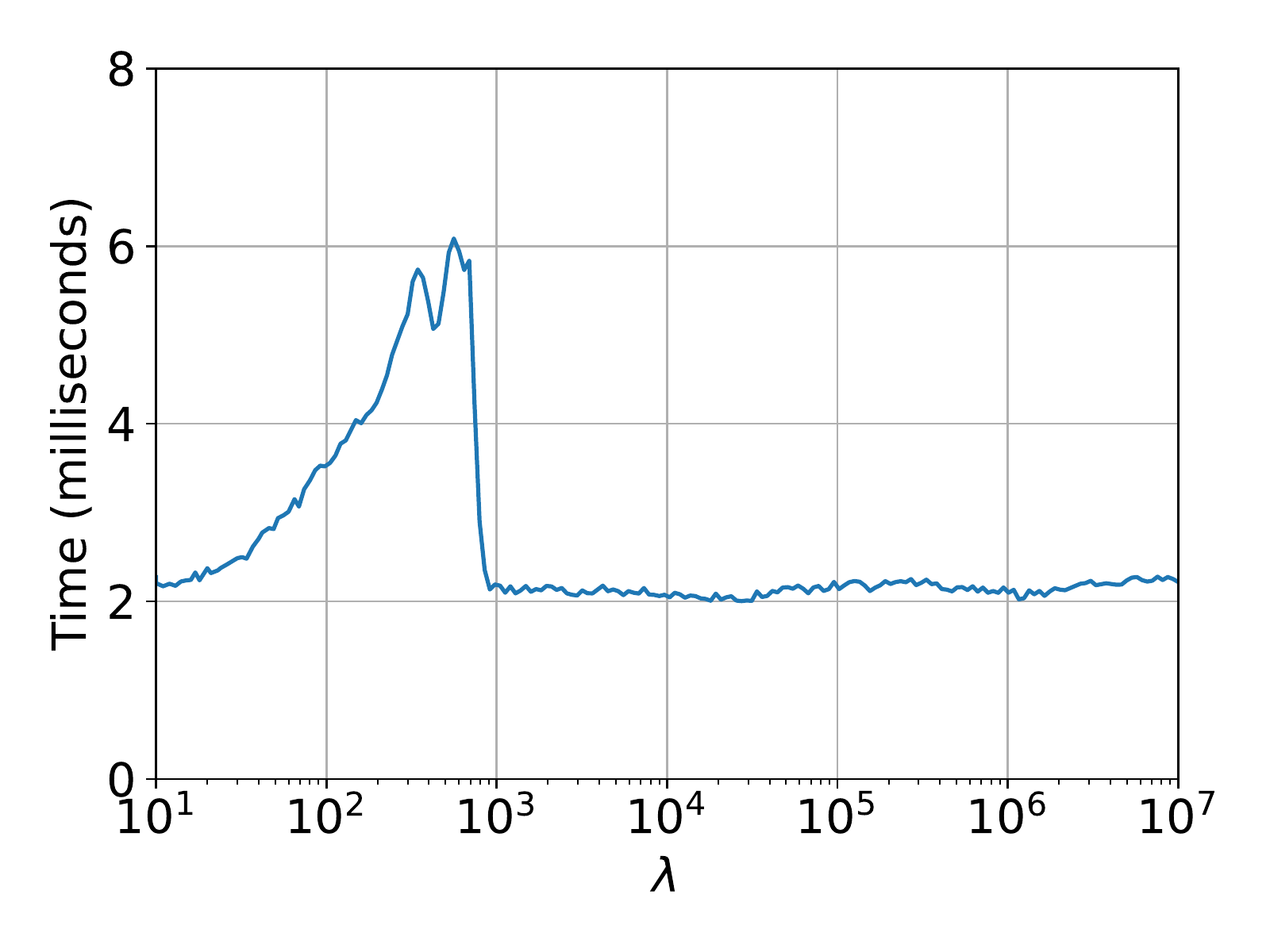}
\hfil

\hfil
\includegraphics[width=.40\textwidth]{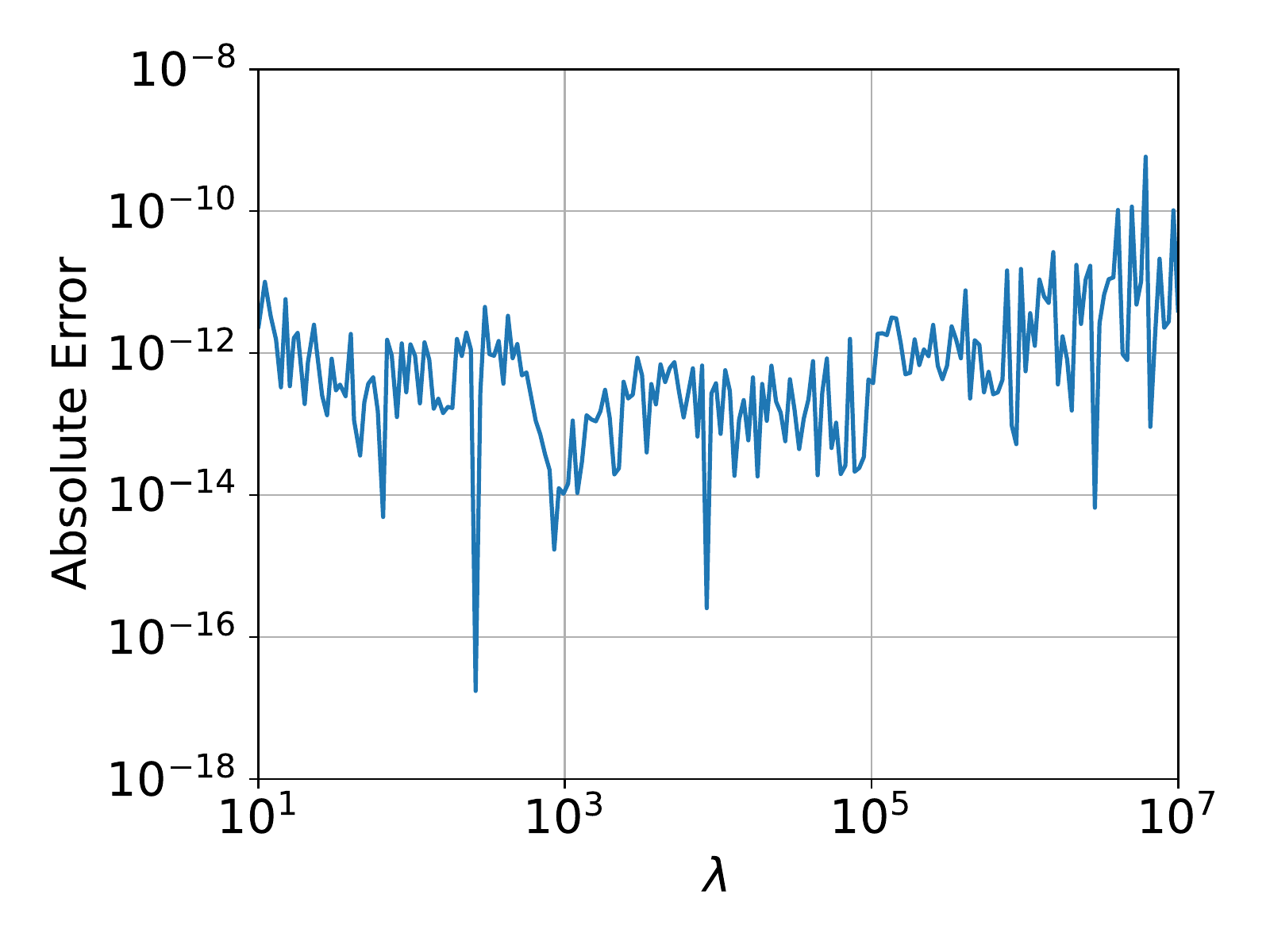}
\hfil
\includegraphics[width=.40\textwidth]{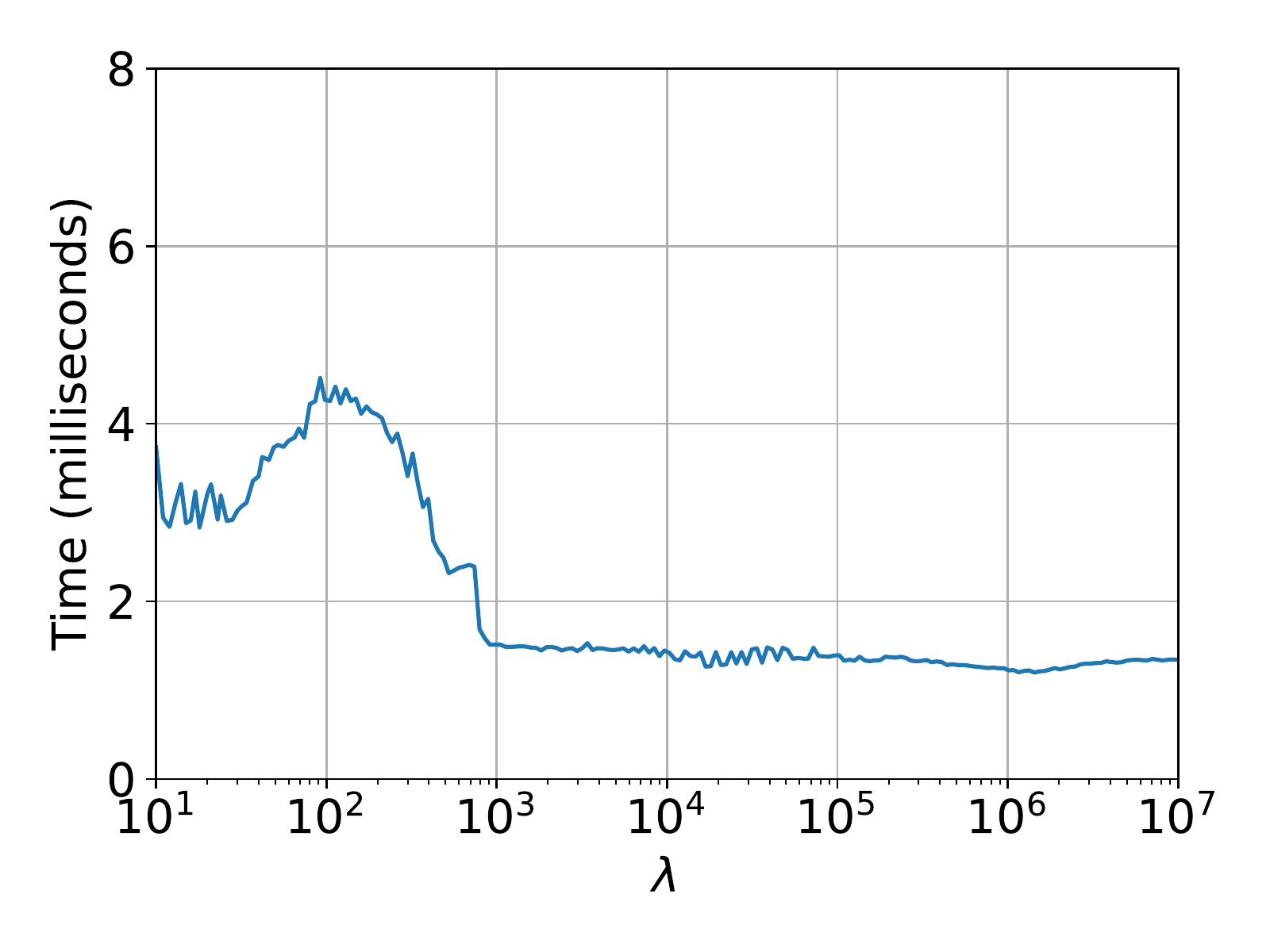}
\hfil

\hfil
\includegraphics[width=.40\textwidth]{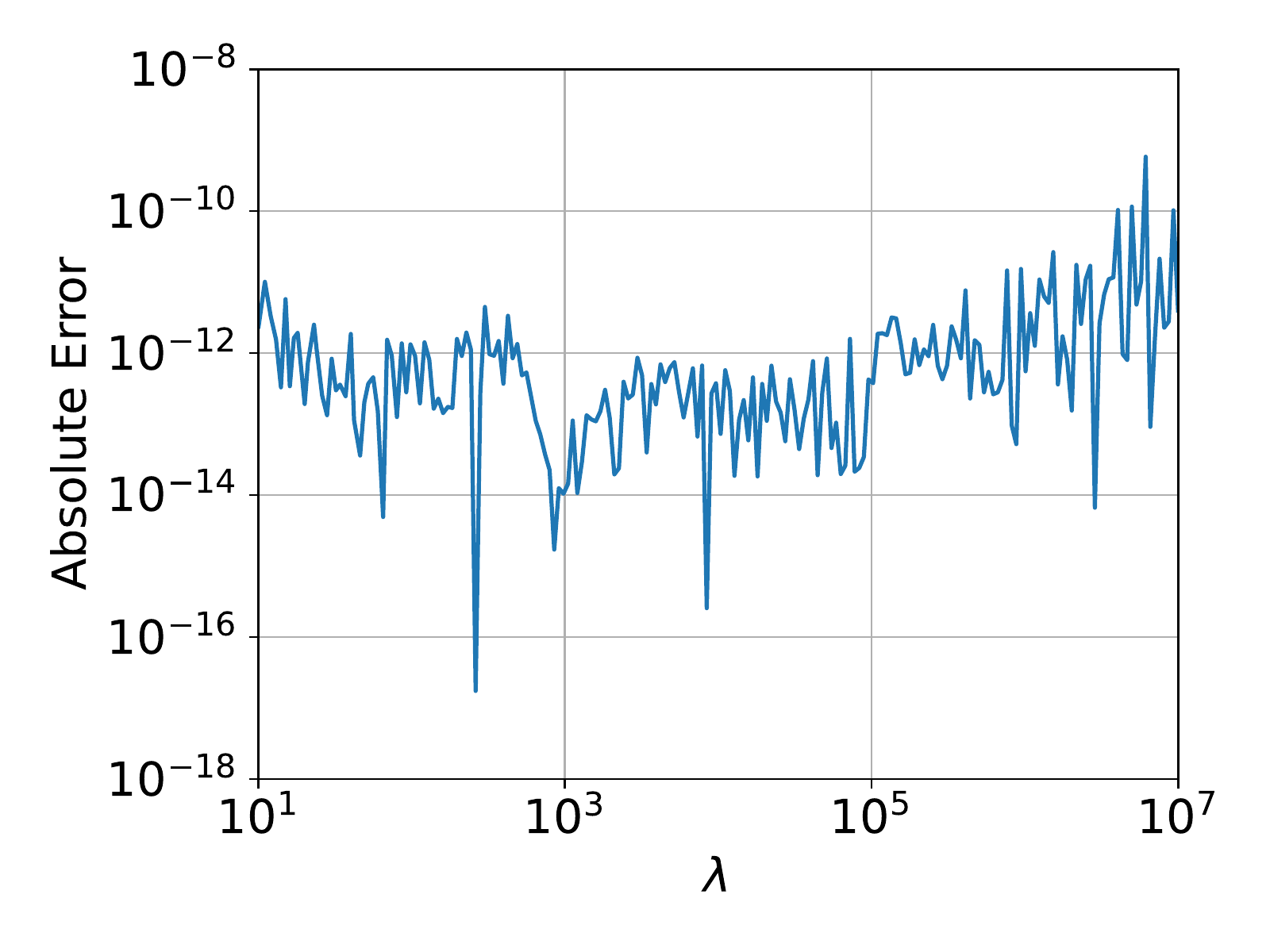}
\hfil
\includegraphics[width=.40\textwidth]{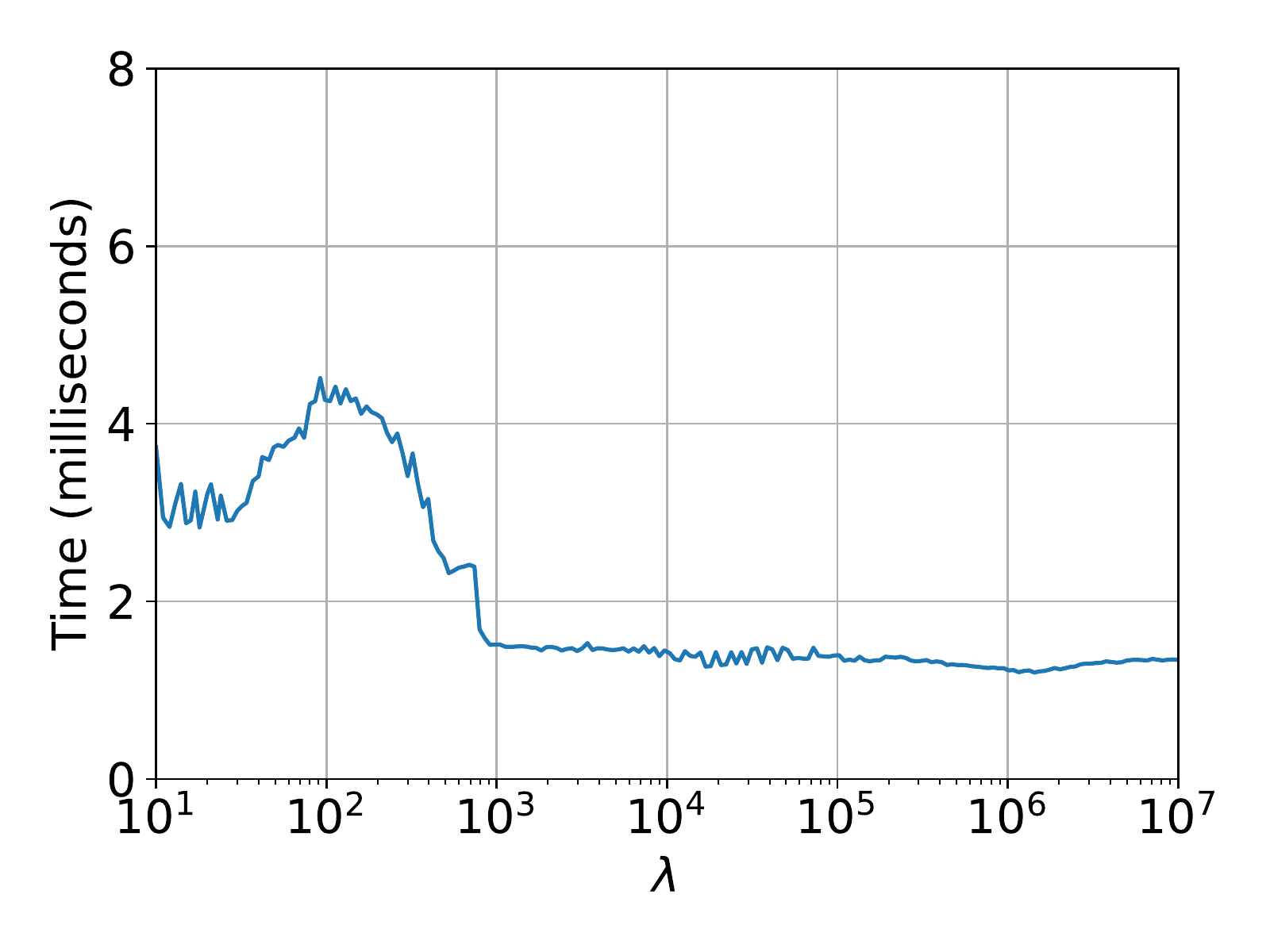}
\hfil

\caption{\small
The results of the experiments of Section~\ref{section:experiments:herm}.
The first row pertains to   $I_{18}(\lambda)$, the second concerns $I_{19}(\lambda)$
and the third pertains to $I_{20}(\lambda)$.
In each row, the plot on the left gives the absolute error in the calculation
of the integral as a function of $\lambda$ and the plot on the right shows
the total time required to compute the integral via the adaptive Levin method
and to construct any necessary phase functions as a function of $\lambda$.
}
\label{figure:hermplots1}
\end{figure}

\end{subsection}

\begin{subsection}{Evaluation of Modal Green's functions}
\label{section:experiments:modal}

\begin{figure}[!h]

\hfil
\includegraphics[width=.40\textwidth]{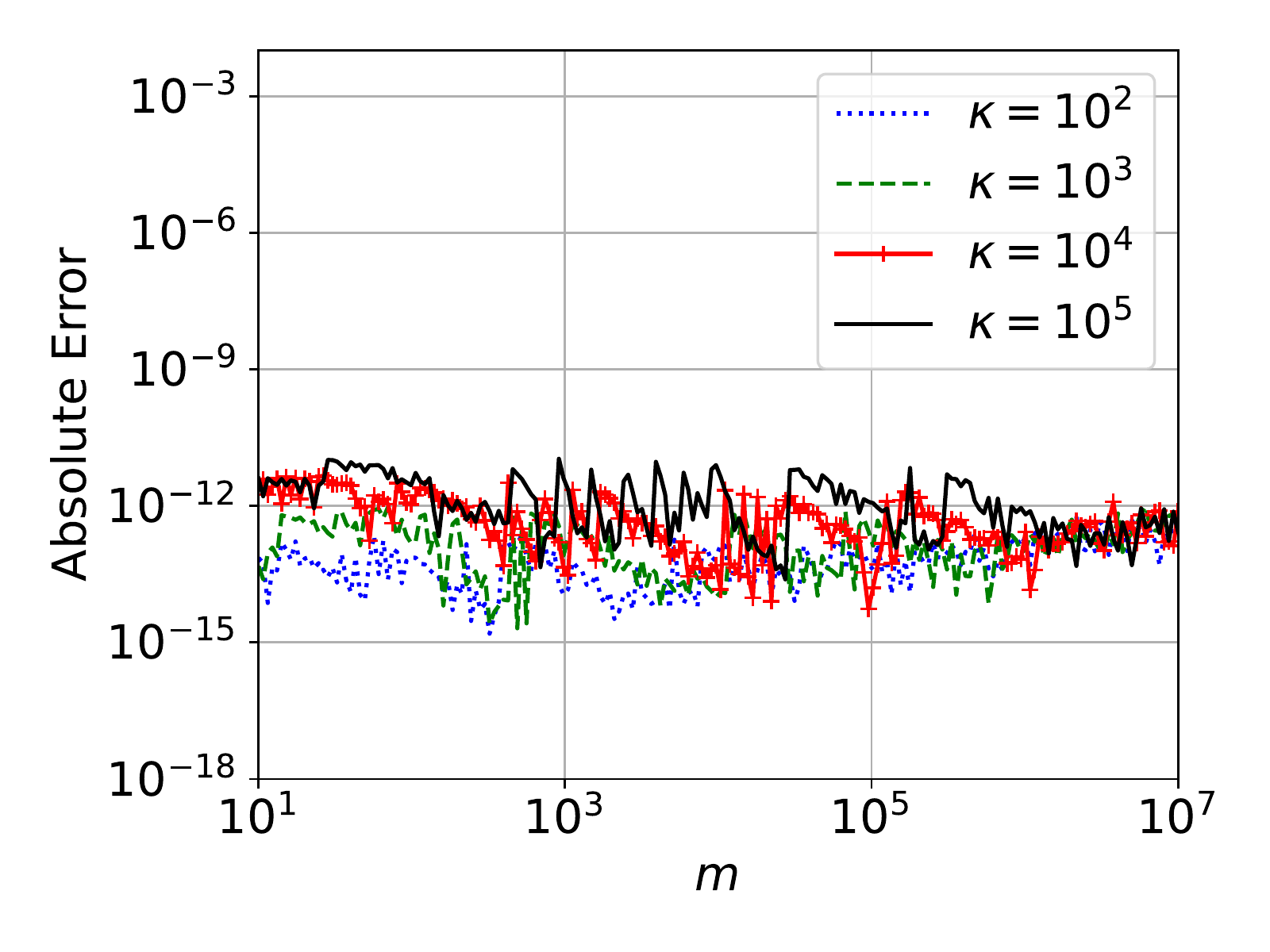}
\hfil
\includegraphics[width=.40\textwidth]{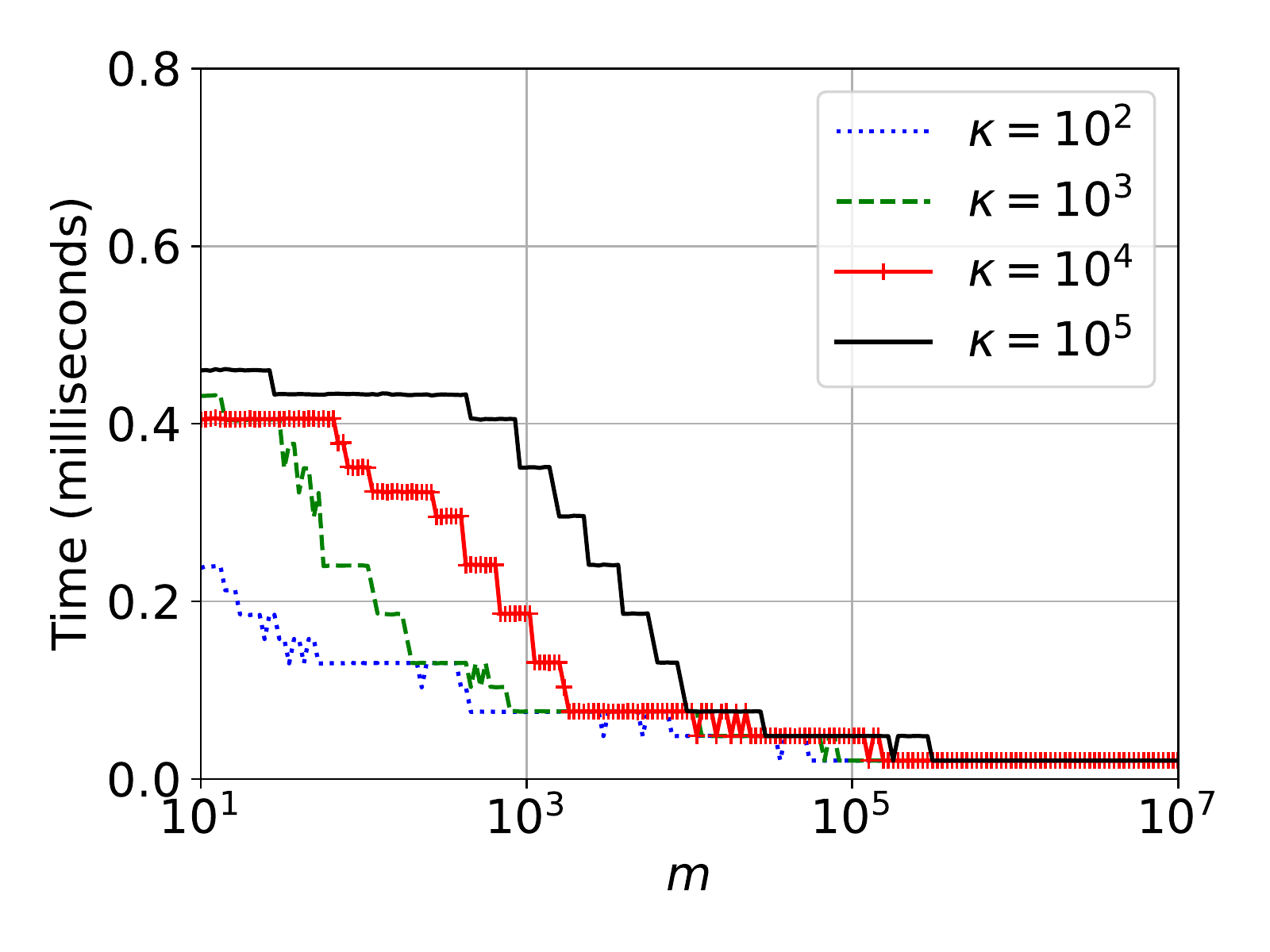}
\hfil

\hfil
\includegraphics[width=.40\textwidth]{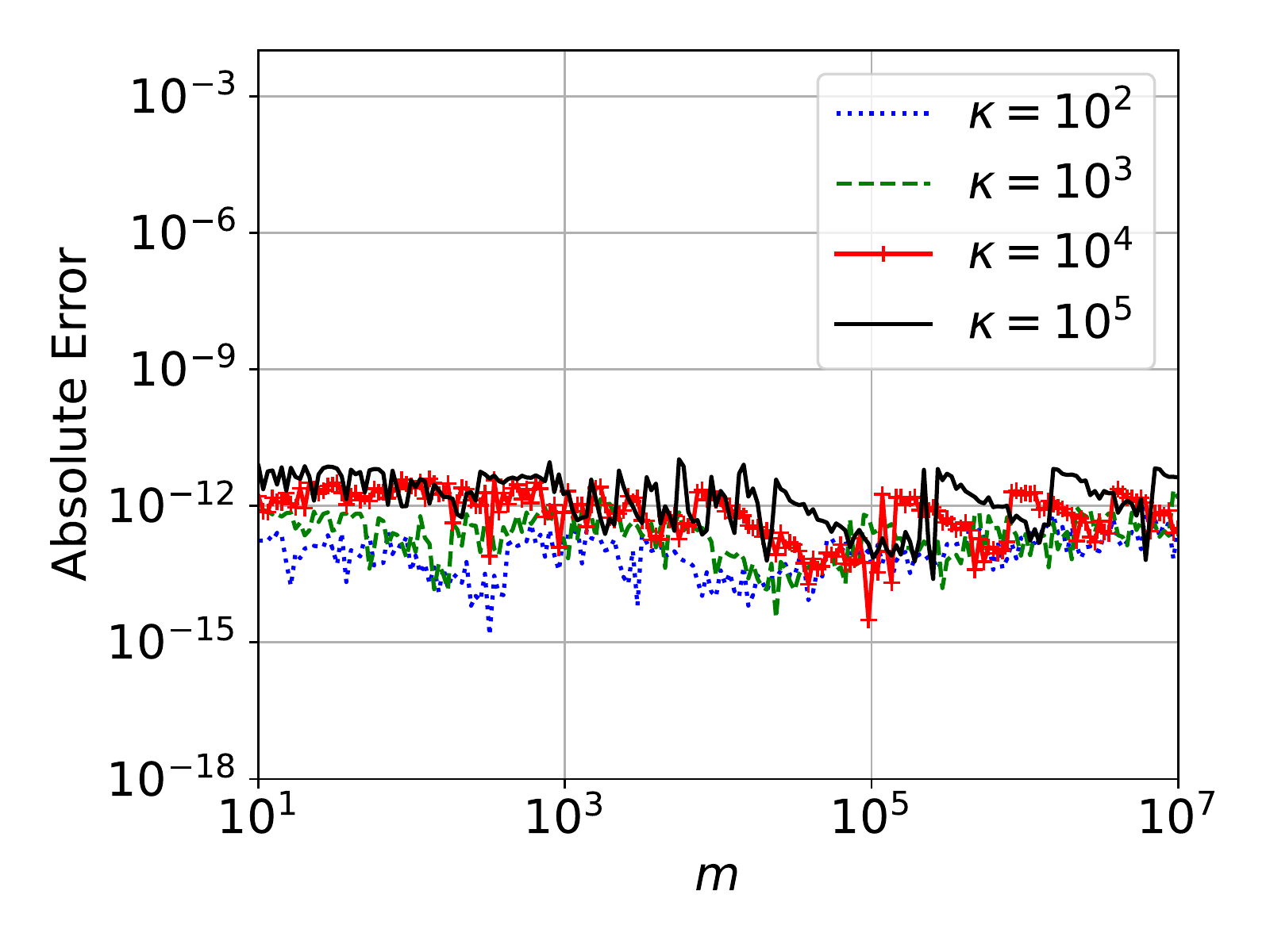}
\hfil
\includegraphics[width=.40\textwidth]{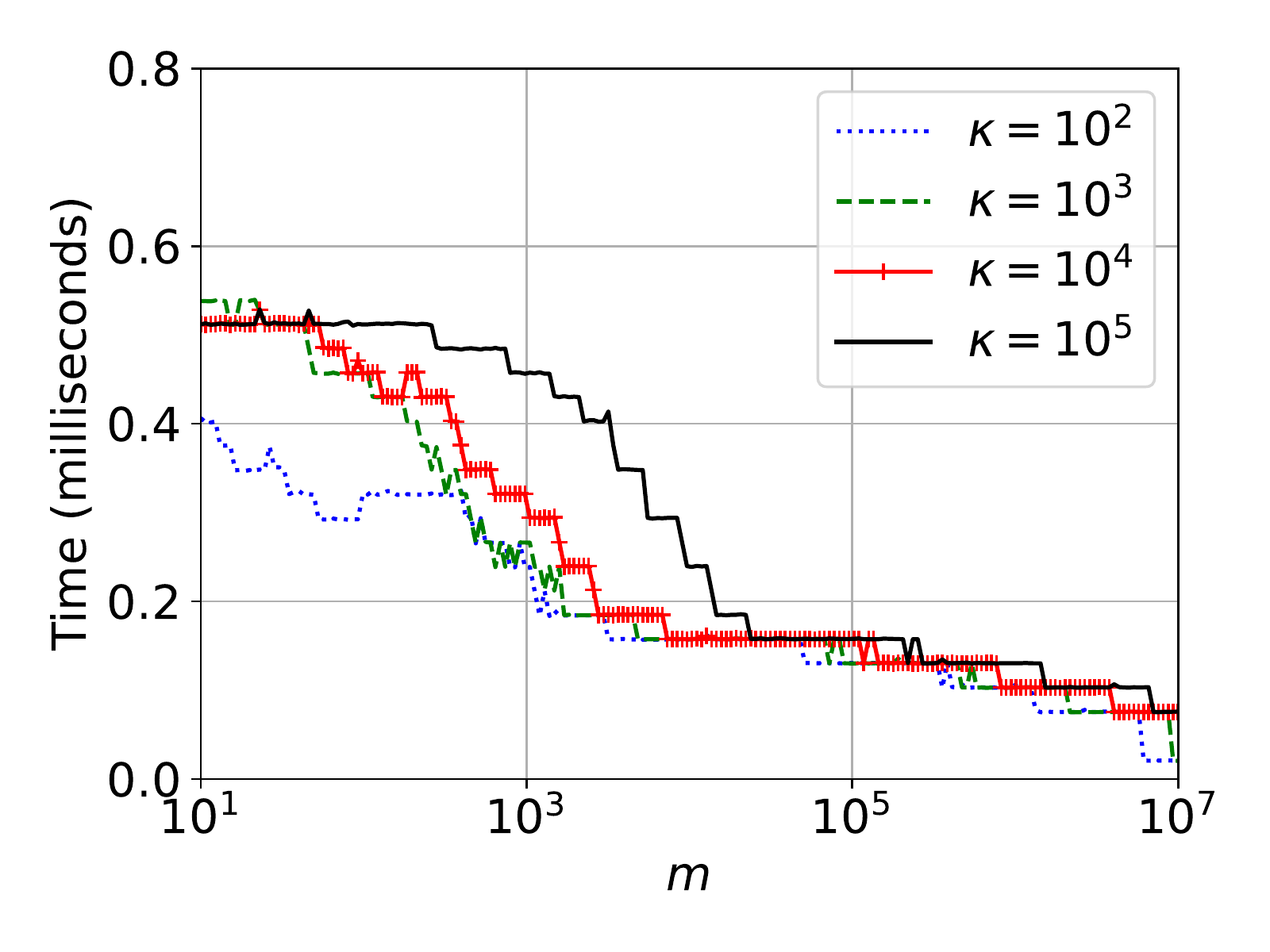}
\hfil

\caption{\small
The results of the first experiment of Section~\ref{section:experiments:modal}.
The plot on the left-hand side of  the first row gives the absolute error in the calculated value of  $I_{21}(\kappa,m,\alpha)$
as a function of $m$ for four values of $\kappa$ in the case $\alpha=0.5$.
The plot on the right-hand side of  the first row gives the time required to calculate  $I_{21}(\kappa,m,\alpha)$
via the adaptive Levin method as a function of $m$ for four values of $\kappa$ in the case $\alpha=0.5$.
The plots in the second row give analagous data in the event $\alpha=0.99$.
}
\label{figure:modalplots1}
\end{figure}

In this set of experiments, we used the adaptive Levin method to evaluate
the azimuthal Fourier components of the Green's function
\begin{equation}
G(x,x') = \frac{\exp(ik|x-x'|)}{4\pi |x-x'|}
\end{equation}
for the Helmholtz equation in three dimensions.  These functions, which are known
as the modal Green's functions for the Helmholtz equation, are given by 
\begin{equation}
\frac{1}{2\pi} \int_{-\pi}^\pi \frac{\exp(ik|x-x'|)}
{4\pi|x-x'| } \exp(-im\theta)\, d\theta.
\label{experiments:modal:gm0}
\end{equation}
By introducing cylindrical coordinates $x = (r,\theta,z)$, 
$x' = (r',\theta',z')$ and letting $\phi = \theta-\theta'$,
we can rewrite (\ref{experiments:modal:gm0}) as
\begin{equation}
I_{21}(\kappa,m, \alpha) =
\frac{1}{4\pi^2} \int_{-\pi}^\pi 
\frac{\exp\left(-i\kappa\sqrt{1-\alpha\cos(\phi)}\right)}{\sqrt{1-\alpha\cos(\phi)}}\cos(m\phi)\,d\phi,
\end{equation}
where
\begin{equation}
\begin{aligned}
\kappa &= kR_0,\ \  \alpha &= \frac{2rr'}{R_0},\ \ \mbox{and}\ \  R_0^2 = r^2 + (r')^2+(z-z')^2.
\end{aligned}
\end{equation}

\begin{figure}[!h]
\hfil
\includegraphics[width=.40\textwidth]{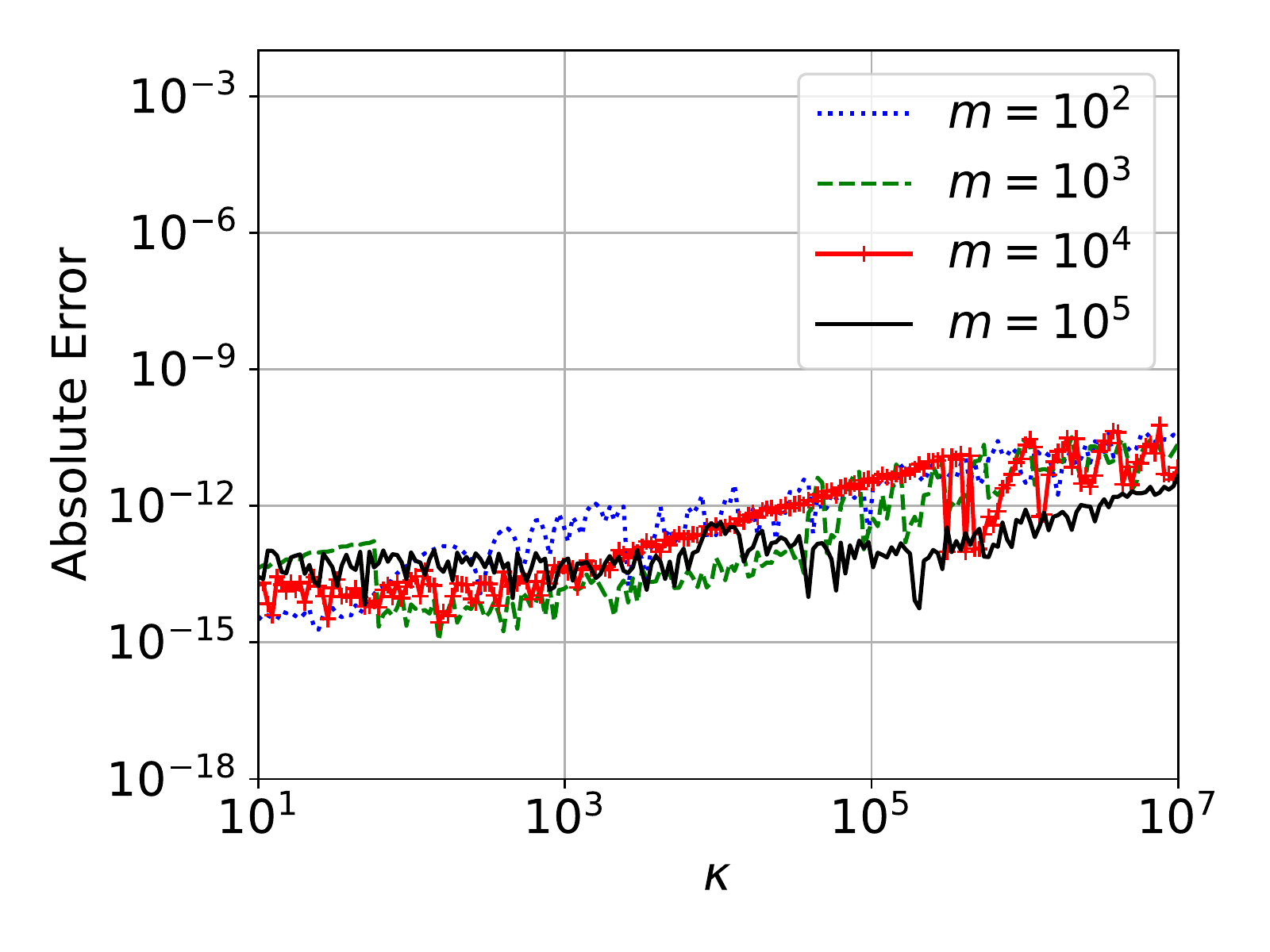}
\hfil
\includegraphics[width=.40\textwidth]{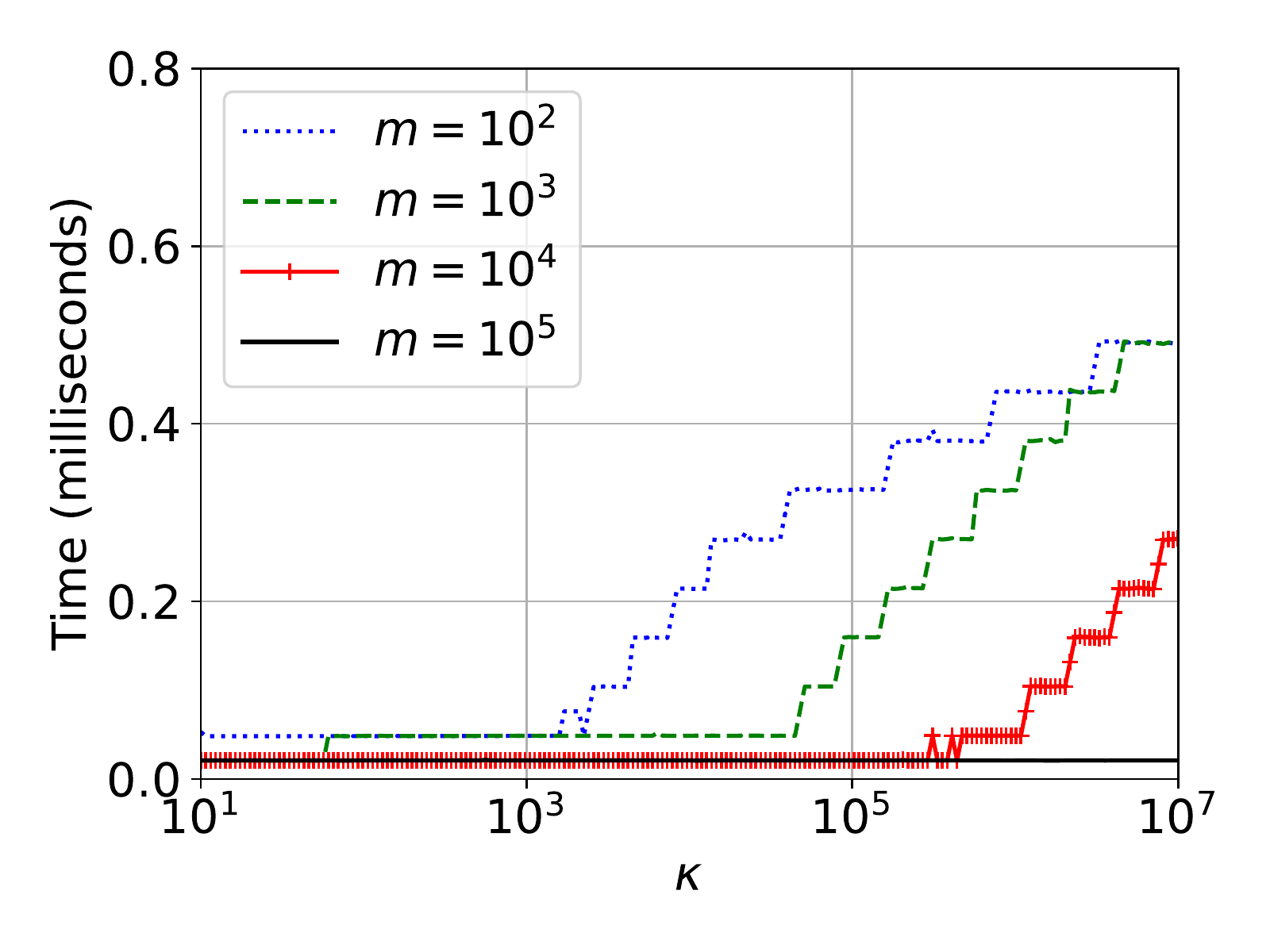}
\hfil

\hfil
\includegraphics[width=.40\textwidth]{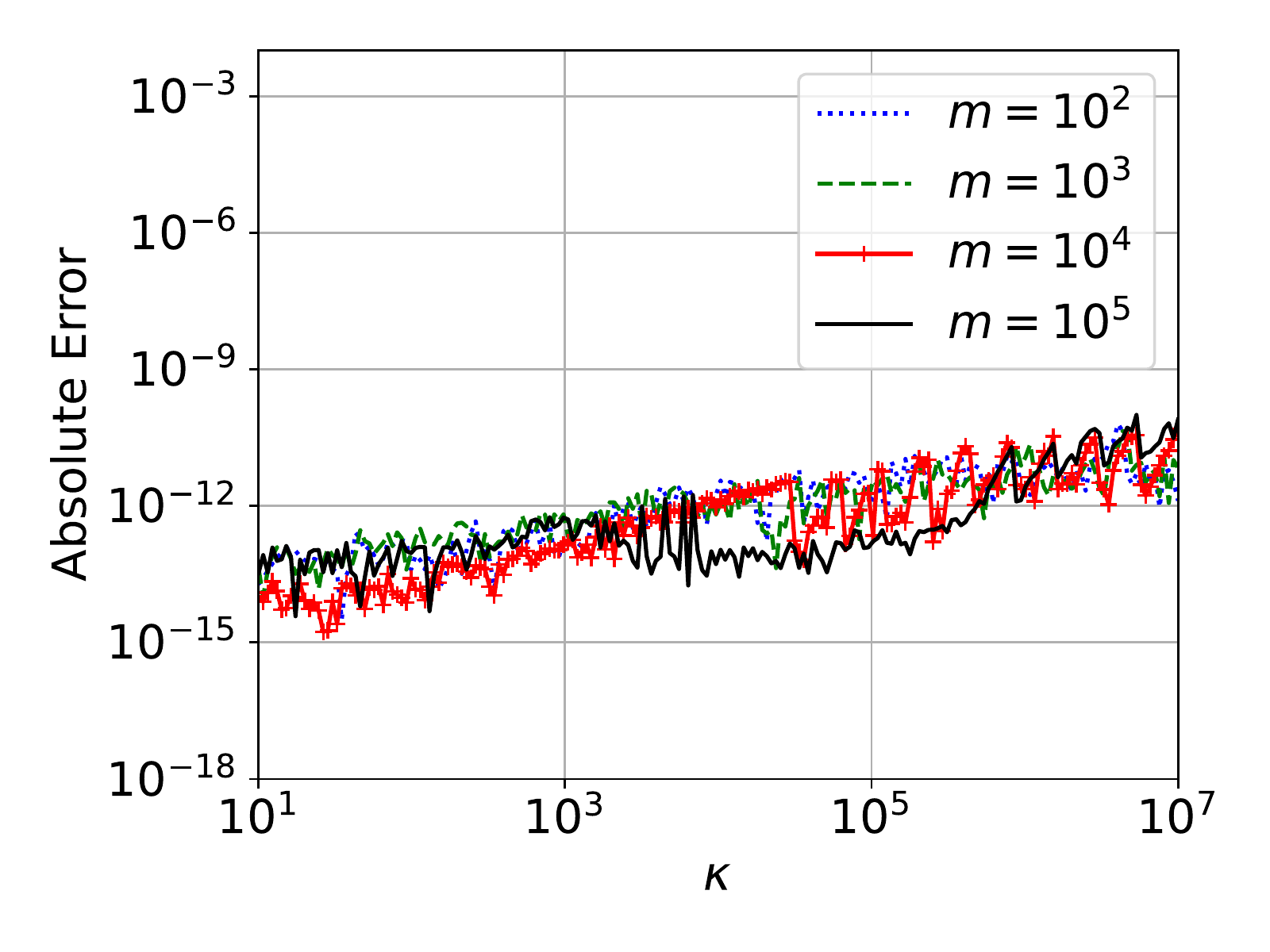}
\hfil
\includegraphics[width=.40\textwidth]{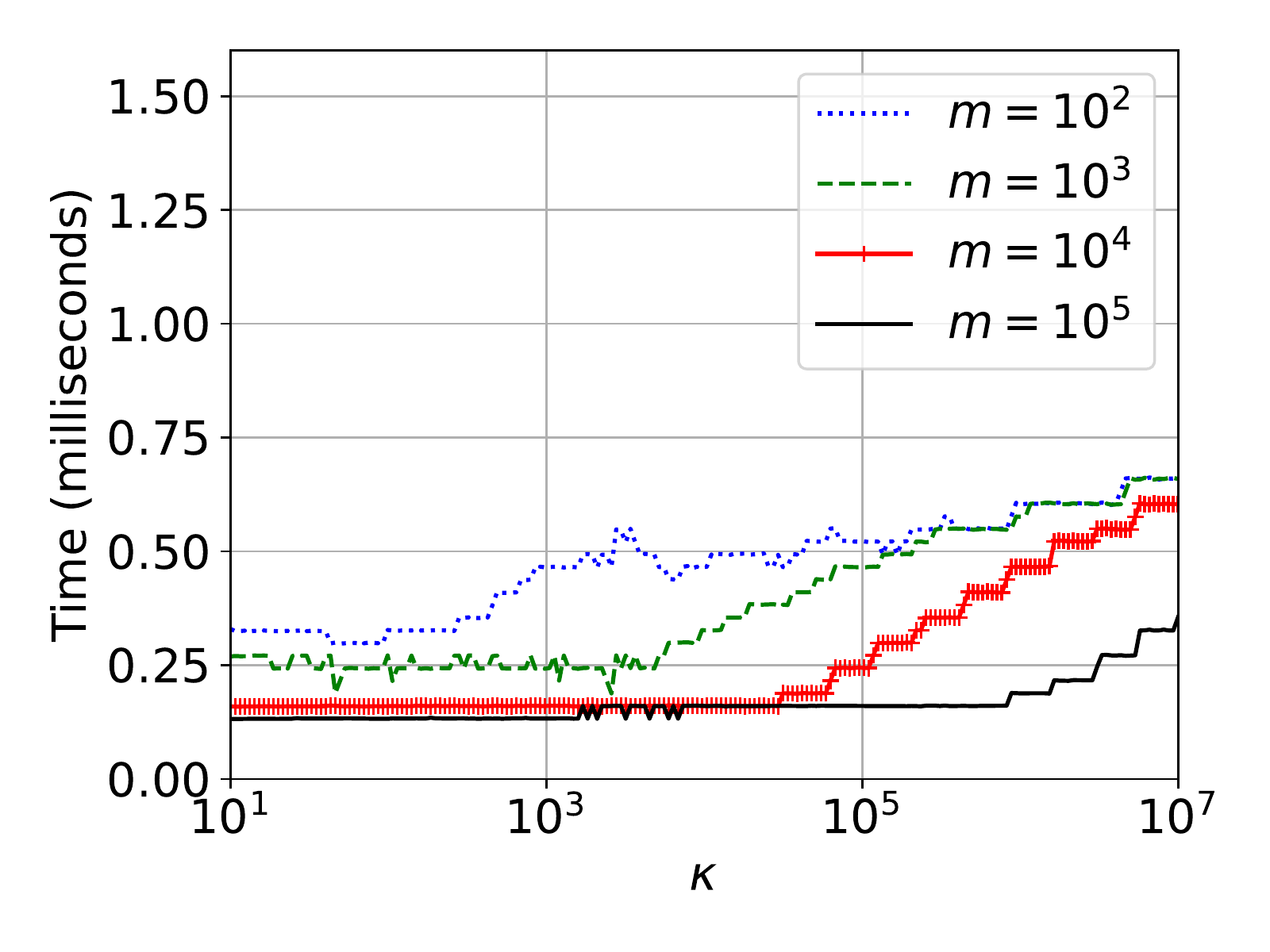}
\hfil

\caption{\small
The results of the second experiment of Section~\ref{section:experiments:modal}.
The plot on the left-hand side of  the first row gives the absolute error in the calculated value of  $I_{21}(\kappa,m,\alpha)$
as a function of $\kappa$ for four values of $m$ in the case $\alpha=0.5$.
The plot on the right-hand side of  the first row gives the time required to calculate  $I_{21}(\kappa,m,\alpha)$
via the adaptive Levin method as a function of $\kappa$ for four values of $m$ in the case $\alpha=0.5$.
The plots in the second row give analagous data in the event $\alpha=0.99$.
}
\label{figure:modalplots2}
\end{figure}

In the first experiment, we sampled 
$l=200$ equispaced  points $x_1,\ldots,x_l$ 
in the interval $[1,7]$. Then, for each $m = 10^{x_1},\ldots,10^{x_l}$,
  $\kappa=10^2, 10^3, 10^4, 10^5$ and $\alpha=0.5,0.99$ we 
evaluated $I_{21}(\kappa,m,\alpha)$ using the adaptive Levin
method and via adaptive Gaussian quadrature. 
Figure~\ref{figure:modalplots1}
gives the results.  

In a second experiment,  we sampled $l=200$ equispaced  points $x_1,\ldots,x_l$ 
in the interval $[1,7]$.  Then, for each $\kappa = 10^{x_1},\ldots,10^{x_l}$,
  $m=10^2, 10^3, 10^4, 10^5$ and $\alpha=0.5, 0.99$
we evaluated $I_{21}(\kappa,m,\alpha)$ using the adaptive Levin
method and via adaptive Gaussian quadrature.  
In this experiment, the tolerance parameters for both the adaptive Levin method
and our adaptive Gaussian quadrature code were
 set to $\epsilon_0 \sqrt{\kappa}$, where $\epsilon_0$ is machine zero
for IEEE double precision arithmetic (about $2.22\times 10^{-16}$).
The value of $I_{21}(\kappa,m,\alpha)$ decreases, but not at a sufficient
rate to completely compensate for the growth in its condition
number, particularly when $\alpha$ is close to $1$.
Hence, the need to allow the tolerance parameter to increase with $\kappa$.
Figure~\ref{figure:modalplots2}
gives the results.  

\begin{figure}[h!!]

\hfil
\includegraphics[width=.39\textwidth]{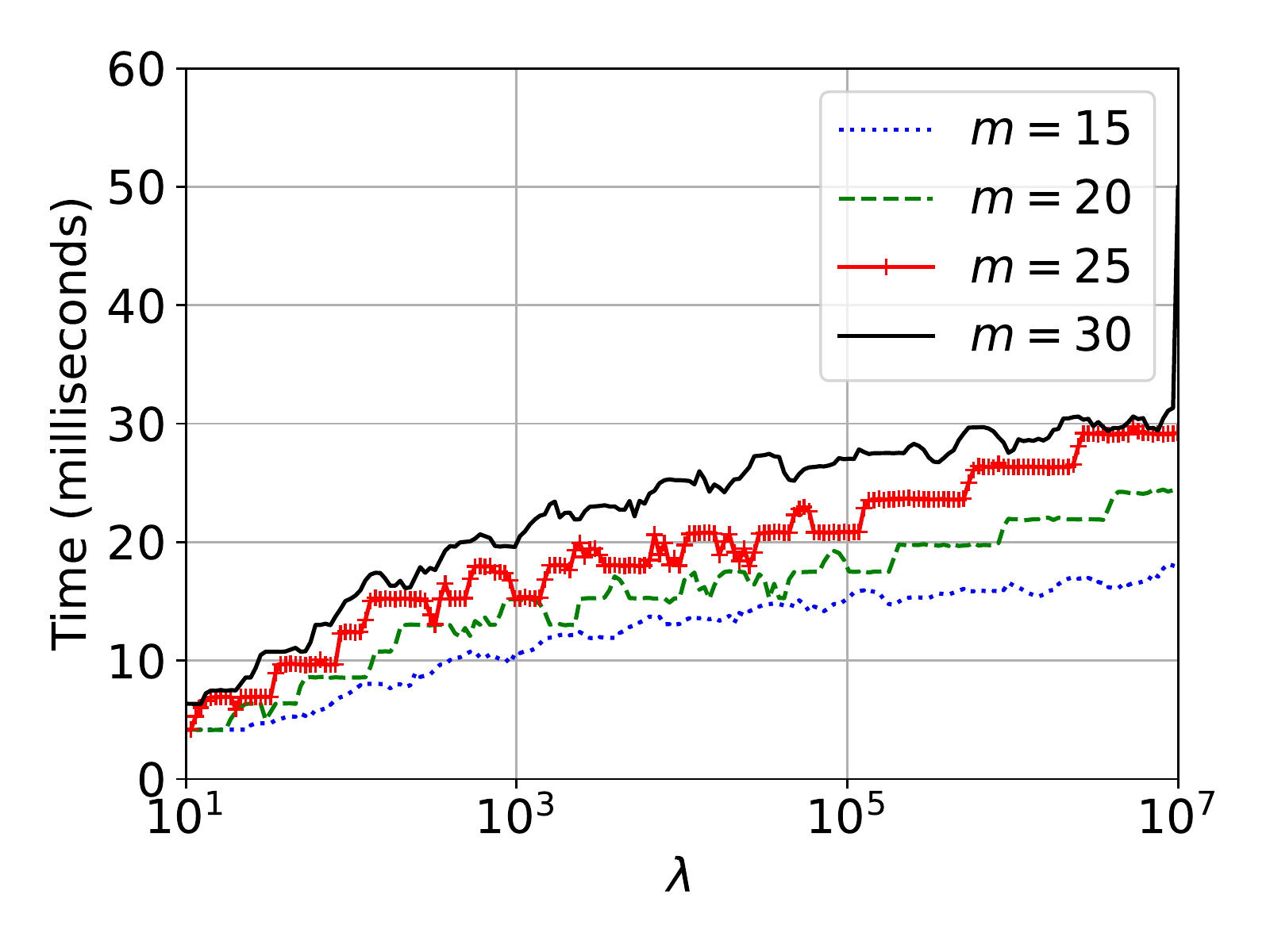}
\hfil
\includegraphics[width=.39\textwidth]{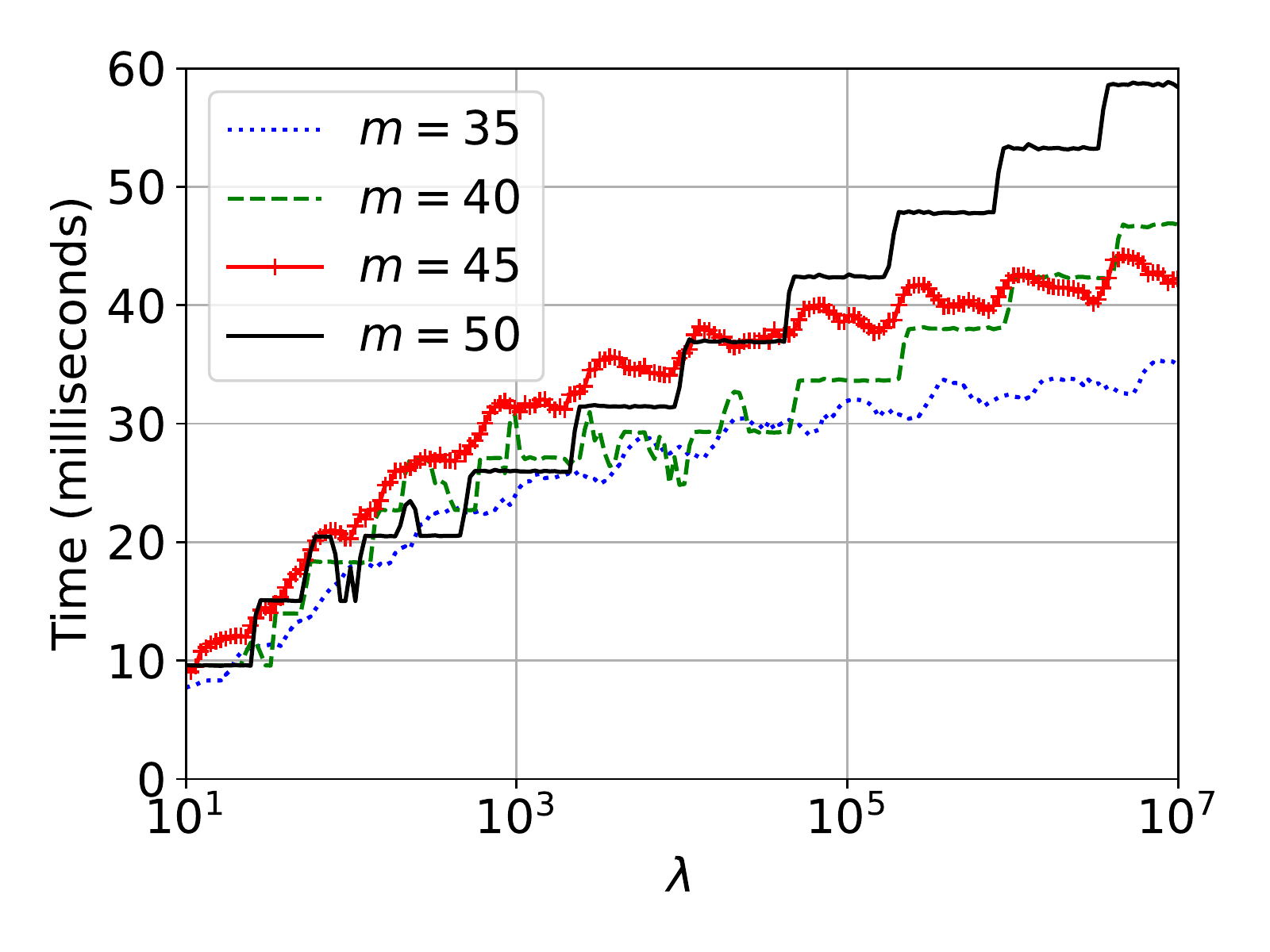}
\hfil

\hfil
\includegraphics[width=.39\textwidth]{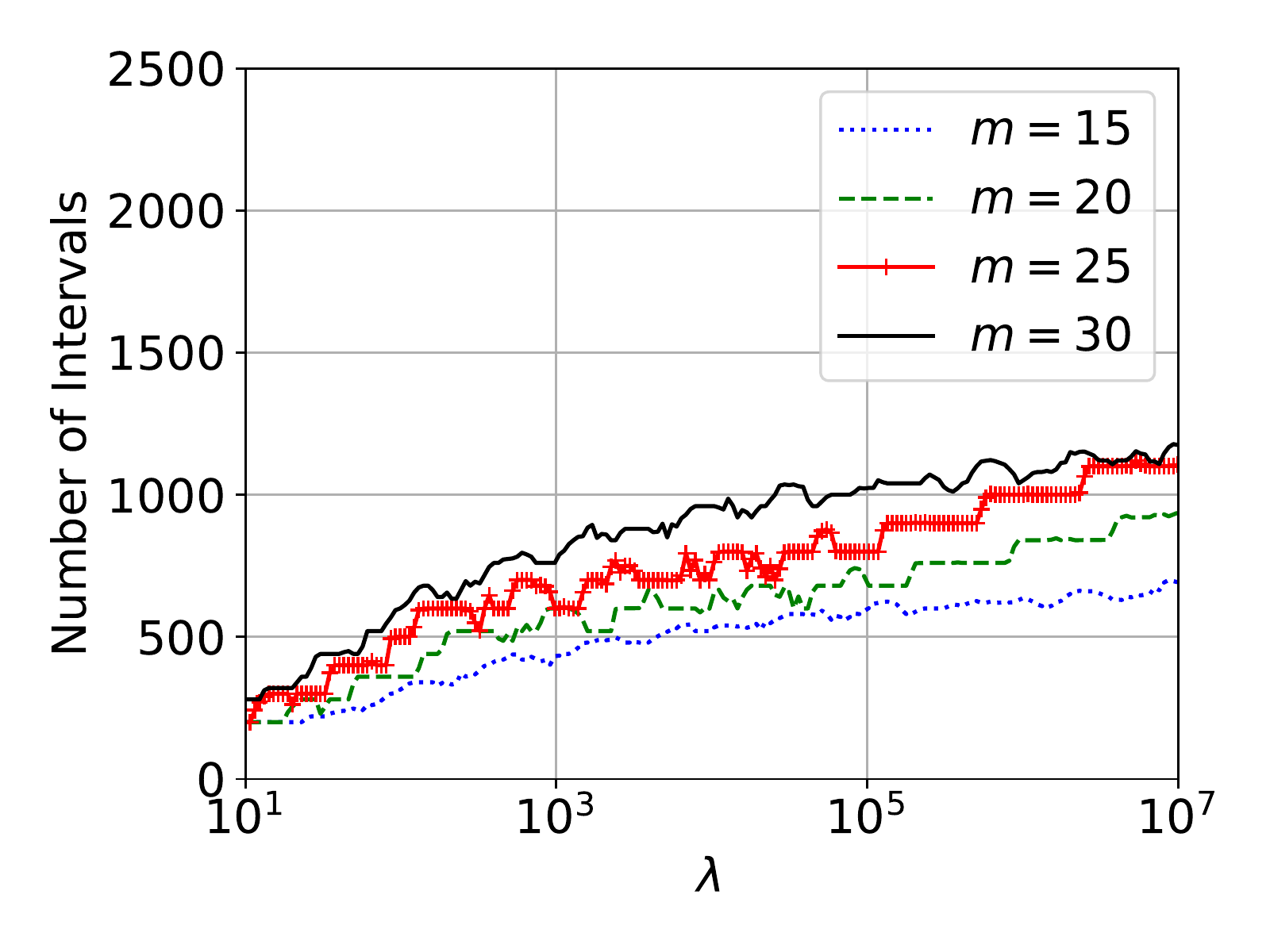}
\hfil
\includegraphics[width=.39\textwidth]{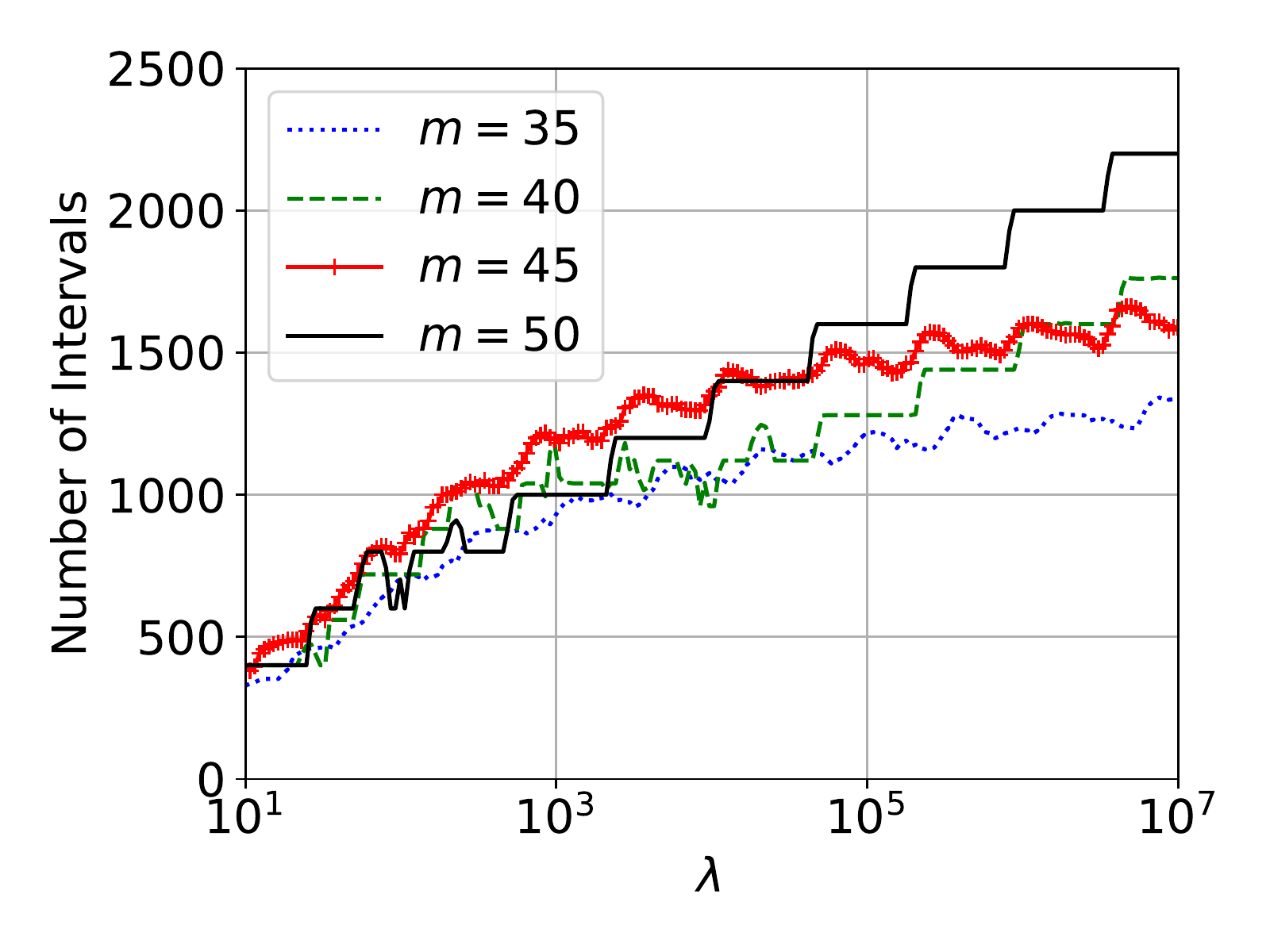}
\hfil

\hfil
\includegraphics[width=.39\textwidth]{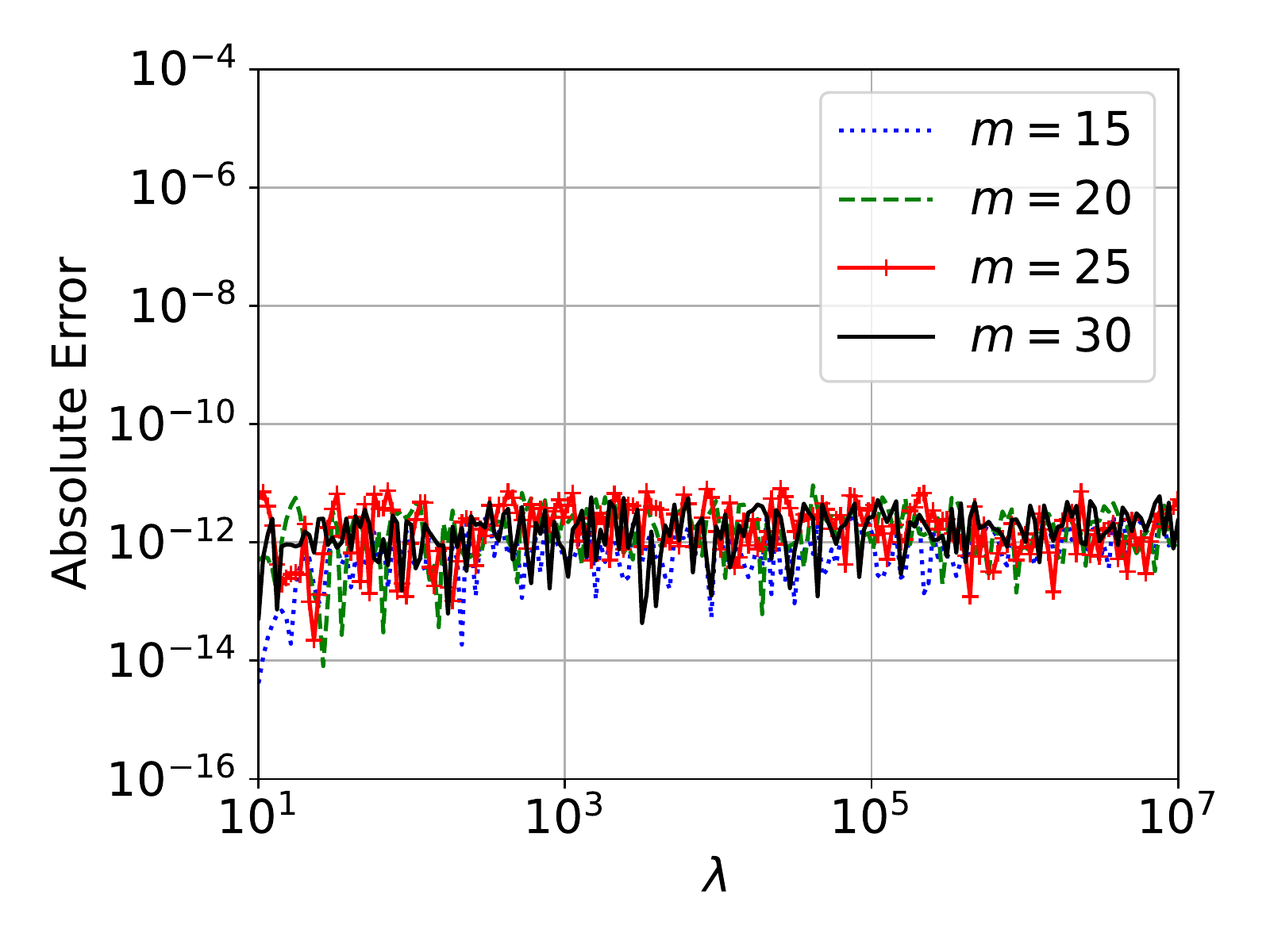}
\hfil
\includegraphics[width=.39\textwidth]{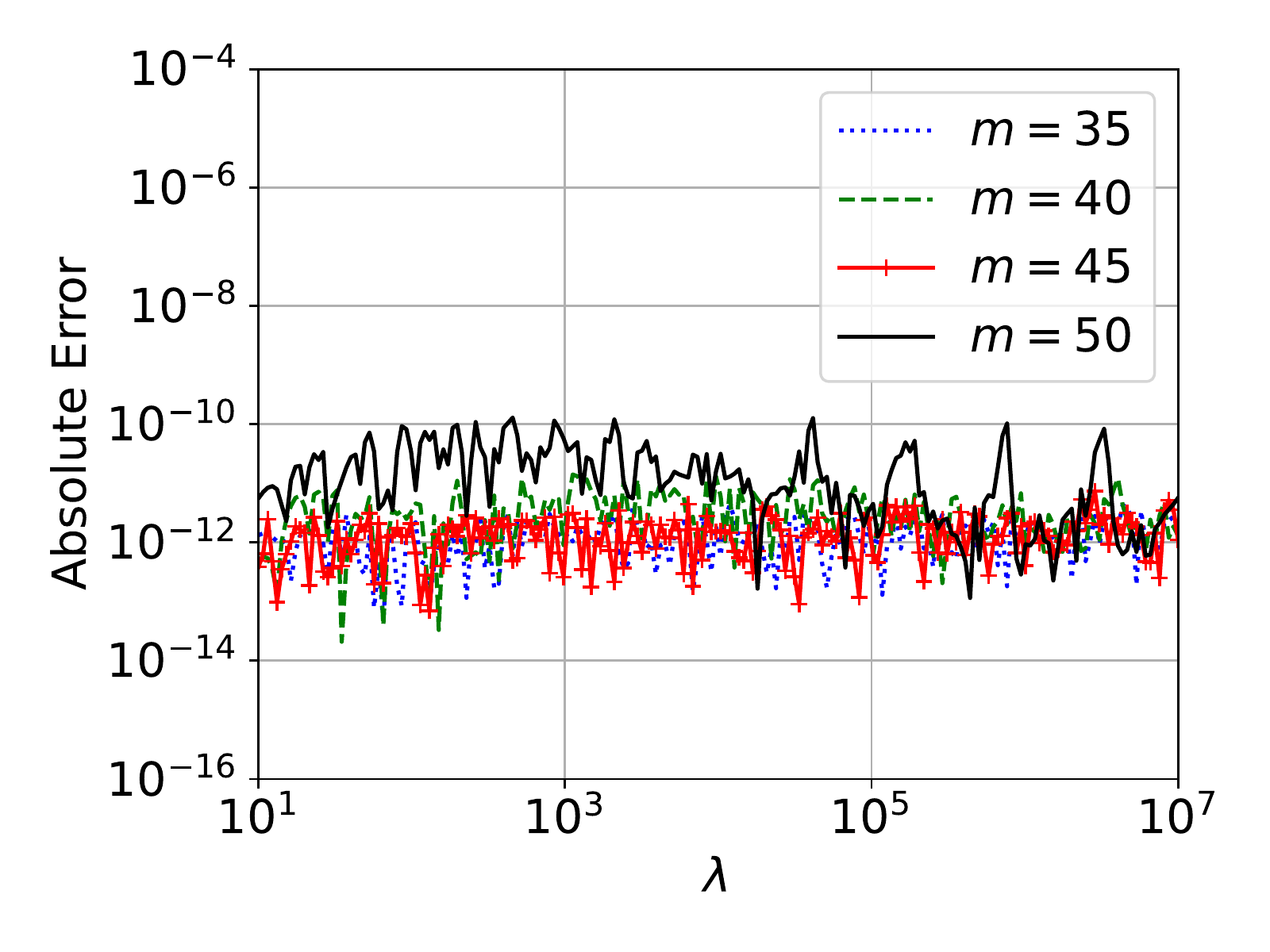}
\hfil

\caption{\small
The results of the experiment of Section~\ref{section:experiments:saddle2}.
The plots in the first row give
the   time taken by the adaptive Levin method as a function of $\lambda$ for the various values
of $m$ considered.    Those in the second row give the 
number of intervals  in the adaptively determined subdivision of $[-1,1]$ used to compute the integral
$I_{22}$ as a function of $\lambda$ for the values of $m$ considered.
The third row of plots give
the error in the calculated value of  $I_{22}(\lambda,m)$ 
 as a function of $\lambda$ for $m=15,20,25,30,35,40,45,50$.
}
\label{figure:saddleplots2}
\end{figure}

A state-of-the-art method for evaluating the modal Green's functions
is discussed in \cite{modal}.  It is more efficient than the 
adaptive Levin method in all cases, and it is much faster
when $\alpha$ is close to $1$.
However, the adaptive Levin method is surprisingly competitive given that 
it is a general-purpose approach to evaluating oscillatory integrals  
and the algorithm of \cite{modal} is highly-specialized.

\end{subsection}

\begin{subsection}{An integral with many stationary points}
\label{section:experiments:saddle2}

In this final  experiment, we considered the integral
\begin{equation}
I_{22}(\lambda, m) = \int_{-1}^1 \exp\left(i \lambda \cos^2\left( \frac{\pi}{2} m x\right)\right) \frac{1}{1+x^2}\,dx,
\end{equation}
which has $m$ stationary points.   We sampled $l=200$ equispaced points
$x_1,\ldots,x_l$ in the interval $[1,7]$ and, for each 
$\lambda=10^{x_1},10^{x_2},\ldots,10^{x_l}$
and
$m=15,20,25,30,35,40,45,50$,  we evaluated  $I_{22}(\lambda,m)$ 
using the adaptive Levin method.  

Figure~\ref{figure:saddleplots2} presents the results.
We see that, when applied to $I_{22}$,  the running time of the adaptive Levin method
grows sublogarithmically with $\lambda$ for all of the values of $m$ 
considered.  The cost also increases quite mildly as a function of the number of stationary points.
It is notable that the method is able to evaluate
$I_{22}\left(10^7,20\right)$ --- a  highly oscillatory integral with $50$ stationary points ---
to 11 digit accuracy in approximately 50 milliseconds.

\end{subsection}

\end{section}

\begin{section}{Conclusions}
\label{section:conclusions}

We have shown that the Levin method does not suffer from numerical breakdown
when the magnitude of $g'$ is small or when $g'$ has zeros,
which explains the effectiveness of adaptive Levin methods.
We have also presented numerical experiments indicating that 
the adaptive Levin method of this paper can accurately and rapidly
evaluate a large class of oscillatory integrals, including
many with singularities and stationary points.     We have further demonstrated that combining the adaptive Levin method
with the algorithm of \cite{BremerPhase,Scalar1} allows  for the efficient evaluation of many 
integrals involving the solutions of differential equations. This class includes integrals involving most of the classical special functions,
as well combinations of such functions.

As the experiments of Section~\ref{section:experiments:modal} indicate, 
specialized techniques designed for particular narrow classes of oscillatory
integrals are often faster than the adaptive Levin scheme of this paper.
However, the numerical experiments described in  Section~\ref{section:experiments} show that the adaptive Levin
method provides  an efficient general-purpose mechanism for evaluating
a huge class of  oscillatory integrals. 
It appears to have the roughly the same behavior when applied to
oscillatory integrals as adaptive Gaussian quadrature does when used to evaluate smoothly varying integrands.

We note that there are obvious implications of the  adaptive Levin method for the rapid application of special
function transforms and the solution of second order linear inhomogeneous differential equations
that should be thoroughly investigated.

\end{section}

\begin{section}{Acknowledgements}
KS was supported in part by the NSERC Discovery 
Grants RGPIN-2020-06022 and DGECR-2020-00356.  
JB was supported in part by NSERC Discovery grant  RGPIN-2021-02613.
\end{section}

\bibliographystyle{sn-mathphys}
\bibliography{levin.bib}

\end{document}